\documentclass[letterpaper,12pt,leqno,oneside]{article}
\usepackage{float}
\usepackage{caption}
\usepackage{enumitem}
\usepackage{color}
\usepackage{ifsym}
\usepackage{bm}
\usepackage{exscale}
\usepackage{amsmath}
\usepackage{amsfonts}
\usepackage{wasysym}
\usepackage{stmaryrd}
\usepackage{amscd}
\usepackage{graphicx}
\usepackage{pagecolor,lipsum}
\usepackage{xcolor}
\usepackage{marvosym}
\usepackage{textcomp}
\usepackage{amsxtra}
\usepackage{amssymb}
\usepackage{theorem}
\usepackage{enumitem}
\usepackage{mathtools}
\usepackage{esint}
\usepackage[hidelinks]{hyperref}
\usepackage[final]{epsfig}
\usepackage{eqnarray}
\usepackage{hyperref}
\usepackage[textsize=tiny]{todonotes}
\newtheorem{proposition}{Proposition}[subsection]
\newtheorem{definition}[proposition]{Definition}

\newtheorem{lemma}[proposition]{Lemma}

{\theorembodyfont{\rmfamily}\newtheorem{remark}[proposition]{Remark}}
\newtheorem{theorem}[proposition]{Theorem}
\newtheorem{corollary}[proposition]{Corollary}

{\theorembodyfont{\rmfamily}\newtheorem{example}[proposition]{Example}}

\newcommand{\mr}{\mathrm}

\newfont{\abc}{cmtt10 scaled 1200}

\newcommand\ve{\varepsilon}

\newcommand\ra{\rightarrow}
\newcommand\cs{\symbol{35}}
\newcommand\p{\partial}
\newcommand\qed{\hfill $\Box$ \\}

\newcommand\mm{\mbox}
\renewcommand\v{= \emptyset}
\newcommand\n{\neq \emptyset}
\newcommand\D{\ensuremath{\mathbf{ID}}}

\newcommand\bp{\langle A \rangle}

\newcommand\si{\ensuremath{\mathcal{S}}} 
\newcommand\sima\si

\newcommand\R{\mathbb{R}}

\newcommand\Z{\mathbb{Z}}

\newcommand\U{\mathbb{U}}
\renewcommand\P{\mathbb{P}}

\newcommand\TP{\mathbb{TP}}

\newcommand\Vol{\mr{Vol}}
\newcommand\diam{\mr{diam}}
\newcommand\dist{\mr{dist}}
\DeclareMathOperator{\tr}{tr}
\DeclareMathOperator{\scal}{scal}
\DeclareMathOperator{\Ric}{Ric}

\def\Xint#1{\mathchoice
{\XXint\displaystyle\textstyle{#1}}%
{\XXint\textstyle\scriptstyle{#1}}%
{\XXint\scriptstyle\scriptscriptstyle{#1}}%
{\XXint\scriptscriptstyle\scriptscriptstyle{#1}}%
\!\int}
\def\XXint#1#2#3{{\setbox0=\hbox{$#1{#2#3}{\int}$ }
\vcenter{\hbox{$#2#3$ }}\kern-.6\wd0}}

\def\fint{\Xint-}

\begin{document}

\vspace*{0.1cm}
\begin{center}\Large{\bf{Scalar Curvature Splittings I: Minimal Factors}}\\
\medskip
\large{\bf{Joachim Lohkamp}}\\
\end{center}
\noindent Mathematisches Institut, Universit\"at M\"unster, Einsteinstra\ss e 62, Germany\\
{\emph{e-mail: j.lohkamp@uni-muenster.de}}
{\footnotesize  {\center  \tableofcontents}}
\setcounter{section}{1}
\renewcommand{\thesubsection}{\thesection}

\subsection{Introduction} \label{introduction}

Lower scalar curvature bounds, $\scal \ge \kappa$ for some $\kappa \in \R$, on a manifold $M^{n+1}$ can be studied using splitting techniques where we inductively consider related scalar curvature constraints on suitable subspaces $V^k \subset M^{n+1}$, $k \le n$, also called the \textbf{splitting factors} of $M^{n+1}$. This strategy was introduced by Hawking \cite{H} and Schoen and Yau \cite{SY1} in the 70ties. Since then the scope of this approach has grown considerably, in particular, due to work by Gromov, Lawson, Galloway and Schoen in \cite{GL},\cite{GS} and \cite{G}.
Typical splitting factors are area minimizing (and some types of almost minimizing) hypersurfaces of $M^{n+1}$. It is a characteristic of this approach that the inherited scalar curvature constraints on  the splitting factors become accessible only after balancing global conformal deformations.
\subsubsection{Minimal Splitting Factors}\label{in}

This strategy depends on the control we have on the deformed hypersurfaces. In low dimensions $\le 7$ they are smooth manifolds and it is easy to inductively descend to lower dimensions. This changes in higher dimensions. An almost area minimizer $H^n$ can have a complex singular set $\Sigma_H \subset H^n$. When we approach $\Sigma_H$ the regular part $H^n \setminus \Sigma_H$ degenerates and the obligatory conformal deformations of $H^n \setminus \Sigma_H$ diverge. From Martin theory \cite{L2} on $H^n \setminus \Sigma_H$ we can tell that completing such conformally deformed  $H^n \setminus \Sigma_H$ can yield still more complicated singular spaces $X^n$ with a new singular set $\Sigma_X$.\\
This said, we show that there are well-controlled \textbf{minimal $\boldsymbol{\scal>0}$-model} geometries $X^n$ on $H^n$ supporting an accessible geometric analysis on metric measure spaces in the sense of Ambrosio, Cheeger and others, for instance, in \cite{A},\cite{C},\cite{CK},\cite{He},\cite{H-T},\cite{B-T}.

  For the sample case of a compact area minimizer in a $\scal>0$-manifold $M^{n+1}$ with induced Riemannian metric $g_H$ on $H^n\setminus \Sigma_H$ we get:\\

\textbf{Minimal Splitting Factors}\, \emph{There is a conformally deformed metric $\Phi^{4/n-2} \cdot g_H$, for some $C^{2,\alpha}$-regular $\Phi >0$, $\alpha \in (0,1)$, which we call the \textbf{minimal factor metric}, so that:
\begin{itemize}[leftmargin=*]
  \item  the metric completion $(X^n,d_X)$ of $(H^n\setminus \Sigma_H,\Phi^{4/n-2} \cdot g_H)$  is compact and \textbf{homeomorphic} to the completion $(H^n,d_H)$ of $(H^n,g_H)$.
\item $\scal(\Phi^{4/n-2} \cdot g_H)>0$ and in any $p \in \Sigma_X$,  $X^n$ has $\boldsymbol{\scal>0}$\textbf{-curved tangent cones}. This permits an inductive asymptotic analysis of $X$ near $\Sigma_X$ similar to the case of area minimizers. In particular, we get that  $\Sigma_X \subset X$ has \textbf{Hausdorff codimension} $\boldsymbol{\ge 7}$.
\item  $X^n$ can be augmented to a \textbf{metric measure space}  that is \textbf{Ahlfors $n$-regular}, in particular,  \textbf{doubling},  and that admits \textbf{Poincar\'{e} inequalities}. This implies further regularity properties of $X$ like the presence of \textbf{isoperimetric inequalities}.
\end{itemize}}

The defining property of $X^n$ is that $\Phi$ has minimal growth towards $\Sigma$, compared to other admissible deformations. Details and extensions are explained in Ch.\ref{over} below.\\ In part II \cite{L5} we use the symmetries of $\scal>0$-curved tangent cones and the control from the Ahlfors $n$-regularity and the isoperimetric inequality to inductively remove $\Sigma$.  A typical  application, we cite from  \cite{L5}, is a splitting scheme with built-in stepwise regularization:\\

\textbf{Partial Regularization}\, {\itshape Let $H^{n}$, $n \ge 2$, be a compact area minimizer, with singular set $\Sigma$, in a $\scal >0$-manifold $M^{n+1}$.
Then there are arbitrarily small neighborhoods $U$ of $\Sigma$, so that $H \setminus U$ is \textbf{conformal} to a $\boldsymbol{\scal>0}$\textbf{-manifold} $(X_U,g_X)$ with \textbf{minimal} boundary $\p X_U$.}\\

 The point is that although, in general, the boundary $\p X_U$ is also singular, its singular set $\Sigma_{\p X_U}$ will have a \emph{lower dimension} than $\Sigma_H$. The minimality of $\p X_U$, in the $\scal>0$  ambient manifold $X_U$\footnote{Actually $\p X_U$ is two-sided minimal in a slightly larger (non-complete) $\scal>0$-manifold $Y_U \supset \overline{X_U}$.}, then allows us to iteratively shift singular problems to lower dimensions before they disappear in dimension $7$, cf.\ the introduction of \cite{L5} and the survey \cite{L4}. \\
In turn, Schoen and Yau \cite{SY2} have described an alternative strategy using nestings of singular minimizers. Minimal splitting factors may be of use in such a setting as well.\\

\textbf{Comparison with the Classical Approach}\,  To get minimal splitting factor we use a setup that differs from the traditional approach, as used for instance in \cite{SY1} or \cite{GL},  for a \emph{regular} compact area minimizing hypersurface $H^n$  in a  $\scal>0$-manifold $M^{n+1}$. In that classical case one considers the first eigenfunction $f_H>0$ of the conformal Laplacian  $L_H = -\Delta  +\frac{n-2}{4 (n-1)} \cdot \scal_H$, i.e. $L_H f_H=\lambda_H \cdot f_H$. Then the stability of $H$ implies that the first eigenvalue $\lambda_H$ is positive and  the  transformation law
\begin{equation} \label{trl} \scal(f_H^{4/(n-2)} \cdot g_H) \cdot f_H^{\frac{n+2}{n-2}} = L_H(f_H) = \lambda_H \cdot f_H > 0,\end{equation}
shows that $\scal(f_H^{4/(n-2)} \cdot g_H) >0$, since $\lambda_H >0$ and $f_H>0$. Turning to the singular case $\Sigma_H \n$ we make some simple but essential modifications. To explain them let $A_H$ be the second fundamental form on $H \setminus \Sigma  \subset M$, $|A|$ is its norm and we use a fixed \si-transform $\bp>0$, cf.\ Ch.\,\ref{bcn}. For the present we may think of $\bp$ as a revamped version of  $|A|$.
\begin{itemize}[leftmargin=*]
  \item In place of the ordinary eigenvalue equation $L_H f_H=\lambda_H \cdot f_H$ we consider an \emph{$\bp$-weighted} eigenvalue equation, i.e. the eigenvalue equation for ${\bp}^{-2} \cdot L$:
\begin{equation}\label{weee}
L_H (u_\lambda) = \lambda \cdot \bp^2 \cdot u_\lambda \,\mm{ on } H \setminus \Sigma, \mm{ for some } \lambda >0.
\end{equation}\noindent
We note that (\ref{weee}), and in particular $\lambda$, remains \emph{invariant under global scalings} of $H$.
  \item Different from the regular case $\Sigma_H \v$ we have many different positive solutions also for any subcritical eigenvalue $\lambda  < \lambda^{\bp}_{H}$, where $\lambda^{\bp}_{H}$ is the principal eigenvalue of ${\bp}^{-2} \cdot L$. One can show that $\lambda^{\bp}_{H}>0$. We \emph{choose} some \emph{subcritical} eigenvalue $\lambda \in (0,\lambda^{\bp}_{H})$. The point is that for these $\lambda$ the potential theory of $L_H - \lambda \cdot \bp^2$ is  particularly well-controlled near $\Sigma$ even without knowing any structural detail of $\Sigma$.
  \item We \emph{choose} a (super)solution $\Phi>0$ of (\ref{weee}) with \emph{minimal growth} towards $\Sigma_H$ to conformally deform $g_H$ to  the \emph{minimal factor metric} $\Phi^{4/(n-2)} \cdot g_H$. This keeps the Hausdorff dimension of the new singular set small  and, combined with the scaling invariance of  (\ref{weee}), this choice yields $\scal>0$-curved tangent \emph{cones} as the exclusive blow-up geometries in singular points.
\end{itemize}
 In the regularizations of \cite{L5} we  employ the asymptotic geometry  of $(H^n\setminus \Sigma_H,\Phi^{4/n-2} \cdot g_H)$ near $\Sigma$, resembling that of  $\scal>0$-cones, to construct surgery style deformations bending $(H \setminus U,\Phi^{4/(n-2)} \cdot g_H)$ to a $\scal >0$-manifold with minimal $\p U$ as described above. In turn, the isoperimetry of the minimal factor geometry is used to validate these properties.
\subsubsection{Statement of Results}\label{over}

The results hold for area minimizing hypersurfaces and for broader classes of almost minimizers $H \in {\cal{G}}$, provided we have $\lambda^{\bp}_H >0$. We recall the definition of $\cal{G}$ and other basics in Ch.\,\ref{bcn} below. For the sake of consistency we state our results for any $H \in {\cal{G}}$. For regular $H$ the results remain valid but they oftentimes become trivial, e.g., when subcritical eigenfunctions do not exist. With this caveat in mind, we consider, for the rest of the paper, a fixed pair of subcritical and principal eigenvalue (see Ch.\,\ref{bcn}.D.3):
\begin{equation}\label{lam}
\boldsymbol{0 < \lambda < \lambda^{\bp}_H}, \mm{ for the given } H \in {\cal{G}}.
\end{equation}
We need $\lambda>0$ to get conformal deformations to $\scal>0$-metrics whereas $\lambda^{\bp}_H - \lambda>0$ is crucial for the validity of the potential theoretic arguments from  \cite{L1}--\cite{L3}. The actual values are immaterial for the qualitative aspects of the theory.\\
Now we turn to the definition of our basic metrics.  Due to the locally Lipschitz regular coefficients of $L_{H,\lambda} := L_H - \lambda \cdot \bp^2$, solutions  of $L_{H,\lambda} \, \phi=0$ are $C^{2,\alpha}$-regular, for any $\alpha \in (0,1)$. This suggests the following regularity assumptions.
\begin{definition} \emph{\textbf{(Minimal Factor Metrics)}}  \label{msge}  For $H \in {\cal{G}}$ let  $\Phi>0$ be a $C^{2,\alpha}$-supersolution  of $L_{H,\lambda} \phi =0$ on $H\setminus\Sigma_{H}$ so that in the case
\begin{itemize}
  \item  $H \in {\cal{G}}^c$: $\Phi$  is a solution in a neighborhood of $\Sigma$ with \textbf{minimal growth} towards $\Sigma$.
  \item $H \in {\cal{H}}^{\R}_n$:  $\Phi$ is a solution on $H\setminus \Sigma_{H}$ with \textbf{minimal growth} towards $\Sigma$.
\end{itemize}
We call the $\scal>0$-metrics  $\Phi^{4/(n-2)} \cdot g_H$ on $H \setminus \Sigma$ the \textbf{minimal factor metrics}.
\end{definition}

\begin{remark} Minimal factor metrics are naturally assigned to any  $H \in {\cal{G}}$. This is owing to the boundary Harnack inequality \ref{mbhsq} (\cite[Theorem 3.4 and 3.5]{L2}). It shows that for any two such supersolutions $\Phi_1$, $\Phi_2$ on $H\setminus \Sigma_{H}$ we have some constant $c \ge 1$ so that  $c^{-1}  \cdot \Phi_1 \le \Phi_2 \le c \cdot \Phi_1$ near $\Sigma$. For $H \in {\cal{H}}^{\R}_n$, the boundary Harnack inequality even shows that $\Phi$ is uniquely determined up to multiples, i.e. $\Phi_2 \equiv  c \cdot \Phi_1$. Then constants in estimates depending on $(H,\Phi)$ will only depend on $H$ since $c$ equally appears on both sides of the respective inequality.\\
Since $\lambda < \lambda^{\bp}_H$, we do \emph{not} have regular positive solutions with minimal growth towards all of $\widehat{\Sigma}$ but we have a \textbf{minimal Green's function} $G(x,y)$ for $L_{H,\lambda}$. This is a function $G:H\setminus \Sigma \times H\setminus \Sigma \to(0,\infty]$ that is finite and $C^{2,\alpha}$-regular outside the diagonal $\{(x,x)\:|\:x \in H\setminus \Sigma \}$ and satisfies the equation $L_{H,\lambda}\,G(\cdot,y)=\delta_y$ in a distributional sense, where $\delta_y$ is the Dirac delta function with basepoint $y$ and  $G(\cdot ,y)$ has minimal growth towards $\widehat{\Sigma}$.\\
 This minimal Green's function is uniquely determined. $G(\cdot,y)$ is a supersolution, for any $y \in H\setminus \Sigma_{H}$. Throughout this paper we \textbf{exclusively} use \textbf{minimal} Green's functions.
\end{remark}

\begin{example}\label{exa} With $G(\cdot ,y)$  one may construct supersolutions which are proper solutions of minimal growth near $\Sigma$. For  any open set $V$ with compact closure $\overline{V} \subset H \setminus \Sigma$ and a smooth function $f$ on $H \setminus \Sigma$ with $f \equiv 0$ on  $(H \setminus \Sigma) \setminus V$ and $f> 0$ on $V$, we set
\begin{equation}\label{suppa}
\mathbf{S}(x)=\mathbf{S}[H,\lambda, V, f](x):=\int_{H \setminus \Sigma} G(x,y) \, f(y) \, dV(y).
\end{equation}
This is a smooth positive supersolution of $L_{H,\lambda} \phi=0$ on $H \setminus \Sigma$ with $\mathbf{S} \in H^{1,2}_{\bp}(H \setminus \Sigma)$, \cite[Lemma 3.11 and Prop.\,3.12]{L2}, and it solves $L_{H,\lambda} \phi=0$ away from $V$ with minimal growth towards $\widehat{\Sigma}$. The Riesz decomposition theorem shows that any regular supersolution that is a proper solutions outside $V$ and has minimal growth near $\Sigma$ can be written in the form (\ref{suppa}).
\end{example}

Minimal factor geometries share many fundamental geometric properties with the original (almost) minimizing geometry on $H$. This is closely tied to the minimal growth condition for $\Phi$. For general conformal deformations the results of this paper become invalid.

\begin{theorem}\emph{\textbf{(Singular Sets)}}\label{idi} For any $H \in {\cal{G}}$ we have:
\begin{itemize}
\item The metric completion $(X^n,d_X)$ of $(H \setminus \Sigma,\Phi^{4/(n-2)} \cdot g_H)$ is a geodesic metric space and it is \textbf{homeomorphic} to $(H,d_H)$, $(X^n,d_X) \cong (H,d_H)$ and, hence,  $\Sigma_X \cong \Sigma_H$. Thus, we can write it as $(H,d_{\sima}(\Phi))$ or briefly $(H,d_{\sima})$ and, for $H \in {\cal{G}}^c_n$,  $(H,d_{\sima})$ is compact.
\item The \textbf{Hausdorff dimension} of the singular set $\Sigma_X$ of $(X^n,d_X)=(H,d_{\sima})$ is $\le n-7$.
\end{itemize}
We call  $(H,d_{\sima})$ a minimal spitting factor, briefly a \textbf{minimal factor}, of its ambient space $M$ and $d_{\sima}$ the completed \textbf{minimal factor metric} extending the definition \ref{msge} on $H \setminus \Sigma$.
\end{theorem}
The second assertion is a non-trivial refinement of $\Sigma_X \cong \Sigma_H$ since the identity map $id_H: (H,d_H) \ra (H,d_{\sima})$ is \textbf{not} Lipschitz regular. What we show is that the upper dimensional bound for $\Sigma_X$ is again $n-7$ but we do not know whether the Hausdorff dimension, in particular of lower dimensional pieces of $\Sigma_H$, remains unchanged cf. Remark \ref{sma} below.\\
The following two results show that minimal factors form a \textbf{blow-up invariant} class of spaces and, in singularities, any of these spaces admits $\boldsymbol{\scal>0}$\textbf{-tangent cones}. This will be used  to inductively study the scalar curvature geometry near $\Sigma$, in particular in \cite{L5}.

\begin{theorem}\emph{\textbf{(Blow-Ups)}}\label{blooin}
For $H \in {\cal{G}}$ we consider $(H, d_{\sima}(\Phi_H))$, any $p \in \Sigma_H$ and any tangent cone $C$ in $p$. Then we get the following \textbf{blow-up invariance:}\\
Any sequence  $(H, \tau_i \cdot d_{\sima}(\Phi_H))$ scaled by a sequence $\tau_i \ra \infty$, $i \ra \infty$, around $p$, subconverges\footnote{This is a convergence of the underlying minimizers and the conformal deformation, cf.~Ch.~\ref{bcn}.E.1.} and the limit of any converging subsequence is $(C, d_{\sima}(\Phi_C))$ for some tangent cone $C$.
\end{theorem}
 \begin{theorem}\emph{\textbf{(Euclidean Factors)}} \label{euclfin}  For any non--totally geodesic $H \in {\cal{H}}^{\R}_n$ there is a \textbf{unique}$^\cs$  space  $(H, d_{\sima}(\Phi_H))$, i.e. unique up to global scaling. For $C \in \mathcal{SC}_{n}$ the associated space  $(C, d_{\sima}(\Phi_C))$ is invariant under scaling around $0 \in C$, that is, it is again a cone.
 \end{theorem}
There is no such blow-up invariant scheme for \emph{principal} eigenvalues. The $\bp$-weighted principal eigenvalues are scaling invariant and  we have $\lambda^{\bp}_H \le \lambda^{\bp}_C$, and typically $\lambda^{\bp}_H <\lambda^{\bp}_C$, cf.~\cite[Lemma 3.9]{L3}. This means that, in general, $\lambda^{\bp}_H$ is a \emph{non-principal} eigenvalue on $C$. In \cite{L4}, we actually use this as a degree of freedom. We \emph{choose} a $\lambda$ that is much smaller than $\lambda^{\bp}_H$  to get lower bounds on the growth rate of solutions towards the singular set.\\

$(H, d_{\sima}(\Phi_H))$ admits a canonical augmentation to a \textbf{metric measure space} from an extension of $\Phi^{2 \cdot n/(n-2)}\cdot \mu_H$ on $H  \setminus \Sigma$,  $\mu_H$ is the $n$-dimensional Hausdorff measure on $(H^n,d_H)$.

\begin{definition}\emph{\textbf{(Minimal Factor Measures)}}\label{mms} For any $H \in {\cal{G}}_n$ equipped with a minimal factor metric
$\Phi^{4/(n-2)} \cdot g_H$, we define the \textbf{minimal factor measure} $\mu_{\sima}$ on $H$ by
\begin{equation}\label{meas}
\mu_{\sima}(E):=\int_{E \setminus \Sigma_H} \Phi^{2 \cdot n/(n-2)}\cdot d\mu_H, \mm{ for any Borel set } E \subset H.
\end{equation}
\end{definition}
$\mu_{\sima}$ is a \textbf{Borel measure} on $(H,d_{\sima})$, cf.~\cite[pp.~62--64]{H-T}. This  uses the Hausdorff dimension estimate for $\Sigma \subset (H, d_{\sima}(\Phi_H))$, Theorem \ref{idi}, and  the Ahlfors regularity, Theorem  \ref{dvintro}, below.

\begin{theorem}\emph{\textbf{(Ahlfors Regularity)}} \label{dvintro}  For $H \in {\cal{G}}_n$, the space $(H,d_{\sima},\mu_{\sima})$ is \textbf{Ahlfors $\boldsymbol{n}$-regular}: there are constants $A(H,\Phi),B(H,\Phi)>0$ so that for any $q\in H$:
\begin{equation}\label{ahl}
 A \cdot r^n \le \mu_{\sima}(B_r(q),d_{\sima}) \le B \cdot r^n, \mm{ for any } r \in [0,\diam(H,d_{\sima})).
\end{equation}
For $H \in {\cal{H}}^{\R}_n$ the constants only depend on the dimension, that is, we have $A(n),B(n)>0$.
\end{theorem}
\begin{corollary}\emph{\textbf{(Doubling Properties)}} \label{dvint} For any $H \in {\cal{G}}_n$ there is a  $C(H,\Phi)>0$,  and $C(n)>0$ for $H \in {\cal{H}}^{\R}_n$, so that  $(H,d_{\sima},\mu_{\sima})$ has the following properties:
\begin{itemize}
  \item $\mu_{\sima}$ is \textbf{doubling}: for any $q \in H$ and $r \in [0,\diam(H,d_{\sima}))$:
  \begin{equation}\label{dou}
  \mu_{\sima}(B_{2 \cdot r}(q),d_{\sima}) \le C \cdot \mu_{\sima}(B_{r}(q),d_{\sima}).
  \end{equation}
  \item  For balls $B_* \subset B \subset H$ we have a \textbf{relative lower volume decay} of order $n$:
  \begin{equation}\label{volgro}
  \diam(B_*)^n/\diam(B)^n \le C \cdot \mu_{\sima}(B_*)/\mu_{\sima}(B).
  \end{equation}
  \item For  $H \in {\cal{G}}^c_n$, the total volume relative to $\mu_{\sima}$ is finite: $\mu_{\sima}(H) < \infty$.
\end{itemize}
\end{corollary}

$(H,d_{\sima},\mu_{\sima})$ admits \textbf{Semmes families of curves}. The presence of such families is essential to establish Poincar\'{e} inequalities, cf.~\cite{Se, He, H-T}, including the following version.

\begin{theorem} \emph{\textbf{(Poincar\'{e} inequality)}}\label{sobb}
For any $H \in {\cal{G}}$, there are $C_0(H,\Phi) >0$,  $\gamma_0(H,\Phi) \ge 1$, depending only on $n$ for $H \in {\cal{H}}^{\R}_n$, so that for concentric balls $B\subset \gamma_0 \cdot B\subset (H,d_\sima)$, for any function $u$ on $H$, integrable on bounded balls, and every upper gradient $w$ of $u$  we have:
\begin{equation}\label{poinm0}
\fint_B |u-u_B| \,  d \mu_{\sima} \le C_0 \cdot \diam(B) \cdot \fint_{\gamma_0 \cdot B} w \, d \mu_{\sima}, \mm{ where } f_B:=\fint_B f  \,  d \mu_{\sima} := \int_B f \, d \mu_{\sima}/\mu_{\sima}(B).
\end{equation}
\end{theorem}

A measurable function $w \ge 0$ on  $(H , d_{\sima})$ is an \emph{upper gradient} of $u$ if  $|u(x)- u(y)| \le  \int_c w(s) ds$ for all rectifiable curves $c$ joining $x$ to $y$, for any pair $x,y \in H$. An example, for $u \in Lip_{loc}(\Omega)$, is $|\nabla u|(x) :=\liminf_{\varrho \ra 0} \sup_{y \in\overline{ B _{\varrho}(x)}}  |u(x)-u(y)|/\varrho$, cf.\cite[p.982]{M}.  Using Cor.\ref{dvint} we can improve Theorem \ref{sobb} using \cite[Theorem 9.1.15]{H-T},  \cite[Theorem 4.5, Remark 4.6]{M}.

\begin{corollary} \emph{\textbf{(Sobolev Inequality)}}\label{sobb0}
For any $H \in {\cal{G}}$, there is a constant $C_1(H,\Phi)>0$, depending only on $n$ for $H \in {\cal{H}}^{\R}_n$,  so that for some open ball $B \subset H$, an $L^1$-function $u$ on $B$ and every upper gradient $w$ of $u$ on $B$, we have
\begin{equation}\label{ii20}
\Big(\fint_B |u-u_B|^{n/(n-1)} \,  d \mu_{\sima}\Big)^{(n-1)/n}  \le C_1 \cdot \diam(B) \cdot \fint_B  w \, d \mu_{\sima}.
\end{equation}
\end{corollary}

In turn, this and  \ref{dvintro} improve \ref{sobb} to the case where $\gamma_0=1$, cf.\cite[Remark 9.1.19]{H-T}. Theorems \ref{dvintro} and \ref{sobb} and work of Ambrosio and Miranda, in \cite{A} and \cite{M}, show that $(H,d_{\sima},\mu_{\sima})$ admits a proper BV(=bounded variations) theory. For a perimeter concept $\mu^{n-1}_{\sima}$ adapted to $(H,d_{\sima},\mu_{\sima})$ we have isoperimetric inequalities for minimal splitting factors.

\begin{theorem} \emph{\textbf{(Isoperimetric Inequality)}}\label{iip}
For $H \in {\cal{G}}$ there is a constant $\gamma(H,\Phi)>0,$ depending only on $n$ when $H \in {\cal{H}}^{\R}_n$, so that for any Caccioppoli set, i.e. a Borel set with locally finite perimeter,  $U \subset H$:
\begin{equation}\label{iinn}
\min  \{ \mu_{\sima}(B_{\rho} \cap U),  \mu_{\sima} (B_{\rho} \setminus U)\}^{(n-1)/n} \le \gamma \cdot \mu^{n-1}_{\sima}(B_{\rho} \cap \p U), \mm{ for any } \rho>0.
\end{equation}
\end{theorem}
As a consequence of the Ahlfors regularity and the isoperimetric inequality we have
\begin{corollary} \emph{\textbf{(Volume Growth)}}\label{agr} For $(H,d_{\sima},\mu_{\sima})$, some open  subset $\Omega \subset H$ and an oriented minimal boundary $L^{n-1} \subset \Omega$ bounding an open set $L^+ \subset \Omega$ there are constants $\kappa, \kappa^+(H,\Phi)  >0$, so that for any $p \in L$:
\begin{equation}\label{est}
\kappa \cdot r^{n-1}  \le \mu^{n-1}_{\sima}(L \cap B_r(p)) \,\mm{ and }\, \kappa^+\cdot r^n  \le  \mu_{\sima}(L^+ \cap B_r(p)),
\end{equation}
for $r \in [0, (A/B)^{1/n} \cdot \dist(p,\p \Omega)/4)$ and where $0<A<B$ are the Ahlfors constants. For $H \in {\cal{H}}^{\R}_n$, $\kappa, \kappa^+  >0$ depend only on $n$.
\end{corollary}

\begin{remark}  The methods and results of this paper carry over to Plateau problems, that is, to \textbf{(almost) minimizers with boundary}. They equally admit hyperbolic unfoldings and the associated potential theory. The needed regularity assumptions for the hypersurfaces and adaptedness properties of the operators are specified in  \cite[Remark 3.10]{L2}.
\end{remark}

\textbf{Organization of the Paper} \, The main results are the Ahlfors regularity and the Poincar\'{e} inequality for $(H,d_{\sima},\mu_{\sima})$. They broadly use the conformal Gromov hyperbolic structure we have on area minimizing (and on almost minimizing) hypersurfaces, cf.~\cite{L1}--\cite{L3}. In Ch.~2 we define canonical Semmes families of curves on the original (almost) area minimizer $H$ using hyperbolic geodesics on $H$. Then we use Doob transforms to derive controls for the length of curves in such families under the transition to $(H,d_{\sima})$. A first application is that the metric completion of $(H \setminus \Sigma,\Phi^{4/(n-2)} \cdot g_H)$ is homeomorphic to $(H,d_H)$. From a geometric reinterpretation of work from \cite{L2} and \cite{L3} we get the existence of $\scal>0$--tangent cones. We combine this with the length control to estimate the Hausdorff dimension of $\Sigma$ relative to $(H,d_{\sima})$. In Ch.~3 we estimate the eccentricity of balls $(H,d_{\sima})$ relative to balls in $(H,d_H)$. This becomes an ingredient for estimates of the volume of distance balls in the metric measure space $(H,d_{\sima},\mu_{\sima})$ to derive the Ahlfors regularity of $(H,d_{\sima},\mu_{\sima})$. Finally we verify that the canonical Semmes families satisfy the Semmes axioms also relative to $(H,d_{\sima},\mu_{\sima})$.  These results imply the validity of a Poincar\'{e}  inequality on $(H,d_{\sima},\mu_{\sima})$.

\begin{remark}   The present paper, together with its second part \cite{L5},  extend our earlier (unpublished) lecture notes \cite{L6} to improve their accessibility and to broaden the range of applications. In \cite{L4} we survey the approach and the vital  r\^{o}le of hyperbolic geometry.
\end{remark}

\subsubsection{Basic Concepts}\label{bcn}%
We summarize some basic notations, concepts and results from \cite{L1}--\cite{L3} we use in this paper.\\

\textbf{A. Basic Classes} of \emph{integer multiplicity rectifiable currents} of dimension $n\ge 2$  with connected
support inside some complete, smooth Riemannian manifold $(M^{n+1},g_M)$
{\small \begin{description}
  \item[${\cal{H}}^c_n$:] $H^n \subset M^{n+1}$ is  compact locally mass minimizing without boundary.
  \item[${\cal{H}}^{\R}_n$:] $H^n \subset\R^{n+1}$ is a complete hypersurface in flat Euclidean space $(\R^{n+1},g_{eucl})$ with $0\in H$ that is an oriented minimal boundary of some open set in $\R^{n+1}$.
  \item[${\cal{H}}_n$:] ${\cal{H}}_n:= {\cal{H}}^c_n \cup {\cal{H}}^{\R}_n$ and ${\cal{H}} :=\bigcup_{n \ge 1} {\cal{H}}_n$. We briefly refer to $H \in {\cal{H}}$ as an \textbf{area minimizer}.
  \item[$\mathcal{C}_{n}$:] $\mathcal{C}_{n} \subset {\cal{H}}^{\R}_n$ is the space of area minimizing $n$-cones in $\R^{n+1}$ with tip in $0$.
  \item[$\mathcal{SC}_{n}$:] $\mathcal{SC}_n \subset\mathcal{C}_n$ is the subset of cones which are at least singular in $0$.
  \item[${\cal{G}}^c_n$:] $H^n \subset M^{n+1}$ is a compact \textbf{almost minimizer}, cf.~\ref{bcn}.D. below.  We set ${\cal{G}}^c :=\bigcup_{n \ge 1} {\cal{G}}^c_n.$
  \item [${\cal{G}}_n$:] ${\cal{G}}_n := {\cal{G}}^c_n \cup {\cal{H}}^{\R}_n$ and ${\cal{G}} :=\bigcup_{n \ge 1} {\cal{G}}_n$.
    \item[$\mathcal{K}_{n-1}$:] For any area minimizing cone $C \subset \R^{n+1}$ with tip $0$, we get the (non-minimizing) minimal hypersurface $S_C:= \p B_1(0) \cap C \subset S^n \subset  \R^{n+1}$
  and we set ${\cal{K}}_{n-1}:= \{ S_C\,| \, C \in {\mathcal{C}_{n}}\}$. We write ${\cal{K}}= \bigcup_{n \ge 1} {\cal{K}}_{n-1}$ for the space of all such hypersurfaces $S_C$.
\end{description}}
We note that each of the classes ${\cal{H}}_n$, ${\cal{H}}^{\R}_n$ and ${\cal{G}}_n$ is \emph{closed under blow-ups}.\\

\textbf{B. One-Point Compactifications} of hypersurfaces $H \in {\cal{H}}^{\R}_n$ are denoted by $\widehat{H}$. For the singular set $\Sigma_H$ of some $H \in {\cal{H}}^{\R}_n$ we \emph{always} add $\infty_H$   to $\Sigma$ as well, even when $\Sigma$ is already compact, to define $\widehat\Sigma_H:=\Sigma_H \cup \infty_H$. On the other hand, for $H \in{\cal{G}}^{c}_n$ we set $\widehat H=H$ and $\widehat\Sigma=\Sigma$.\\

\textbf{C.1. \si-Structures}  An \si-transform $\bp$ is a distance measure to singular and highly curved parts of an almost minimizer. There are several ways to define such an $\bp$, but they all share some simple properties: an assignment $\bp$ which associates with any $H \in {\cal{G}}$ a locally Lipschitz function $\bp_H:H \setminus \Sigma_H\to\R^{\ge 0}$ is an  \textbf{\si-transform}, more precisely a Hardy \si-transform, provided it satisfies the following axioms.
\begin{itemize}
  \item  If $H \subset M$ is totally geodesic, then $\bp_H \equiv 0$. Otherwise  we have $\bp_H>0, \bp_H \ge |A_H| $, $\bp_H(x) \ra \infty,$ for  $x \ra p \in \Sigma_H$
and  $\bp_{\lambda \cdot H} \equiv \lambda^{-1} \cdot  \bp_{H}$ for any $\lambda >0$.
  \item  If $H$ is not totally geodesic, and thus $\bp_H>0$, we define the \textbf{\si-distance}  $\delta_{\bp_H}:=1/\bp_H$.
  $\delta_{\bp_H}$ is
  $L_{\bp}$-Lipschitz regular for some constant $L_{\bp}=L(\bp,n)>0$, i.e.,
  \[
  |\delta_{\bp_H}(p)- \delta_{\bp_H}(q)|   \le L_{\bp} \cdot d_H(p,q) \mm{ for any } p, q \in  H \setminus \Sigma \mm{ and any } H \in {\cal{G}}_n.
  \]
  We may choose $\bp_H$ so that $\boldsymbol{L_{\bp}=1}$. Throughout this paper, and the second part \cite{L5}, we make this choice to simplify our computations.
  \item
  If $H_i \in {\cal{H}}_n$, $i \ge 1$, is a sequence converging to the limit space $H_\infty \in {\cal{H}}_n$,
  then $\bp_{H_i}\overset{C^\alpha}  \longrightarrow {\bp_{H_\infty}}$ for any $\alpha \in (0,1)$. For general $H \in {\cal{G}}_n$,
  this \textbf{naturality} holds for blow-ups:  $\bp_{\tau_i \cdot H} \overset{C^\alpha}  \longrightarrow {\bp_{H_\infty}}$, for any sequence $\tau_i \ra \infty$ so that  $\tau_i \cdot H \ra H_\infty \in {\cal{H}}^\R_n$.
    \item For any compact $H \in {\cal{G}}^c$ and any $C^{\alpha}$-regular $(2,0)$-tensor $B$, $\alpha \in (0,1)$, on the ambient space $M$ of $H$ with
    $B|_H \not\equiv -A_H$ there exists a constant $k_{H;B} > 0$ such that
  \[
  \int_H|\nabla f|^2 + |A+B|_H|^2 \cdot f^2 dV \ge k_{H;B}  \cdot \int_H \bp^2 \cdot f^2 dV \ge k_{H;B} \cdot \int_H \frac{f^2}{\dist_H(x, \Sigma)^2}dV.
  \]
\end{itemize}

\noindent  It is worthy to recall that a \emph{totally geodesic} almost minimizer $H$, where  $|A| \equiv 0$, is automatically \emph{regular} since $|A|$ diverges when we approach hypersurface singularities.\\
Any  \si-transform $\bp$  admits a  $C^\infty$-\textbf{Whitney smoothing} $\bp^*$ that still satisfies these axioms except for a slightly weaker form of naturality:
\begin{equation}\label{smot}
c_1 \cdot \delta_{\bp}(x) \le \delta_{\bp^*}(x)  \le c_2 \cdot \delta_{\bp}(x)\quad\mm{and}\quad |\p^\beta \delta_{\bp^*}  / \p x^\beta |(x) \le c_3(\beta) \cdot \delta_{\bp}^{1-|\beta|}(x)
\end{equation}
for constants $c_i(H,\beta) >0$, with $c_i(n,\beta) >0$ for ${\cal{H}}^{\R}_n$, $i=1,\,2,\,3$. Here, $\beta$ is a multi-index for derivatives with respect to normal
coordinates around $x \in H \setminus \Sigma$. Throughout this paper we choose \emph{one fixed pair of a \si-transform and an associated Whitney smoothing } $\bp$ and $\bp^*$ . The precise choices are immaterial for the sequel as the results will not depend on the concrete \si-transform or Whitney smoothing.\\

\textbf{C.2. \si-Pencils}  We can use $\bp$ is to quantify a non-tangential way of approaching $\Sigma$. We define \textbf{\si-pencils} to describe an inner cone condition viewing $\Sigma_H$ as the boundary of $H \setminus \Sigma_H$:
\begin{equation}\label{pen}
\P(z,\omega):= \{x \in H \setminus \Sigma \,|\,\omega \cdot d_H(x,z) < \delta_{\bp} (x)\}
\end{equation}
pointing to $z \in \Sigma$, where $\omega>0$. The angle $\arctan(\omega^{-1})$ is some kind of aperture of $\P(z,\omega)$ relative to $z$. When $H$ is a cone $C$ singular in $0$, we write $C \setminus \{0\} \cong S_C \times \R^{>0}$, for $S_C:=\p B_1(0) \cap C$. Then the pencil $\P(0,\omega)$ is just a subcone $C(U) \subset C$ over some open set $U \subset S_C$.\\
It will be useful to also define the \textbf{truncated \si-pencils} $\TP$. Compared to the \si-pencils $\P$, the $\TP$ are in controllable distance to the singular set and $\D$-maps easily extend to these sets.
\begin{equation}\label{trpen}
\TP(z,\omega,R,r)=\TP_H(z,\omega,R,r):=B_{R}(z) \setminus B_{r} (z) \cap  \P(z,\omega)\subset H.
\end{equation}

\textbf{C.3. \si-Sobolev spaces}  We use dedicated Sobolev spaces, the Hilbert spaces $H^{1,2}_{\bp}(H \setminus \Sigma)$. We recall from \cite[Ch.~5.1]{L1}:
\begin{itemize}
  \item The \textbf{$H^{1,2}_{\bp}$-scalar product}:\, $\langle f,g \rangle_{H^{1,2}_{\bp}(H \setminus \Sigma)}:= \int_{H \setminus \Sigma} \langle \nabla f, \nabla g \rangle + \bp^2 \cdot f \cdot g \, dV $, for $C^2$-functions $f,g$. The $H^{1,2}_{\bp}$-norm associated to this scalar product is written $|f|_{H^{1,2}_{\bp}(H \setminus \Sigma)}$.
  \item The \textbf{\si-Sobolev space} $H^{1,2}_{\bp}(H \setminus \Sigma)$ is the $H^{1,2}_{\bp}(H \setminus \Sigma)$-completion of the subspace of  functions in  $C^2(H \setminus \Sigma)$ with finite $H^{1,2}_{\bp}(H \setminus \Sigma)$-norm. $H^{1,2}_{\bp}(H \setminus \Sigma)$ is a Hilbert space.
\end{itemize}
For non--totally geodesic hypersurfaces $H\in\cal{G}$, we have a vital \textbf{compact approximation} result
\begin{equation}\label{csa}
H^{1,2}_{\bp}(H \setminus \Sigma) \equiv H^{1,2}_{\bp,0}(H \setminus \Sigma):= H^{1,2}_{\bp}\mm{-completion of }C^2_0(H \setminus \Sigma)
\end{equation}
where $C^2_0(H \setminus \Sigma)$ is the space of smooth functions with compact support on $H \setminus \Sigma$.\\

\textbf{D.1. Almost Minimizers} \,  An \textbf{almost minimizer} $H^n$ is a possibly singular
hypersurface looking more and more like an area minimizer the closer we zoom into it. \\
More precisely,  the volume of a ball $B_r(p) \subset H^n$ of radius $r>0$ exceeds that of the area minimizer with the
same boundary by at most $c_H \cdot r ^{n+ 2 \cdot \beta}$, for some constant $c_H >0$. Such an almost minimizer $H^n$ is a $C^{1, \beta}$-hypersurface except for some singular set $\Sigma_H$ of Hausdorff-dimension $\le n-7$.   Sequences of scalings $H_i= \tau_i \cdot H$ of $H$, for some sequence $\tau_i \ra \infty$, for $i \ra \infty$, around a
 given singular point $x \in \Sigma \subset  H^n$,  \emph{flat norm subconverge}\footnote{Saying a sequence \emph{subconverges} means it converges after some possible selection of a subsequence.} to area minimizing tangent cones $C^n \subset \R^{n+1}$.\\

\textbf{D.2. Tameness and $\D$-Maps} \,   In the case of Euclidean area minimizers we know that $H \setminus \Sigma_H$ is smooth and there is even a
 $C^k$-approximation by tangent cones $C$  for any $k \in \Z^{\ge 0}$ in the following sense. For $B_R(q) \cap C \setminus \Sigma_C$, $R>0$, we have from Allard theory:
for any $k \in \Z^{\ge 1}$ and large $i$, $B_R(q_i) \cap H_i$, for suitable $q_i \in H \setminus \Sigma$, is a local $C^k$-section
 \[\Gamma_i :B_R(q) \cap C \ra  B_R(q_i) \cap H_i\subset\nu \mm{ of  the normal bundle } \nu \mm{ of } B_R(q)\cap C\] up to minor adjustments near $\p B_R(q)$ and,
  for $i \ra \infty$, $\Gamma_i$ converges in $C^k$-norm to the zero section, which we identify with $B_R(q) \cap C$. We call  the $C^k$-section $\D := \Gamma_i$
the \emph{asymptotic identification map} or $\D$\textbf{-map} for short. We will later briefly say the $H_i$ compactly \textbf{$\D$-map-converge} to $C$ and write $\D$ when all other details are known from the context. With $\D$-maps $C^k$-functions on $B_R(p_i) \cap H_i$ become comparable to $C^k$-functions on $B_R(p) \cap C$ from an $\D$-map pull-back to $C$ (or $H_\infty$). We use this to specify the class $\cal{G}$ of almost minimizers with a sufficient degree of regularity for the purposes of this series of papers:
an almost minimizer $H$ belongs to $\cal{G}$ provided the following  $C^{k,\gamma}$-\emph{tameness} properties, for some $k \ge 2$, $\gamma \in (0,1)$, hold:
\begin{enumerate}
  \item The (generalized) mean curvature is locally bounded.
  \item $H \setminus \Sigma$ and the $\D$-maps $\Gamma_i$ are $C^{k,\gamma}$-regular with $|\D - id_{C}|_{C^{k,\gamma}(B_R(q) \cap C )}  \ra 0$, for $i \ra \infty$,
  for any given tangent cone $C$,  $p \in \Sigma_H$ and $q \in C$ with $\overline{B_R(q)} \subset  C \setminus \sigma_C$, for some $R>0$.
\end{enumerate}
Throughout this paper we assume $k=5$ (and drop $\gamma$) to be on the safe side. Besides area minimizers, $\cal{G}$ covers cases we typically encounter in scalar curvature geometry or  general relativity like hypersurfaces of prescribed mean
curvature but also cases \emph{not} arising from variational principles, like marginally outer trapped surfaces (= horizons of black holes).\\

 \textbf{D.3. Principal Eigenvalues} \, Let $H\in\cal{G}$ be a non--totally geodesic almost minimizer. Then we use potential theory, from \cite{L2} and \cite{L3}, to study the conformal Laplacian near $\Sigma$.

\begin{itemize}
  \item There exists a finite constant $\tau = \tau(H)>-\infty$ such that for any $C^2$-function $f$ compactly supported in $H\setminus \Sigma$:
 $\int_H  f  \cdot  L_H f  \,  dV \, \ge \, \tau \cdot \int_H \bp^2\cdot f^2 dV$. The largest  such $\tau \in \R$   is the \textbf{principal eigenvalue} $\lambda^{\bp}_{H}$ of $\delta_{\bp}^2 \cdot L$. For any $H\in\cal{G}$, we have $\lambda^{\bp}_H > -\infty$.
\item If, in addition, $\scal_M \ge 0$ and $H \in {\cal{H}}$, then $\lambda^{\bp}_{H}>0$. This is a proper upgrade to the classically known positivity of the ordinary principal eigenvalue of $L$.
\end{itemize}
Throughout this paper we express the eigenfunctions as solutions of the equations $L_{H,\lambda} \, \phi=0$ for  the following   \emph{shifted conformal Laplacian}:
 \begin{equation}\label{v1}
L_{H,\lambda}:=L_H - \lambda \cdot \bp^2, \mm{ for } \lambda  < \lambda^{\bp}_H ,
\end{equation}
for the principal eigenvalue $\lambda^{\bp}_H \mm{ of }{\bp}^{-2} \cdot L$.   For $\lambda = \lambda^{\bp}_{H}$, there is a positive solution of $L_{H,\lambda} \, \phi=0$ which is the counterpart to the first eigenfunction in the smooth compact case. However, the potential theory of $L_{H,\lambda^{\bp}_H}$ is less well-controlled than that of $L_{H,\lambda}$ for \textbf{subcritical eigenvalues} $\lambda < \lambda^{\bp}_{H}$. The condition  $\lambda < \lambda^{\bp}_{H}$ is essential to derive the fundamental asymptotic control over solutions of $L_{H,\lambda} \, \phi=0$.
In the potential theory of these conformal Laplacians $L_{H,\lambda}$ is called \textbf{\si-adapted} when $\lambda  < \lambda^{\bp}_H$. In particular, $L_H$ is
 is \si-adapted when  $\lambda^{\bp}_H>0$, whereas $L_{H,\lambda^{\bp}_H}$ is never \si-adapted.\\

 \textbf{E.1. Induced Solutions} \, We can derive estimates for eigenfunctions on $H$ from  \textbf{induced eigenfunctions} on $C$. We briefly recall how this works:
while we scale by increasingly large $\tau >0$, we observe that $\tau \cdot  \TP_H(p,\omega,R/\tau,r/\tau)$, cf.~C.2., is better and better $C^k$-approximated by the corresponding  truncated \si-pencil in the given tangent cone. This carries over to the analysis on $\P$. When we choose any $p \in \Sigma_H$ then we have for any $\ve > 0$, $R > 1 > r >0$ and some $\omega \in (0,1)$ and any solution $u>0$ of $L(H) \,\phi = 0$ there exists some $\tau^*(L, u,\ve,\omega, R , r, p)>0$ such that for \emph{any} $\tau\ge \tau^*$ there is some tangent cone $C^\tau_p$ with $|\D - id_{C^\tau_p}|_{C^{k,\gamma}(\TP_{C^\tau_p}(0,\omega,R,r))} \le \ve $ and a solution $v>0$, of $L(C^\tau_p) \,\phi = 0$, that can be chosen \textbf{independently}  of $\ve,\omega, R,r$, with
\begin{equation}\label{ffefua}
|u \circ \D / v-1|_{C^{2,\alpha}(\TP_{C^\tau_p}(0,\omega,R,r) )} \le \ve.
\end{equation}
We call such a solution $v$ on $C^\tau_p$ an \textbf{induced solution}.\\

 \textbf{E.2. Minimal Growth} \, A (super)solution $u \ge 0$ of $L_{H,\lambda}\, \phi=0$ has
\textbf{minimal growth}  towards $p \in \widehat{\Sigma}$,  if there is a supersolution $w >0$,  such that $(u/w)(x) \ra 0$, for $x \ra p$, $x \in H \setminus \Sigma$. We say $u$ has \textbf{minimal growth}  towards $W \subset \widehat{\Sigma}$ if $u$ has minimal growth towards any point $p \in W$.
An important result, \cite[Th.3]{L4}, which we use extensively, is that \textbf{minimal growth is inherited} under convergence of the underlying spaces, in particular under blow-ups. We consider three cases:
\begin{itemize}
  \item \textbf{Tangent Cones} \, Let $u >0$ be a supersolution of $L_{H,\lambda}\phi =0$, $\lambda  < \lambda^{\bp}_H $ that is a solution on a  neighborhood $V$ of some $z \in \Sigma$ with \emph{minimal growth towards $V \cap \Sigma$}. Then, if $C$ is a tangent cone of $H$ in $z$, any solution induced on $C$ has \emph{minimal growth towards $\Sigma_C$.}
  \item \textbf{General Blow-Ups} \,  More generally, let $u >0$ be a supersolution of $L_{H,\lambda}\phi =0$, $\lambda  < \lambda^{\bp}_H $ that is a solution on a  neighborhood $V$ of $z \in \Sigma$ with \emph{minimal growth towards} $V \cap \Sigma$ and consider a sequence $s_i \ra \infty$ of scaling factors and a sequence of points $z_i \ra z$ in $\Sigma_{H}$ such that $(s_i \cdot H, z_i)$ subconverges to a limit space $(H_\infty, z_\infty)$ with $H_\infty \in {\cal{H}}^{\R}_n.$  Then the induced solutions on $H_\infty \setminus \Sigma_{H_\infty}$  have \emph{minimal growth towards }$\Sigma_{H_\infty}$.
  \item \textbf{Tangent Cones at Infinity} \,  For $H \in {\cal{H}}^{\R}_n$ the argument of \cite[Th.3]{L3} equally applies to  \emph{tangent cones at infinity}. Recall that they are the possible limit spaces we get from scaling $H$ by some sequence $\tau_i >0$ with  $\tau_i \ra 0$ for $i \ra \infty$. We assume $u >0$ solves  $L_{H,\lambda}\phi =0$, $\lambda  < \lambda^{\bp}_H $, with \emph{minimal growth towards $\Sigma_H$} (but not towards infinity). Then, for any tangent cone $C$ of $H$ at infinity, any solution induced on $C$ has \emph{minimal growth towards }$\Sigma_C$.
\end{itemize}

 \setcounter{section}{2}
\renewcommand{\thesubsection}{\thesection}
\subsection{Hyperbolic Unfoldings} \label{hyd}

The Gromov hyperbolic metric $d_{\bp}$  on hypersurfaces $H \in {\cal{G}}$ is used throughout this paper. We oftentimes use all three geometries $d_H$, $d_{\sima}$ and $d_{\bp}$ in one argument. Here we discuss some basics.

\subsubsection{Canonical Semmes Families} \label{uhy}

\begin{definition}\emph{\textbf{(Gromov Hyperbolicity)}}\label{grh}  A metric space $X$ is \textbf{geodesic}, when any two points can be joined by a geodesic, i.e., a path that
is an isometric embedding of an interval. A geodesic metric space $X$ is $\mathbf{\delta}$\textbf{-hyperbolic,} if all its geodesic triangles are $\mathbf{\delta}$\textbf{-thin} for
some $\delta \ge 0$, that is, the $\delta$-neighborhood of any two sides of such a triangle contains the third side.
 The space $X$ is called \textbf{Gromov hyperbolic} when it is $\delta$-hyperbolic for some $\delta$.
\end{definition}

A \textbf{generalized geodesic ray} $\gamma: I \ra X$ is an isometric embedding of the interval $I \subset \R$ into $X$, where either $I = [0,\infty)$, then $\gamma$ is a proper
geodesic ray, or $I = [0,R]$, for some $R \in (0,\infty)$. Then $\gamma$ is a geodesic arc. When we fix a base point $p \in X$ we can use the hyperbolicity to canonically
identify any $x \in X$ with a (properly) generalized ray $\gamma_x$ with endpoint $\gamma(R) = x$. We extend the definition of such a ray to $I = [0,\infty]$ setting $\gamma(t) = \gamma(R)$ when $t \in [R,\infty]$.

\begin{definition} \emph{\textbf{(Gromov Boundary)}}   Let $X$ be a complete Gromov hyperbolic space. The \textbf{Gromov boundary} $\p_G X$ of $X$ is the set of equivalence classes $[\gamma]$ of geodesic rays, starting from a base point $p \in X$, with two rays being equivalent if they have finite Hausdorff distance.
 \end{definition}

 $\p_G X$  does not depend on the choice of the base point $p$. To topologize $\overline{X}_G = X \cup \p_G X$, we say $x_n \in \overline{X}$ \emph{converges} to $x \in \overline{X}$ if there exist
generalized rays $c_n$ with $c_n(0) = p$ and $c_n(\infty) = x_n$ subconverging (on compacta) to a generalized ray $c$ with $c(0) = p$ and $c(\infty) = x$. The canonical map $X \hookrightarrow \overline{X}_G$ is a homeomorphism onto its image, $\p_G X$ is closed and $\overline{X}_G$  is compact and called the \textbf{Gromov compactification} of $X$.\\

Now we turn to almost minimizers $H \in {\cal{G}}$.  The locally Lipschitz Riemannian metric $\bp^2 \cdot g_H$ and  its $C^\infty$-Whitney smoothing  $(\bp^*)^{2} \cdot g_H$, both defined on $H \setminus \Sigma$, induce the distance functions $d_{\bp_H}$ and $d_{\bp^*_H}$, where we drop the index $H$, and write $d_{\bp}$ and $d_{\bp^*}$, when $H$ is known from the context. We recall  \cite[Theorem 1.11, Cor.~3.6 and Prop.~3.11]{L1}.
\begin{theorem} \emph{\textbf{(Hyperbolic Unfoldings)}} \label{hu} For any non--totally geodesic $H \in {\cal{G}}$,
$(H \setminus \Sigma, d_{\bp})$ and the quasi-isometric $(H \setminus \Sigma, d_{\bp^*})$ are \textbf{complete Gromov hyperbolic spaces} with \textbf{bounded geometry} and we call them \textbf{hyperbolic unfoldings} of  $(H \setminus \Sigma, g_H)$. We have:
\begin{itemize}
  \item The identity map on $H \setminus \Sigma$ extends to \textbf{homeomorphisms}
\[
\widehat{H}\cong\overline{(H \setminus \Sigma,d_{\bp})}_G \cong \overline{(H \setminus \Sigma,d_{\bp^*})}_G \,\mm{ and } \, \widehat{\Sigma} \cong\p_G(H \setminus \Sigma,d_{\bp}) \cong \p_G(H \setminus \Sigma,d_{\bp^*}),
\]
  \item the assignment $(H \setminus \Sigma, d_{\bp})$ to  $H \in {\cal{G}}$ is \textbf{natural}: $d_{\bp_H}$ commutes with the compact convergence of the regular portions of a sequence of underlying $H_k\in {\cal{G}}$ with limit $H_\infty$,
  \item for $H\in{\cal{H}}^\R_n$, $(H \setminus \Sigma, d_{\bp})$ is $\delta(n)$-hyperbolic with $(\sigma_n, \ell_n)$-bounded geometry and $(H \setminus \Sigma, d_{\bp^*})$ is $\delta^*(n)$-hyperbolic of  $(\sigma_n, \ell_n)$-bounded geometry. $\sigma_n $ and $\ell_n$ depend only on $n$.
\end{itemize}
\end{theorem}
\begin{remark}\leavevmode 1. The condition for bounded geometry is this: For a global Lipschitz constant $\ell \ge 1$ and a radius $\varrho>0$, there exists around any  $p\in H\setminus\Sigma$ an $\ell$-bi-Lipschitz chart $\phi_p:B_\varrho(p) \ra U_p$ from the ball $B_\varrho(p)$ in $(H\setminus\Sigma, d_{\bp})$ to an open set $U_p \subset(\R^n,g_{\R^{n}})$. We shall always assume that $0 \in U_p$ and $\phi_{p}(p) = 0$. We say that $(H\setminus\Sigma, d_{\bp})$  has  \emph{$(\varrho,\ell)$-bounded geometry}. For the smooth manifold $(H \setminus \Sigma, d_{\bp^*})$ we have a smooth bounded geometry with a uniform bound on the Riemann tensor and its covariant derivatives up to order $k$. \\
2. For $H \in {\cal{H}}^{\R}_n$, the boundary $\p_G(H \setminus \Sigma,d_{\bp})$ always contains the point at infinity even when $H$ is regular, whereas for regular $H \in {\cal{G}}^c_n$, we have  $\p_G(H \setminus \Sigma,d_{\bp}) \v$. \\
3. We frequently use scaling and blow-up arguments from geometric measure theory. We note that the scalings $\lambda \cdot H$, for some $\lambda >0$, of the original (almost) minimizing geometries $(H,d_H)$ are compensated by $\bp_{\lambda \cdot H} \equiv \lambda^{-1} \cdot  \bp_{H}$, that is, $H$ and $\lambda \cdot H$ have the same hyperbolic unfolding.  \qed
\end{remark}

We start with a geometric application of hyperbolic unfoldings. For this we recall that to handle analysis on rather general metric measure spaces, Semmes has decompiled the classical proof of Poincar\'{e} inequality on $\R^n$ where, in an important step, one uses the presence of uniformly distributed families of curves linking any two given points \cite{Se, He, H-T}.
The abstracted concept is that of \textbf{thick families of curves}, also called \textbf{Semmes families}, satisfying the conditions (i) and (ii) in \ref{sem0} below. Under some reasonable assumptions on the metric space, the presence of Semmes families implies the validity of a Poincar\'{e} inequality. On almost minimizers we have Poincar\'{e}, Sobolev and isoperimetric inequalities \cite{BG}. In this case the presence of Semmes families comes hardly as a surprise, but there is a little extra we get from hyperbolic unfoldings. The hyperbolic geodesics give us \emph{canonically} defined Semmes families on  $(H, d_H)$. The interesting point is that they are still Semmes families relative to $(H, d_{\sima})$ and the (yet to define) minimal factor measure $\mu_{\sima}$. This is proved in Ch.~\ref{sobo} where it is used to derive Poincar\'{e} inequalities for $(H, d_{\sima},\mu_{\sima})$.

\begin{proposition} [Canonical Semmes Families on $H \in {\cal{G}}_n$]\label{sem0} For any $H \in {\cal{G}}_n$, there are constants $C=C(H)>0$, with $C=C(n)>0$ for $H \in {\cal{H}}^{\R}_n$,
and families $\Gamma_{p,q}$ of rectifiable curves $\gamma: I_\gamma \ra H$, $I_\gamma \subset \R$,  joining any two $p,q \in H$ so that:
\begin{enumerate}
 \item For any $\gamma \in \Gamma_{p,q}$: $l(\gamma|_{[s,t]}) < C \cdot d(\gamma(s),\gamma(t))$ for  $s,t \in I_\gamma$.
 \item Each family $\Gamma_{p,q}$ carries a probability measure $\sigma_{p,q}$ so that for any Borel set $A \subset X$, the assignment $\gamma \mapsto l(\gamma \cap A)$  is $\sigma$-measurable with
 \small \begin{equation}\label{tcu1}
 \int_{\Gamma_{p,q}} l(\gamma \cap A) \, d \sigma(\gamma) \le  C \cdot \int_{A_{C,p,q}} \left(\frac{d(p,z)}{\mu(B_{d(p,z)}(p))} + \frac{d(q,z)}{\mu(B_{d(q,z)}(q))}\right) d \mu(z)
 \end{equation}
 \normalsize
 for $A_{C,p,q}:=(B_{C \cdot d(p,q)}(p) \cup B_{C \cdot d(p,q)}(q))\cap A$.
\end{enumerate}
\end{proposition}

\begin{remark}
The family $\Gamma_{p,q}$ is a fibration of a twisted double cone in $(H,d_H)$, directly obtained from \si-uniformity of $H$, which surrounds a geodesic $\gamma_{p,q}$ in the hyperbolic unfolding, see Step~2 below. We call $\gamma_{p,q}$ the \textbf{core} of $\Gamma_{p,q}$. There is a degree of freedom
in the definition of $\Gamma_{p,q}$ we exploit in Ch.\ref{sobo}: we can prescribe the thickness $d$ of the double cone and get families of curve families $\Gamma_{p,q}[d]$ and of measures $\sigma_{p,q}[d]$ for sufficiently small $d >0$, depending only on $H$, only on $n$ when $H \in {\cal{H}}^{\R}_n$. The smaller $d$ becomes the better we can control the elliptic analysis on  $\Gamma_{p,q}[d]$, but the larger $C(H)$ and $C(n)$ have to be chosen. \qed
\end{remark}

\textbf{Proof} \, We show that on $\R^n$ such thick families of curves can be constructed explicitly. Then we use hyperbolic unfoldings to transfer these families to  $(H,d_H)$.\\

\textbf{Step 1} \textbf{(Euclidean Model)} For $x, y \in \R^n$ we consider the hyperplane $L^{n-1}(x,y)$ orthogonal to the line segment $[x, y] \subset \R^n$ passing through the midpoint $m(x,y)$ of $[x, y]$. We consider a ball $B_r:=B_r^{n-1}(m(x,y)) \subset L^{n-1}(x,y)$ of radius $r  \in (0,d(x,y)]$. For any $z\in B_r$, let $\gamma_{z}$ be the unit speed curve from $x$ to $y$ we get when we follow the line segments $[x, z]$ and $[z,y]$. The path space
\begin{equation}\label{ga}
\textstyle\Gamma^{\R^n}_{x,y}=\Gamma^{\R^n}_{x,y}(r):=\{\gamma_{z}\,|\,z\in B_r\}
\end{equation}
is supported on the double cone $\bigcup \{\gamma | \gamma \in \Gamma^{\R^n}_{x,y} \} \setminus \{x,y\}$ with edge $\p B_r$.
There is a canonical $\sqrt{2}$-bi-Lipschitz map  $Q_{\gamma_z}: [x, y] \ra \gamma_z$, $z\in B_r$. The orthogonal projection to $[x, y]$ is the inverse map.
We also define the following \textbf{enveloping rounded double cones} $\textbf{E}^{\R^n}_{x,y}(\varsigma)$:
\begin{equation}\label{ga2}
\textstyle \textbf{E}^{\R^n}_{x,y}(\varsigma) := \bigcup_{z \in [x,y]} B_{\varsigma \cdot \min \{|x-z|,|y-z|\}}(z), \mm{ for some } \varsigma >0.
\end{equation}
Since $\bigcup \{\gamma | \gamma \in \Gamma^{\R^n}_{x,y} \} \setminus \{x,y\} \subset  \textbf{E}^{\R^n}_{x,y}(\varsigma)$ for   $\varsigma \ge 2 \cdot r/d(x,y)$, we henceforth choose $\varsigma := 2 \cdot r/d(x,y)$. We will later use these envelopes to get analytic estimates on $\Gamma^{\R^n}_{x,y}$ by means of Harnack inequalities on $\textbf{E}^{\R^n}_{x,y}(\varsigma)$.
Since $r  \le d(x,y)$ we have some $C_0>0$, independent of $x, y$ and $r$, so that (i) is satisfied for any $\gamma \in \Gamma^{\R^n}_{x,y}$.
Towards (ii), we define a probability measure $\alpha^{\R^n}_{x,y}(r)$ on $\Gamma^{\R^n}_{x,y}(r)$:
\begin{equation}\label{defsem}
\alpha^{\R^n}_{x,y}(W)  :=\mathcal{H}^{n-1}(\{z\in B_r^{n-1}\,|\, \gamma_{z}\in W\})/\mathcal{H}^{n-1}(B_r^{n-1}), \mm{ for any Borel subset } W\subset\Gamma^{\R^n}_{x,y}.
\end{equation}
From the coarea formula \cite[Lemma 2.99]{AFP}, we see that for any Borel set $A\subset \R^{n}$, the function $\gamma\mapsto\ell(\gamma\cap A)$ on $\Gamma^{\R^n}_{x,y}$ is $\alpha^{\R^n}_{x,y}$-measurable and thus expressions on both sides of (\ref{tcu1}) are Borel measures evaluated on $A$.
Now assume that all points in $A$  are closer to $x$ than to $y$. For the distance between corresponding points on any two segments, we have
\begin{equation}\label{dist}
d(s \cdot z_1 + (1-s) \cdot x, s \cdot z_2 + (1-s) \cdot x) \le 2 \cdot r \cdot s, \mm{ for } s \in [0,1] \mm{ and } z_1, z_2 \in B_r.
\end{equation}
From this, the coarea formula gives the following inequality  for the annuli $A_{j} :=A \cap B(x,\ 2^{-j}d(x,\ y))\setminus B(x, 2^{-j-1}d(x, y))$, $j \in \Z^{\ge 0}$:
 \begin{equation}\label{coer}
 \int_{B_r}\ell(\gamma_{z}\cap A)d\mathcal{H}^{n-1}(z)=\sum_{j=0}^{\infty}\ \int_{B_r}\ell(\gamma_{z}\cap A_{j})d\mathcal{H}^{n-1}(z) \le  \sum_{j=0}^{\infty} 2^{(j+2)(n-1)} \mu(A_{j})
 \end{equation}
\[\le 4^{n-1} \cdot \int_{A\cap B(x,d(x,y))}\frac{d(x,z)}{\mu(B_{d(x,z)}(x))} d \mu(z).\]
\medskip

In the last inequality we used that  $\mu(B_{d(x,z)}(x)) =  c_n \cdot d(x,z)^n$, where $c_n >0$ is the volume of the unit ball, to remove the factor $2^{j \cdot (n-1)}$. For $A$ closer to $y$ than to $x$, we argue similarly with the integrand $d(y,z)/\mu(B_{d(y,z)}(y))$ and decomposing a general $A$ into the two components of points closer to $x$ or $y$ we get (\ref{tcu1}) in (ii). For  (i) and (ii)  we choose $C=4^{n-1}/ (c_{n-1} \cdot r^{n-1}) + C_0$.\\

\textbf{Step 2} \textbf{(\si-Uniform Envelopes)}  Since the Gromov compactification $\overline{X}_G$ of  $X=(H \setminus \Sigma,\bp^2 \cdot g_H)$ is homeomorphic to $(H, d_H)$, there is for any two points $p,q \in H$ a hyperbolic geodesic $\gamma_{p,q} \subset X \cup \{p,q\} \subset \overline{X}_G$ that joins these points.
 This choice yields \ref{sem0}(i) for $\gamma_{p,q}$, the core curve in the yet to define family $\Gamma_{p,q}$. Namely, relative to $(H,g_H)$, \emph{each segment }$\widetilde{\gamma}_{x,y} \subset \gamma_{p,q}$ joining two points $x,y \in \gamma_{p,q}$ is \emph{\si-uniform}, more precisely a  $c$-\si-uniform curve for some $c \ge 1$, that is:
  \begin{itemize}
    \item \textbf{Quasi-geodesic:} \, $l(\widetilde{\gamma}_{x,y})  \le c \cdot  d(x,y).$
    \item  \textbf{Twisted double \si-cones:} \, $l_{min}(\widetilde{\gamma}_{x,y})(z) \le c \cdot \delta_{\bp}(z)$ for any $z \in \widetilde{\gamma}_{x,y}$.
  \end{itemize}
where $l_{min}(\widetilde{\gamma}_{x,y})(z):=$ minimum of the lengths of the two subcurves of $\widetilde{\gamma}_{x,y}$ from $x$ to $z$ and from $q$ to $z$. From Prop.~3.11  and  Lemma 3.13 of \cite[Ch.~3.2]{L1}, the constant $c$ depends only on $H$ and only on $n$ for $H \in {\cal{H}}^{\R}_n$. Now we use the Lipschitz continuity of $\delta_{\bp}$ and $\delta_{\bp}\equiv 0$ on $\Sigma$,
\begin{equation}\label{lip}
 |\delta_{\bp_H}(p)- \delta_{\bp_H}(q)|   \le d_H(p,q) \mm{ for any } p, q \in  H \setminus \Sigma \mm{ and any } H \in {\cal{G}}_n.
\end{equation}
This  gives $\delta_{\bp_H}(p)  \le \dist(p, \Sigma)$, i.e., $\bp_H(p)  \ge \dist(p, \Sigma)^{-1}$, and when we multiply (\ref{lip}) by $\bp_H(p) \cdot \bp_H(q)$ we get \cite[Lemma B.2]{L1}:
\begin{equation}\label{lipi}
\bp(q) \le 2 \cdot \bp(p), \mm{ for any } q \mm{ with } d_H(p,q) \le 1/2  \cdot \delta_{\bp_H}(p).
\end{equation}
 Now we use the twisted double \si-cone condition to get a counterpart of $\textbf{E}^{\R^n}_{x,y}(\varsigma) \subset \R^n$, defined in (\ref{ga2}) above, in $(H, d_H)$. The interval $[x,y]  \subset \R^n$ is replaced by the \textbf{core} $\gamma_{p,q} \subset H$. The twisted double \si-cone around $\gamma_{p,q}$ is supported on
 what we will  become the envelope of our Semmes family:
\begin{equation}\label{tube}
\textbf{E}_{p,q}[d]  := \bigcup_{z \in \gamma_{p,q} \setminus \{p,q\}} B_{d \cdot l_{min}(\gamma_{p,q}(z))/c}(z).\\
\end{equation}

We observe that for  $d  \in (0,1)$, we have $\textbf{E}_{p,q}[d] \subset H \setminus \Sigma$ and $\overline{\textbf{E}_{p,q}[d]} \cap \Sigma \subset \{p,q\}$ since
\begin{equation}\label{beg}
\dist(z,\Sigma) \ge \delta_{\bp_H}(z) \ge c^{-1} \cdot  l_{min}(\gamma_{p,q}(z)).
\end{equation}
For any $z \in \gamma_{p,q}$, the balls $B_{d \cdot l_{min}(\gamma_{p,q}(z))/c}(z) \subset B_{d \cdot \delta_{\bp}(z)}(z)$ scaled to unit size admit common bounds on the geometry independent of $z,p,q$. Concretely, when $B_{d \cdot l_{min}(\gamma_{p,q}(z))/c}$ is scaled by $({d \cdot l_{min}(\gamma_{p,q}(z))/c})^{-1}$, we get a ball of radius $1$.
For the present we write it as $(B^\circ_1(z),g^\circ)$. The twisted double \si-cone condition and $\bp_{\lambda \cdot H} \equiv \lambda^{-1} \cdot  \bp_{H}$ show  that  $\bp(z) \le d$ and thus, from (\ref{lipi}) and further restricting $d  \in (0,1/2)$, we get $\bp(x) \le 2 \cdot d$, for any $x  \in B^\circ_1(z)$. \\

Since  $|A| \le \bp$ and $|A|$ bounds the norm of the principal curvatures of $B^\circ_1(z)$  in its scaled ambient space, the Gauss formulas give us uniform bounds on the sectional curvature. In turn, when $H \in {\cal{H}}^{\R}_n$, Allard regularity \cite[Theorem 24.2]{Si1} shows that for any $\ve >0$ and any integer $k \ge 0$ there is a small $d  \in (0,1/2)$ so that the exponential map $\exp_z: B_1(0) \ra B^\circ_1(z)$, for  $B_1(0) \subset T_z(({d \cdot l_{min}(\gamma_{p,q}(z))/c})^{-1} \cdot H)$ in the tangent space in $z$, is a $C^k$-diffeomorphisms and
\begin{equation}\label{pull}
|\exp^*_z(g^\circ) - g_{Eucl}|_{C^k(B_1(0))} \le \ve, \mm{ for the pull-back } \exp^*_z(g^\circ) \mm{ of } g^\circ.
\end{equation}
This carries over to more general $H \in {\cal{G}}_n$, $k=3$, using the tameness assumption of D.2.
In particular,  for $\ve>0$ small enough, the exponential map in any such $z$ is a local bi-Lipschitz map with constant $L \le 2$ on $B_{d \cdot l_{min}(\gamma_{p,q}(z))/c}(0) \subset T_zH$. Note that the Lipschitz constant $L$ does not change under common rescaling of source and image.\\
For the next step, we additionally choose $d \le c/2$. Summarizing, and using \cite[Prop.~2.7 (iv)]{L1} and \cite[Lemma 3.13]{L1}, we see that for given $\ve>0$ and $k\ge 0$, $d$ can be chosen independently of $p,q \in H$ and for $H \in {\cal{H}}^{\R}_n$ the parameter $d$ only depends on $n$.\\

\textbf{Step 3} \textbf{(Semmes Families in $(H,d_H)$)}  The envelope $\textbf{E}_{p,q}[d]$ contains the total space of a fibration by curves that defines the desired Semmes families $\Gamma_{p,q}$ on $(H, d_H)$.\\
Let $\ell$ be the length of $\gamma_{p,q} \subset (H, d_H)$. We parameterize $\gamma_{p,q}$ by arc-length, choose a point $x \in \R^n$ with $|x| =\ell/2$ and the family $\Gamma_{-x,x}(r)$ in $\R^n$ with $r=d \cdot \ell/c\le \ell/2$. We start with the isometry, similar to a parametrization by arc length,
\begin{equation}\label{path}
P_{\gamma_{p,q}}: [-x,x] \ra \gamma_{p,q}, \mm{ with } P_{\gamma_{p,q}}(-x)=p \mm{  and  } P_{\gamma_{p,q}}(x)=q.
\end{equation}
With these choices, we consider \emph{Fermi coordinates} along  $[-x,x] \setminus \{-x,x\}$ on $\textbf{E}^{\R^n}_{-x,x}(\varsigma)$ and $\gamma_{p,q}  \setminus \{p,q\}$ on  $\textbf{E}_{p,q}[d]$ and use them to extend $P_{\gamma_{p,q}}$ to a smooth and $2 \cdot L$-bi-Lipschitz map $\textbf{P}_{\gamma_{p,q}}$ and
we define
\begin{equation}\label{semfam}
\Gamma_{p,q}[d]:=\{\textbf{P}_{\gamma_{p,q}}\circ\gamma \,|\, \gamma \in \Gamma^{\R^n}_{x,y}(r)\}
\end{equation}
with the measure $\alpha_{p,q}[d]$ on $\Gamma_{p,q}[d]$ induced by the measure  $\alpha^{\R^n}_{x,-x}(r)$ on $\Gamma^{\R^n}_{x,-x}(r)$. For the core curve, the quasi-convexity property (i) follows from the $c$-\si-uniformity (and thus the c-quasi-geodesy of any of the  subcurves). This extends to the other curves in $\Gamma_{p,q}[d]$ from the bi-Lipschitz maps $Q_{\gamma_z}: [x, y] \ra \gamma_z$, $z\in B_r$ we used in the definition of the Euclidean model families and
the bi-Lipschitz maps $\textbf{P}_{\gamma_{p,q}}$. For property (ii) we note that the left and the right hand side of (\ref{coer}) are merely changed by (powers of) the bi-Lipschitz constants and they remain to be Borel measures.\qed

\subsubsection{\si-Doob Transforms and Upper $d_{\sima}$-Estimates} \label{bhpdo}

Now we turn to analytic applications of hyperbolic unfoldings.
The main reason for our interest in these unfoldings is that the potential theory of $\delta_{\bp^*}^2 \cdot L_{H,\lambda}$ on $(H \setminus \Sigma, d_{\bp^*})$ is nicely structured towards the Gromov boundary from the work of Ancona, cf.~\cite{An1}, \cite{An2} and \cite{KL}. \\
To formulate associated boundary regularity results for $(H \setminus \Sigma, d_H)$ along $\Sigma$, we briefly look at the standard conformal deformation of the flat unit disk $D$ to the hyperbolic plane.  In this \emph{Poincar\'{e} metric}, the intersections $B \cap D$ of flat Euclidean discs $B$ centered in boundary points $\p D$ become hyperbolic halfspaces.\\ This suggests the following generalization: we consider hyperbolic halfspaces $\textbf{U} \subset (H \setminus \Sigma, d_{\bp})$. They are generally not quite conformally equivalent to the distance balls in $(H,d_H)$, but it turns out that the $\textbf{U}$, and \emph{not} the distance balls in $(H,d_H)$, are the adequate choice for analytic estimates on $(H,d_H)$ towards $\Sigma_H$. To define these halfspaces, we recall that
for any  $p,y, z \in (H \setminus \Sigma,d_{\bp^*})$,  the \textbf{Gromov product} of $y$ and $z$ with respect to $p$ is defined as
\begin{equation}\label{gro}
(y \cdot z)_p:=\frac 12\cdot \left(d_{\bp^*}(p, y) + d_{\bp^*}(p,z)-d_{\bp^*}(y,z)\right).
\end{equation}

This product measures how long two geodesic rays from $p$ to $y$ and from $p$ to $z$  travel together before they diverge (and, adding $\gamma_{y,z}$, form a $\delta$-thin triangle), cf.~\cite[Lemme 2.17]{GH}:
\begin{equation}\label{grrr}
(y \cdot z)_p  \le d_{\bp^*}(p, \gamma_{y,z}) \le (y \cdot z)_p + 4 \cdot \delta.
\end{equation}

We use the product to describe some neighborhood bases \cite[Lemma III.H.3.6]{BH}:

\begin{lemma}[Neighborhood Bases via Gromov Product]\label{nbgp} We define the following subsets for some $z \in \p_G (H \setminus \Sigma,d_{\bp^*})$
and $a>0$:
\begin{itemize}[leftmargin=*]
  \item $\mathrm{U}(z,a):= \{ x \in \p_G (H \setminus \Sigma,d_{\bp^*}) \,|\,\mm{ there are geodesic rays } \gamma_1,\, \gamma_2$ from $p$ to $\gamma_1(\infty) =z,\, \gamma_2(\infty) = x$ such that $\liminf_{t \ra \infty} (\gamma_1(t)\cdot \gamma_2(t))_p \ge a\} \subset \p_G (H \setminus \Sigma,d_{\bp^*})$.
  \item $\mathrm{\mathbf{U}}(z,a):= \{x \in (H \setminus \Sigma,d_{\bp^*}) \,|\, \mm{ there is a geodesic ray } \gamma \mm{ starting from }p \mm{ to }
\gamma(\infty)= z \mm{ such that }\\\liminf_{t \ra \infty} (\gamma(t) \cdot x)_p \ge a\}  \subset (H \setminus \Sigma,d_{\bp^*}).$
\item $\U(z,a):= \mathrm{\mathbf{U}}(z,a) \cup \mathrm{U}(z,a) \subset \overline{(H \setminus \Sigma,d_{\bp^*})}_G$
\end{itemize}
Then the $\mathrm{U}(z,a)$ and $\U(z,a)$ are neighborhood bases of $z$ in their respective space.
 \end{lemma}

There are several essentially equivalent ways to define such halfspaces, cf.~\cite[Lemma 8.3]{BHK}, \cite[Ch.~5.2]{An2}.  \cite{L2} uses the more direct definition $\{x \in H \setminus \Sigma\,|\, \dist\big(x,\gamma([t,a))\big)<\dist\big(x,\gamma((0,t]) \big)\}$. Here we rather use $\textbf{U}(z,a)$. This makes it easier to derive estimates for the minimal factor metrics.\\

Now we reach the central boundary regularity results \cite[Theorem 3.4 and 3.5]{L2}, considering the singular set $\widehat{\Sigma}$  as a \emph{boundary} of its regular complement $H \setminus \Sigma$. We get them by transferring the theory on the hyperbolic unfolding back to $L_{H,\lambda}$ on $(H \setminus \Sigma, g_H)$. As in \cite[Theorem 3.4 and 3.5]{L2} we select uniformly spaced neighborhood bases around points $z \in \widehat{\Sigma}$:
\begin{equation}\label{ui}
\textbf{U}_k(z):=\textbf{U}(z,c \cdot k) \mm{ and }   \U_k(z):=\U(z,c \cdot k) \quad\text{for $k=1,2,\dots$}
\end{equation}
where $c>0$ depends only on the hyperbolicity constant $\delta$ of $(H \setminus \Sigma, d_{\bp^*})$. Since the hyperbolicity constant $\delta$ is already determined from Prop.~\ref{hu} we no longer mention it in our discussions.
\begin{theorem}[Boundary Harnack Inequality]\label{mbhsq}
For  $H \in {\cal{G}}$  and $L_{H,\lambda}$ \si-adapted,  let  $u$, $v>0$ be two supersolutions of $L_{H,\lambda}\, \phi= 0$ on $H \setminus \Sigma$
both solving $L_{H,\lambda}\, \phi= 0$ on $\emph{\textbf{U}}_k(z) \subset  H \setminus \Sigma$, around some $z \in \widehat{\Sigma}$ and for some $k=1,2,\dots$, with minimal growth towards $\mathrm{U}_k(z)$. Then there is a constant $C=C(H)>0$, and $C=C(n)>0$ when $H\in{\cal{H}}^\R_n$, independent of $k$, so that we have:
\begin{equation}\label{fhepq}
u(x)/v(x) \le C \cdot  u(y)/v(y), \mm{ for all } x,\, y \in \emph{\textbf{U}}_{k+1}(z).
\end{equation}
The same inequality holds for $u,v$ considered as solutions of $\delta_{\bp^*}^2 \cdot L_{H,\lambda} \, \phi=0$ on $(H \setminus \Sigma, d_{\bp^*})$.
\end{theorem}

\begin{remark}[Inequalities along Arcs] The validity of (\ref{fhepq}) does not depend on the interpretation of $\textbf{U}_k(z)$ as neighborhoods of the singular point $z$. Technically, it is the relative position of  $\p \textbf{U}_{k+1}(z)$ to $\p \textbf{U}_k(z)$ in a hyperbolic space that yields (\ref{fhepq}). Therefore   (\ref{fhepq}) equally holds for any hyperbolic geodesic arc $\gamma:[0, (m+2) \cdot c+\eta] \ra (H \setminus \Sigma,d_{\bp^*})$, for some $\eta >0$, $m \in \Z^{>0}$, of finite length from $p=\gamma(0)$ to $z= \gamma((m+2) \cdot c+\eta)$ when $p,z \in H \setminus \Sigma$, cf.~\cite[Th\'{e}or\`{e}me 6.1]{An2}, \cite[Lemma 2.4 and Theorem 2.12]{L2}.
The extended description of the halfspaces $\textbf{U}_k(\gamma)$ now reads:
\begin{itemize}[leftmargin=*]
  \item $\textbf{U}_k(\gamma):= \{x \in (H \setminus \Sigma,d_{\bp^*}) \,|\, (\gamma((k+1) \cdot c) \cdot x)_p \ge k \cdot c \}$, for $k=1,2,\dots m.$
\end{itemize}
We note in passing that, in general, $(H \setminus \Sigma,d_{\bp^*})$ is not a visual metric space. That is, a finite arc $\gamma$ does not necessarily admit extensions to complete geodesic rays or approximations by such rays.  \qed
\end{remark}

Typical \emph{supersolutions} satisfying the assumptions of Thm.~\ref{mbhsq} are the minimal Green's function $G(\cdot, p)$, $p \in H  \setminus \Sigma$, the functions $\mathbf{S_\lambda}$ from Example \ref{exa} (\ref{suppa}) and the functions $\Phi$ from Def.~\ref{msge} of minimal factor metrics.  The boundary Harnack inequality shows that around any $z \in \Sigma \subset H$ there is a small neighborhood so that $\Phi \le c \cdot G(\cdot, p)$ for some constant $c>0$.  We show that $c$ is essentially independent of $p$ and $z$.

\begin{proposition}[Harnack Estimate along Geodesics]\label{gphi}    We consider the two cases where
\[\mm{ \emph{(i)} } H \in {\cal{H}}^{\R}_n \mm{ is non--totally geodesic}, \, \mm{ \emph{(ii)} } H \in {\cal{G}}^c_n \mm{ is singular.}\]
Let $p,z \in H  \setminus \Sigma$  be points with $\bp(p) \le a$, for some $a >0$, and $d_H(p,z) \le 1$.
Then, there is a constant $\xi>0$ with $\xi=\xi(n,a)$ in case (i) and  $\xi=\xi(H,\Phi,a)$ in case (ii)  so that
\begin{equation}\label{oba}
\Phi(z) \le \xi \cdot \Phi(p) \cdot G(z,p), \mm{ for any  } p, z \in H \setminus \Sigma.
\end{equation}
\end{proposition}

\textbf{Proof} \, The idea is to link  $p \in H \setminus \Sigma$ and $z \in H$ with a hyperbolic geodesic arc $\gamma_{p,z}$ in $(H \setminus \Sigma_{H}, d_{\bp^*})$, where we append $z$ when $z \in \Sigma$. Then we use the Gromov product to see that (\ref{oba}) can be derived from a combination of the ordinary and the boundary Harnack inequality. This also exploits that $\gamma_{p,z}$ is an $c$-\si-uniform arc in $(H ,d_H)$, for some $c(H)>0$, with $c(n)$ for  $H \in  {\cal{H}}^{\R}_n$:  $l_{min}(\gamma_{p,z}(x)) \le c \cdot \delta_{\bp}(x)$, for any $x \in \gamma_{p,z}$, and $l(\gamma_{p,z})  \le c \cdot  d(p,z).$
\begin{itemize}[leftmargin=*]
 \item When  $z \in H \setminus \Sigma$ with $\bp(z) \le b$, for some $b >0$,  we use the $\bp$-bounds in $p$ and $z$ and a suitable collection ${\cal{B}}(p,z)$ of balls from the envelope $\textbf{E}_{p,z}[d]$ to build a \emph{finite} ball cover of $\gamma_{p,z}$
  \begin{equation}\label{fbc}
 B_{d \cdot \delta_{\bp}(p)}(p),  B_{d \cdot \delta_{\bp}(z)}(z), B_{d \cdot l_{min}(\gamma_{p,z}(x_j))/c}(x_j),  j=1,...m,
  \end{equation}
for some integer $m(a,b,c)$ that merely depends on $a,b,c$ but not on $p,z$ or $H$. Since $d  \in (0,1/2)$, the balls $B_{3/2 \cdot d \cdot \delta_{\bp}(p)}(p)$,... etc.\ still belong to $H \setminus \Sigma$. Thus we have a Harnack constant $\kappa(H,a,b,c,d)>0$ resp. $\kappa(n,a,b,c,d)>0$ that applies to any solution $v>0$ of  $L_{H,\lambda} \, \phi=0$ on $H \setminus \Sigma$ on each of the balls in ${\cal{B}}(p,z)$ and we get $\sup_{\bigcup \{B \in {\cal{B}}(p,z)\}}  v \le \kappa^{m+2} \cdot \inf_{\bigcup \{B \in {\cal{B}}(p,z)\}}  v$.\\
In the case of $H  \in  {\cal{G}}^c_n$, the set where $\Phi$ is not a proper solution is included in some compact set $K \subset H \setminus \Sigma$ where we get such a constant $\kappa$ from the continuity of $\Phi$. Thus we always have
\begin{equation}\label{gloh}
\textstyle \sup_{\bigcup \{B \in {\cal{B}}(p,z)\}} \Phi \le \kappa^{m+2} \cdot \Phi(p).
\end{equation}
\end{itemize}
Now we turn to  case (i) and assume that (\ref{oba}) does not hold on ${\cal{H}}^{\R}_n$. Then we have a sequence $H_i \in  {\cal{H}}^{\R}_n$ converging to some $H_\infty \in  {\cal{H}}^{\R}_n$ and pairs of points $p_i, z_i$ with $\bp(p_i) \le a$  and $d_{H_i}(p_i,z_i) = 1$ (since scaling by a constant $\ge 1$ only decreases $\bp(p_i)$) so that the $p_i$ $\D$-map-converge to some $p_\infty \in H_\infty \setminus \Sigma_{H_\infty}$ with $\bp(p_\infty) \le a$. Moreover, we can assume that the $H_i$ are equipped with minimal growth solutions $\Phi_i>0$, with $\Phi_i(p_i)=1$, and minimal Green's function $G_i$ both compactly $\D$-map-converging to
 $\Phi_\infty>0$, with $\Phi_\infty(p_\infty)=1$ and the minimal Green's function $G_\infty$ on $H_\infty$, respectively, with
 \begin{equation}\label{diver}
 \Phi_i(z_i)/ G_i(z_i,p_i) \ra \infty, \mm{ for } i \ra \infty.
 \end{equation}
We may assume that
there are hyperbolic geodesic arcs $\gamma_i$ of length $L_i \in (0,\infty]$ in $(H_i \setminus \Sigma_{H_i}, d_{\bp^*})$ from $p_i$ to $z_i$ $\D$-map-converging to a hyperbolic geodesic arc $\gamma_\infty$ from $p_\infty$ to some $z_\infty$ with  $d_{H_\infty}(p_\infty,z_\infty) = 1$.  Then we have $z_\infty \in \Sigma_{H_\infty}$, otherwise the $z_i$ would subconverge to some $z_\infty \in H_\infty \setminus \Sigma_{H_\infty}$ and hence we would have $\bp(z_i) \le \bp(z_\infty)+1 <\infty$ for $i$ large enough and could infer from (\ref{gloh}) that $\Phi_i(z_i)/ G(z_i,p_i)$ remained bounded for $i \ra \infty$ contradicting (\ref{diver}). Thus
$\delta_{\bp}(z_i) \ra 0$ and $L_i \ra \infty$ for $i \ra \infty$. Now we consider a hyperbolic geodesic ray $\gamma_i^*$  from $p_i$ to $\infty$ and compare it with $\gamma_i$:
\begin{itemize}[leftmargin=*]
  \item  we have $l_{H_i}(\gamma_i^*)=\infty$ and, hence, the $c$-\si-uniformity of $\gamma_i^*$  shows  that  for $\gamma_\infty(0)=p_\infty$,
  \begin{equation}\label{e1}
  l_{H_i}(\gamma_i^*|_{[0,t]}) \le c \cdot \delta_{\bp}(\gamma_i^*(t)) \mm{ for any } t \ge 0,
  \end{equation}
  and, thus, for $d_{H_i}(\gamma_i^*(0),\gamma_i^*(t)) \ge 1/2$ we have $\delta_{\bp}(\gamma_i^*(t)) \ge (2 \cdot c)^{-1}$
  \item from $|\delta_{\bp}(a)- \delta_{\bp}(b)|   \le d_{H_i}(p,q)$, $a, b \in  {H_i} \setminus \Sigma_{{H_i}}$ we get for $\ve >0$,  large enough $i$ and $s \in (0,L_i)$:
      \begin{equation}\label{e2}
      \delta_{\bp}(\gamma_i(s)) \le d_{H_i}(\gamma_i(s),\gamma_i(L_i)) + \ve.
      \end{equation}
\end{itemize}
Since the $\gamma_i$ compactly $\D$-map-approximate the hyperbolic geodesic arc $\gamma_\infty$ from $p_\infty$ to $z_\infty$ outside $z_\infty$,
we get from (\ref{e1}) and (\ref{e2}) for any $\ve>0$ and sufficiently large $i$ that $\gamma_i^*$ diverges from  $\gamma_i$ after some
uniformly upper bounded time $s_0$ so that
\begin{equation}\label{drei}
\delta_{\bp}(\gamma_\infty(s)) \le  d_{H_\infty}(\gamma_\infty(s),z_\infty) \le (2 \cdot c)^{-1} \cdot \eta,\mm{ for some small } \eta>0 \mm{ and } s \ge s_0.
\end{equation}
Namely, we see from (\ref{grrr})  and the relation \cite[Lemma 3.12]{L1}
\begin{equation}\label{riww}
\big|\log\delta_{\bp}(a) - \log\delta_{\bp}(b)\big| \le d_{\bp}(a,b), \mm{ for any } a, b \in  {H_\infty} \setminus \Sigma_{{H_\infty}},
\end{equation}
 that for small $\eta$, there is a $k_0$ so that, for large $i$ and then independent of $i$,  the  $\textbf{U}_k(\gamma_i) = \{x \in (H \setminus \Sigma,d_{\bp^*}) \,|\, (\gamma_i((k+1) \cdot c) \cdot x)_{p_i} \ge k \cdot c \}$, for $k=1,2,\dots m(i)$, $m(i) \ra \infty$ when $i \ra \infty$, satisfy
\begin{equation}\label{end}
\gamma_i^*  \cap \textbf{U}_k(\gamma_i) \v \mm{ for } k \ge k_0.
\end{equation}
Thus we can apply the boundary Harnack inequality (\ref{fhepq}) to $\Phi$ and $G$ on these $\textbf{U}_k(\gamma_i)$ where both functions are solutions of  minimal growth towards the extension of $\textbf{U}_k(\gamma_i)$ to $\Sigma_{H_i}$. The point about the independence of $k_0$ from $i$ is that the relation (\ref{riww}) and the estimate (\ref{grrr}) also show that there is a constant $b >0$ independent of $i$ so that $\bp(q_i) \le b$ for the first point $q_i \in \gamma_i$ we reach in $\gamma_i  \cap \p \textbf{U}_k(\gamma_i)$ starting from $p_i$.\\
Thus following $\gamma_i$ from $p_i$ to $z_i$, we use the chain (\ref{gloh}) of ordinary Harnack estimates until we reach $q_i$. After passing $q_i$ we continue with the boundary Harnack estimate (\ref{fhepq}) on $\textbf{U}_k(\gamma_i)$ to get an upper estimate for $ \Phi_i(z_i)/ G_i(z_i,p_i)$, independent of $i$. This contradiction settles case (i).\\
In case (ii) argue similarly using chains (\ref{gloh}) of Harnack estimates on $H_i=a_i \cdot H$ with $a_i:=d_{H}(p_i,z_i)^{-1}\ge 1$. In place of (\ref{end}) we find linking geodesics $\gamma_i$ and a  $k_0$, independent of $i$, so that
\begin{equation}\label{end2}
K  \cap \textbf{U}_k(\gamma_i) \v \mm{ for } k \ge k_0,
\end{equation}
and derive a contradiction as in case (i). \qed

To derive estimates for the minimal Green's function of an elliptic operator $L$ we use \emph{Doob transforms}, also called \emph{h-transforms} $L^{h}$. They are standard transformations from stochastic analysis for operators $L$ on function spaces: for smooth $h>0$ on $H \setminus \Sigma$, defined as $L^{h}:= h^{-1}Lh$, i.e., for smooth $g$ we set $L^{h}g(x):= h^{-1}(x)\, L (h \cdot g )(x).$ We choose $L= L_{H,\lambda}$ and $h=\delta^{-(n-2)/2}_{\bp^*}$.

\begin{proposition} \emph{\textbf{(\si-Doob Transforms)}}\label{doob}  For any $H \in {\cal{H}}_n$ we define the \textbf{\si-Doob Transform} $L^{\sima}=L^{\sima}_{H,\lambda}$  on
$(H \setminus \Sigma, d_{\bp^*})$ of $L=L_{H,\lambda}$ to be the following Schr\"odinger operator on suitably regular functions (for the present, we consider $C^2$-functions $\phi$):
\begin{equation}\label{lsdef}
L^{\sima} \phi:= \delta_{\bp^*}^2 \cdot  L^{\delta^{-(n-2)/2}_{\bp^*}}_{H,\lambda}\phi = \delta^{(n+2)/2}_{\bp^*} \cdot  L ({\delta^{-(n-2)/2}_{\bp^*}}_{H,\lambda} \cdot \phi) =
\end{equation}
\small \[ \delta_{\bp^*}^2   \cdot  \Big(-\Delta_{{\bp^*}^{^2}\cdot g_H} \phi  +
\Big(\frac{n-2}{4 (n-1)} \scal_H - \lambda \bp ^{2} - \delta_{\bp^*}^{(n-2)/2} \Delta_{g_H} \delta_{\bp^*}^{-(n-2)/2} \Big)\phi \Big).\]
\normalsize
\begin{itemize}[leftmargin=*]
  \item For a function $v>0$, we set $u:=\delta^{-(n-2)/2}_{\bp^*} \cdot v$. Then  $u^{4/(n-2)} \cdot g_H=v^{4/(n-2)} \cdot  {\bp^*}^{^2}\cdot g_H$ and
\begin{equation}\label{3g}
u \mm{  solves } L \phi =0  \Leftrightarrow v \mm{  solves }  L^{\sima} \phi=0 \, , \,
L u \ge c \cdot  \delta^{-2}_{\bp^*} \cdot u \Leftrightarrow L^{\sima} v = \delta^{(n+2)/2}_{\bp^*}\cdot  L u \ge c \cdot  v,
\end{equation}
for some $c>0$. That is, $L$ is  \si-adapted $\Leftrightarrow L^{\sima}$ is adapted and weakly coercive.
  \item The operator $L^{\sima}$ is symmetric  on $(H \setminus \Sigma, d_{\bp^*})$ and, thus, the minimal Green's function $G^{\sima}$ of $L^{\sima}$ on $(H \setminus \Sigma, d_{\bp^*})$ satisfies $G^{\sima}(x,y)=G^{\sima}(y,x), \mm{ for } x \neq y \in H \setminus \Sigma.$
  \item The minimal Green's functions $G$ of $L$ on $(H \setminus \Sigma, d_H)$ and $G^{\sima}$ of $L^{\sima}$ on $(H \setminus \Sigma, d_{\bp^*})$ satisfy
\begin{equation}\label{gre}
G^{\sima}(x,y) =  \delta_{\bp^*}^{(n-2)/2}(x)  \cdot   \delta_{\bp^*}^{(n-2)/2}(y) \cdot G(x,y), \mm{ for } x \neq y \in H \setminus \Sigma.
\end{equation}
\item Assume that there is some $v>0$ so that for some $c>0$: $L^{\sima} v  \ge c \cdot  v$. Then we have:
\begin{enumerate}[leftmargin=*]
    \item for any singular $H \in {\cal{G}}^c_n$ there are constants $\beta(H,c),\alpha(H,c),\sigma(H,c)>0$ and
    \item for any non--totally geodesic $H \in {\cal{H}}^{\R}_n$ there are $\beta(n,c),\alpha(n,c),\sigma(n,c)>0$ so that
\end{enumerate}
\begin{equation}\label{ges}
 G^{\sima}(x,y)\le \beta  \cdot \exp(-\alpha \cdot d_{\bp^*}(x,y)),\mm{ for } x,y  \in H \setminus \Sigma \mm{ and } d_{\bp^*}(x,y)> 2 \cdot \sigma.
\end{equation}
\end{itemize}
\end{proposition}

\textbf{Proof} \, The formula (\ref{lsdef}) for $L^{\sima}$ and the relations  (\ref{3g}) follow from straightforward computations. The equivalence of the \si-adaptedness of $L$ and the adaptedness of $L^{\sima}$ follows easily from the definitions cf.~\cite[Prop.~3.3]{L2} for details.\\
The symmetry of $L^{\sima}$ is seen as follows: for any two functions  $u,v \in C_0^\infty(H \setminus \Sigma,\R)$ we have, writing $dV^{\sima}=\delta^{-n}_{\bp^*} dV$ for the volume element associated to ${\bp^*}^{^2}\cdot g_H$:
\small \[\int_{H \setminus \Sigma} u \cdot L^{\sima} v  \, dV^{\sima} = \int_{H \setminus \Sigma} u \cdot  \delta^{(n+2)/2}_{\bp^*} \cdot  L_{H,\lambda} ({\delta^{-(n-2)/2}_{\bp^*}} \cdot v)  \cdot \delta^{-n}_{\bp^*}  dV =
\int_{H \setminus \Sigma} L_{H,\lambda} ({\delta^{-(n-2)/2}_{\bp^*}} \cdot u) \cdot {\delta^{-(n-2)/2}_{\bp^*}}  \cdot v \, dV.\]
\normalsize

For the transformation law (\ref{gre}), we recall that one may characterize $G^{\sima}(x,y)$ as the unique function that solves
 \[L_x^{\sima} \int_{H \setminus \Sigma}  G^{\sima}(x,y) u(y) dV^{\sima}(y)= u(x) \mm{ for any } u \in C_0^\infty(H \setminus \Sigma,\R).\]
We insert the right hand side to check this identity from the identity for $G$ relative to $L$:
\[\delta^{(n+2)/2}_{\bp^*}(x) L_x \int_{H \setminus \Sigma} \delta^{(n-2)/2}_{\bp^*}(y) \cdot  G(x,y) \cdot u(y) \cdot \delta^{-n}_{\bp^*}(y) dV(y)=\delta^{(n+2)/2}_{\bp^*}(x) u(x) \delta^{-(n+2)/2}_{\bp^*}(x)=u(x).\]
 Finally, we get  (\ref{ges}) from the work of Ancona \cite[Ch.~2]{An1}, carried out in more detail in \cite[Proposition 2.13]{KL}, and $G^{\sima}(q,p)=G^{\sima}(p,q)$. These results apply due to the uniform ellipticity (more precisely described as adaptedness in \cite{An1} and \cite{KL}) of $L^{\sima}$ on  $(H \setminus \Sigma, (\bp^*)^{2} \cdot g_H)$ and the weak coercivity assumption we made for $L^{\sima}$. \qed

\begin{remark}\label{bhp3}\, 1. The Whitney smoothed  $\bp^*$ satisfies the approximative naturality axiom (\ref{smot}). This implies that estimates expressed in terms of $\bp^*$ can equally be used for $\bp$ and vice versa. Analytically, $\bp^*$  is merely used intermediately to unfold $(H \setminus \Sigma, g_H)$ to a \emph{smooth} manifold. This is a common way to bypass regularity issues from the weaker Lipschitz regularity of $\bp$ in the potential theoretic analysis. The appearance of  Whitney smoothings cancels out on the way back  from the unfolding to $(H \setminus \Sigma, g_H)$, cf.~\cite[Ch.~3]{L1} for details.\\
2. The constants $\beta(H,c),\alpha(H,c),\sigma(H,c)>0$, or $\beta(n,c),\alpha(n,c),\sigma(n,c)>0$ when $H \in {\cal{H}}^{\R}_n$, in Prop.~\ref{doob} \eqref{ges} come from the potential theory on the hyperbolic unfolding. They are determined from the hyperbolicity and bounded geometry constants  for the hyperbolic unfolding in Prop.~\ref{hu}.  $H$ and any scaled version $\lambda \cdot H$, $\lambda>0$
have the same unfolding since  $\bp_{\lambda \cdot H} \equiv \lambda^{-1} \cdot  \bp_{H}$. Thus these constants remain unchanged under scalings and blow-ups of $H$.\qed
\end{remark}

\begin{corollary}[Upper $d_{\sima}$-Estimates]\label{glen} For  $H \in \cal{G}$, $p \in H \setminus \Sigma$ and   $B=B_{2 \cdot \sigma}(p)$  measured relative to the $(\sigma,\ell)$-bounded geometry $d_{\bp^*}$,  we have: the distance of any two $x,y  \in H \setminus B$ with respect to $G(\cdot,p)^{4/(n-2)}  \cdot g_H$ is upper bounded by some constant $D_H>0$ and $D_n >0$ when $H \in {\cal{H}}^{\R}_n$.
 \end{corollary}
\textbf{Proof} \, We apply Prop.~\ref{doob} (\ref{ges}) and (\ref{gre}) to estimate $G$. Any two points $x,y  \in (H \setminus \Sigma, d_{\bp^*})$ can be joined by a hyperbolic geodesic arc $\gamma_{x,y}$, with $\gamma_{x,y}(0)=x$. By definition these arcs are parameterized by arc length and we have $d_{\bp^*}(x,y)=l_{d_{\sima}}(\gamma_{x,y})$.  We have $t=d_{\bp^*}(\gamma_{p,y}(t),p)$, $g_H(\dot \gamma(t),\dot \gamma(t)) = \delta_{\bp^*}^{2}(\gamma(t))$ and hence for any $z=\gamma_{x,y}(t)$ with $d_{\bp^*}(p,z) \ge 2 \cdot \sigma$,
\begin{equation}\label{upest}
\; G(z,p)^{\frac{4}{n-2}}  \cdot g_H(\dot \gamma,\dot \gamma) =  (\bp^*)^2(p) \cdot  G^{\sima}(z,p)^{\frac{4}{n-2}}
\end{equation}
\[\le  (\bp^*)^2(p) \cdot c \cdot \left(\exp(-\alpha \cdot d_{\bp^*}(z,p))\right)^{\frac{4}{n-2}}.\]
Since $\left(\exp(-\alpha \cdot t)\right)^{2/(n-2)}$ is integrable on $\R^{\ge 0}$, the integral of $\left(\exp(-\alpha \cdot d_{\bp^*}(\cdot,p))\right)^{2/(n-2)}$ along $\gamma_{x,y} |_{[2 \cdot \sigma,d_{\bp^*}(y,p))}$ is \emph{uniformly upper bounded} for any $y \in H \setminus (\Sigma \cup B)$.\\ For those parts of such a hyperbolic geodesic passing through $B$ we get uniform upper bounds using the bounded geometry on $B=B_{2 \cdot \sigma}(p)$ and the fact that $G^{\sima}(\cdot, p) \le c_0$ on $\p B_\sigma(p)$, for some $c_o>0$ that is independent of $p \in H \setminus \Sigma$. On ${\cal{H}}^{\R}_n$ the constant $c_0$ is even independent of $H$ \cite[Prop.~2.8]{KL}.  \qed

The estimate (\ref{upest}) is used quite frequently throughout this paper. It also gives upper topological bounds for our minimal factors. This is counterbalanced by lower bounds we get from the $\lambda>0$-condition. This condition allows us to use the Bombieri--Giusti Harnack inequality \cite{BG} e.g. in Lemma \ref{umia} and Prop.~\ref{lcmg} below.\\

\begin{lemma}[Diameter Bounds]\label{db} For  $H \in \cal{G}$ we have the following estimates:
\begin{enumerate}
\item For $H \in {\cal{G}}^c$ we have $\diam(H,d_H)<\infty$ and $\diam \big(\widehat{H \setminus \Sigma}, \widehat{d_{\sima}}\big)<\infty$.
\item For $H \in \cal{G}$ and balls $(B_r(z), d_H)$, we have $\diam\big(\widehat{B_r(z) \setminus \Sigma},\widehat{d_{\sima}}\big) \ra 0$ for $r \ra 0$.
\item  There is a continuous map $\widehat{I_H}: (H,d_H) \ra (\widehat{H \setminus \Sigma}, \widehat{d_{\sima}})$ that canonically extends the identity map
$id: (H \setminus \Sigma,d_H) \ra (H \setminus \Sigma, d_{\sima})$.
\end{enumerate}
\end{lemma}

\textbf{Proof} \,  We start with (i). The assertion that, for $H \in {\cal{G}}^c_n$, the \emph{intrinsic} diameter $\diam(H,d_H)$ is finite is \cite[Theorem 1.8]{L1}.
We choose a basepoint $p \in H \setminus \Sigma$ and consider the ball $B=B_{2 \cdot \sigma}(p)$ with radius measured relative to $d_{\bp^*}$, where $\sigma$ is the bounded geometry constant for $d_{\bp^*}$ as used in Prop.~\ref{doob} \eqref{ges}.  From the boundary Harnack inequality \ref{mbhsq} and the compactness of $\Sigma$ we have a constant $a>0$ so that $\Phi(x) \le a \cdot G(x,p)  \mm{ for any } x \in H \setminus (\Sigma \cup B).$
From this,  \ref{glen} implies that $\diam \big((H \setminus \Sigma,d_{\sima})\big) < \infty$ and we infer for the metric completion  $\diam \big((\widehat{H \setminus \Sigma}, \widehat{d_{\sima}})\big)<\infty$.

For (ii) we only need to consider $z \in \Sigma$. The boundary Harnack inequality can be used to upper estimate $\Phi$ by $G(\cdot,p)$ near $z$. The argument of part (i) shows that the integral of $\left(\exp(-\alpha \cdot d_{\bp^*}(\cdot,p))\right)^{2/(n-2)}$ along $\gamma_{x,y} |_{[2 \cdot \sigma,d_{\bp^*}(y,p))}$ is uniformly upper bounded for any $y \in H \setminus (\Sigma \cup B) \cap B_R(z)$ and fixed $R>0$. From this  we have $\diam\big(\widehat{B_r(z) \setminus \Sigma},\widehat{d_{\sima}}\big) \ra 0$ for $r \ra 0$.

For (iii) we use that $\Sigma \subset (H,d_H)$ is homeomorphic to $\p_G (H \setminus \Sigma,d_{\bp^*}) \setminus \{\infty\}$, where  $\p_G$ denotes the \emph{Gromov boundary}, to define a
canonical continuous map
\begin{equation}\label{idd}
I_H: (H,d_H) \ra (\widehat{H \setminus \Sigma}, \widehat{d_{\sima}}).
\end{equation}
From (\ref{upest}) we see that $\diam \big((\textbf{U}(z,a),  d_{\sima})\big) \ra 0$ for  $a  \ra \infty$.
This shows that  $\bigcap_{a >0}\overline{\textbf{U}(z,a)}$, where $\overline{\textbf{U}(z,a)}$ is the completion in $(\widehat{H \setminus \Sigma}, \widehat{d_{\sima}})$, contains exactly one point $z_{\bp} \in \widehat{H \setminus \Sigma} \setminus (H \setminus \Sigma)$, the new singular set, and we define
\begin{equation}\label{ext}
I_H: \Sigma \ra  \widehat{H \setminus \Sigma} \setminus (H \setminus \Sigma), \mm{ by }   I_H(z) :=   z_{\bp},  \mm{ for } z \in \Sigma.
\end{equation}
$I_H$ extends the identity map on $H \setminus \Sigma$ to a map  $\widehat{I_H}: (H,d_H) \ra (\widehat{H \setminus \Sigma}, \widehat{d_{\sima}})$.
From the previous discussion we also see, for a converging sequence $a_k \ra a$ in $(H,d_H)$, that $d_{\sima}(\widehat{I_H}(a_k), \widehat{I_H}(a)) \ra 0$ for $k \ra \infty$. That is, $\widehat{I_H}$ is a \emph{continuous} map. \qed

As a preparation to show that $\widehat{I_H}$ is a \emph{homeomorphism} we derive the following coarse estimate.

\begin{lemma}[Growth Rate for $H \in {\cal{H}}^{\R}_n$]\label{umia} For any non--totally geodesic $H \in {\cal{H}}^{\R}_n$ and any $\lambda \in (0,\lambda^{\bp}_H)$ there is some constant $\beta(H,\lambda)>0$ so that for the Euclidean ball $B_R(0)$ and $R>0$ large enough we have:
\begin{equation}\label{div}
\dist_{d_{\sima}}(\p B_R(0) \cap H \setminus \Sigma, \{0\}) \ge R^\beta.
\end{equation}
 \end{lemma}

In our applications we are not interested in the actual value of $\beta>0$ but in the consequence that, for $R \ra \infty$, the $d_H$-distance spheres of radius $R$ do not have $d_{\sima}$-accumulation points. We need this since, at this stage, $\p_{\sima}$ for $H \in {\cal{H}}^{\R}_n$ could contain points at infinity of $(H \setminus \Sigma, g_H)$. \\

\textbf{Proof} \, In the cone case $H=C$, we know from \cite[Th.~4.4]{L3} that in terms of polar coordinates $ (\omega,r) \in \p B_1(0) \cap C \setminus \Sigma \times \R^{>0}\cong C \setminus \Sigma$:
\begin{equation}\label{sep}
\Phi_C(\omega,r) = \psi(\omega) \cdot r^{\alpha_+}   \mm{ for } \textstyle \alpha_+ = - \frac{n-2}{2} + \sqrt{ \Big( \frac{n-2}{2} \Big)^2 + \mu}<0,
\end{equation}
where $\mu(C,\lambda)>-(n-2)^2/4$ is the ordinary (unweighted) principal eigenvalue of an associated elliptic operator on $\p B_1(0) \cap C $.
Now we recall the inheritance result \cite[Th.~3.13]{L3} for cones in $\mathcal{SC}_n$. Since this space is compact, we get a constant $A[n,\lambda]>0$ depending only on $n$ and $\lambda$ so that
\small \begin{equation}\label{exx}
 -\frac{n-2}{2} <-A[n,\lambda] < \alpha_+(C,\lambda) <0, \mm{ for any } C\in \mathcal{SC}_n,
\end{equation}
\normalsize
that is, $-A[n,\lambda] \cdot 2/(n-2) >-1$ and $(r^{-A[n,\lambda]})^{2/(n-2)}$ is \emph{not integrable} over $[1,\infty)$.\\
For a general $H \in {\cal{H}}^{\R}_n$, we get using tangent cones at infinity, cf. E.2 above,  that for any $\omega >0$, there are a constant $k>0$ and some radius $\rho>0$ so that
 \begin{equation}\label{divv}
r^{-A[n,\lambda]} \le k \cdot \Phi(x) \mm{ for } x \in \P(0,\omega) \setminus B_\rho(0) \mm{ with } r = d_H(x,0).
\end{equation}
Thus we get a constant $k^*$ so that for $r \ge \rho$:
\begin{equation}\label{a}
r^{-n} \cdot r^{n-1} \cdot r^{-A[n,\lambda]} \cdot r \le k^* \cdot  \Vol(B_{r}(z))^{-1} \cdot  \int_{B_{r}(z) \cap H \setminus \Sigma} \Phi \, dV
\end{equation}
Now we note that $\scal_H=-|A|^2$ and since $\lambda>0$, we have
\begin{equation}\label{ml}
\Delta_H \Phi = -\Big(\frac{n-2}{4 (n-1)} \cdot |A|^2 + \lambda \cdot \bp^2\Big) \cdot \Phi \le 0,
\end{equation}
that is, $\Phi$ is superharmonic on  $H \setminus \Sigma$. Therefore, the Bombieri--Giusti Harnack inequality \cite[Th.~6]{BG} applies to $\Phi$. Indeed $\Phi$ meets the requirements of that Harnack inequality since we have  $\Delta_H \Phi \le 0$ not only on $B_{r}(z) \cap H \setminus \Sigma$ but $\Phi$ is also a weak  supersolution on $B_{r}(z) \cap H$. This follows from a standard cut-off argument that uses that the Hausdorff codimension of $\Sigma$ is $>2$. This argument is carried out in detail e.g.~in \cite[p.334]{Si2}.\\
From this Harnack inequality we get a constant $c>0$
independent of $H \in {\cal{H}}^{\R}_n$  and of  $r \ge \rho$ so that
 \begin{equation}\label{bghs}\,
r^{-A[n,\lambda]} \le k^* \cdot  \Vol(B_{r}(z))^{-1} \cdot  \int_{B_{r}(z) \cap H \setminus \Sigma} \Phi \, dV  \le k^* \cdot  c \cdot  \inf_{{B_{r}(z) \cap H \setminus \Sigma}}\Phi.
\end{equation}
This  and (\ref{exx}) yield some $\beta>0$ so that for $R>0$ large enough,
 \begin{equation}\label{divv}
\dist_{d_{\sima}}(\p B_R(0) \cap H \setminus \Sigma, \{0\}) \ge R^\beta.
\end{equation}
 \qed
With these results we can prove the following.
\begin{proposition}[Topology of Minimal Factor Metrics]\label{lcmg} For any $H \in \cal{G}$ we have:
\begin{enumerate}
\item For any $z \in H$ and $r>0$, there are an outer $d_{\sima}$-radius $0<\rho_\mathrm{out}(r,z)<\infty$ and an inner $d_{\sima}$-radius $0<\rho_\mathrm{inn}(r,z)<\infty$ so that
 $(B_{\rho_\mathrm{inn}}(z) , d_{\sima}) \subset (B_r(z), d_H)  \subset (B_{\rho_\mathrm{out}}(z), d_{\sima})$.
 \item  $(\widehat{H \setminus \Sigma}, \widehat{d_{\sima}(\Phi)})$ is homeomorphic to $(H,d_H)$. Thus we can write it as $(H,d_{\sima})$. In particular, we have that the \textbf{singular set} of $(H,d_{\sima})$,
\begin{equation}\label{poh}
\Sigma_{\sima}:=  \widehat{H \setminus \Sigma} \setminus (H \setminus \Sigma) \subset (\widehat{H \setminus \Sigma}, \widehat{d_{\sima}}),
\end{equation}
is homeomorphic to $\Sigma \subset (H,d_H)$.
\end{enumerate}
Writing $(B_a(p), d_H)  \subset (B_b(p), d_{\sima})$ means the set-theoretic inclusion $B_a(p) \subset B_b(p)$ where the radii are measured relative to $d_H$ and $d_{\sima}$, respectively.
\end{proposition}
\textbf{Proof} \, For (i) we first note that $(B_r(z), d_H)  \subset (B_{\rho_\mathrm{out}}(z), d_{\sima})$ is just Prop.~\ref{db}(ii). To show the existence of $\rho_\mathrm{inn}$, we prove that for $r>0$ small enough, there is some $c>0$ so that  $\Phi \ge c$  on $B_r(z)\cap H \setminus \Sigma$.
From this we can choose $\rho_\mathrm{inn}= c^{2/(n-2)} \cdot r$. To get that lower estimate, we recall the Gau\ss-Codazzi equation: $
 |A_H|^2 + 2\Ric_M(\nu,\nu)  =  \scal_M - \scal_H +(\tr A_H)^2,$
  where $\tr A_H$ is the mean curvature of $H$. In the non-trivial case $z\in\Sigma$, there is an $m >0$ and for any  constant $k>0$ there is a small $r>0$ so that:
$|\Ric_M(\nu,\nu)|, |\scal_M|,|\tr A_H| \le m,$ whereas $\bp \ge k \mm{ on } B_{r}(z) \cap H \setminus \Sigma.$
The bound on $|\tr A_H|$ is part of the chosen class of almost minimizers in D.2. Thus for large $k\gg 1$ and hence small $r>0$, any solution $\Phi>0$ of $L_{H,\lambda} \phi =0$  is superharmonic on  $B_{r}(z) \cap H \setminus \Sigma$:
\begin{equation}\label{ml}
\Delta_H \Phi = \Big(\frac{n-2}{4 (n-1)} \cdot \scal_H - \lambda \cdot \bp^2\Big) \cdot \Phi \le 0.
\end{equation}
Almost minimizers share their regularity theory with area minimizers. In particular, we get the same blow-up limits. This implies the validity of the Bombieri--Giusti $L^1$--Harnack inequality
\cite[Theorem~6, p.~39]{BG} for any $H \in \cal{G}$. As in the argument of Lemma \ref{umia}, the Hausdorff codimension of $\Sigma$ is $>2$ and, hence, $\Phi$ meets the requirements of that Theorem.\\
Now we recall that intrinsic and extrinsic distances are equivalent: there is a constant $c(H)\in (0,1)$ such that for any $p,q \in H$: $c \cdot  d_{g_{H^n}}(p,q) \le  d_{g_{M^{n+1}}}(p,q) \le d_{g_{H^n}}(p,q)$ \cite[Corollary 2.10]{L1}. Thus \cite[Theorem~6, p.~39]{BG}, which uses extrinsic distances, applies in the following form: for some small $r_H>0$ and $r \in (0,r_H)$ (we may choose $r_H=\infty$ in the case $H \in {\cal{H}}^{\R}_n$), we have
\begin{equation}\label{bgha}
0< a  \cdot  \int_{B_{r}(z) \cap H \setminus \Sigma} \Phi \, dV=:c  \le \inf_{{B_{r}(z) \cap H \setminus \Sigma}}\Phi,\, \mm{ for some }a = a(L_{H,\lambda},r)>0.
\end{equation}
This shows that the map $I_H: \Sigma \ra \Sigma_{\sima}$ from (\ref{ext})  is \emph{injective}.\\

$I_H$ is also \emph{surjective} since any point $x \in \Sigma_{\sima}$ is the $d_{\sima}$-limit of a sequence of points $x_k \in H \setminus \Sigma$ and we can choose hyperbolic geodesic arcs $\gamma_{p,x_k}$ from $p$ to the $x_k$.
A subsequence of these arcs converges to a hyperbolic geodesic ray that represents some point $\overline{x} \in \Sigma$ and then we see that $\{x\}=\bigcap_{a >0}\overline{\U(\overline{x},a)}$, that is, $I_H(\overline{x})=x$.\\

Finally, we show that $\widehat{I_H}^{-1}$ is \emph{continuous}. For compact $H$, any closed set is mapped onto a compact and thus closed set.  For $H \in {\cal{H}}^{\R}_n $, we use Lemma~\ref{umia} to argue similarly. For any $R>0$ and any closed set $A \subset H$, the set  $\widehat{I_H}(A \cap \overline{B_R(0)})$ is again closed. Since (\ref{div}) shows that the sequence $\widehat{I_H}(\p B_i(0) \cap H)$, $i =1,2,...$, has no accumulation points in $(H,d_{\sima})$, we infer that $\widehat{I_H}(A)$ is also closed. This implies that $\widehat{I_H}$ is a \emph{homeomorphism}.
\qed

\begin{remark} We observe that $(H,d_{\sima})$  is a \emph{geodesic metric space}. For  any two $p,q \in (H,d_{\sima})$ we can find $p_i,q_i  \in H \setminus \Sigma$ with $p_i \ra p$, $q_i \ra q$, for $i \ra \infty$. For small and smoothly bounded neighborhoods $U$ of $\Sigma$ with $p_i,q_i \in H \setminus U$ we have, from \cite[Prop.2.1]{L1}, a path $\gamma_U \subset H \setminus U$ joining  $p_i,q_i \in (H,d_{\sima})$. Relative to the intrinsic metric on $H \setminus U$, we can assume that $length(\gamma_U)=d(p_i,q_i)$. From  \ref{db} and BV-compactness results for curves \cite [Theorem 4 in Section 4.5]{SG} we observe that for neighborhoods $U_k \supset U_{k+1}$ shrinking to $\Sigma$, i.e. $\bigcap_k U_k =\Sigma$, there is a subsequence of the $\gamma_{U_k}$ converging to some curve $\gamma_{p_i,q_i} \subset  (H,d_{\sima})$ with $length(\gamma_{p_i,q_i})=d(p_i,q_i)$ in $(H,d_{\sima})$ and, in turn, for $i \ra \infty$, there is a subsequence of the $\gamma_{p_i,q_i} \subset  (H,d_{\sima})$ converging to the wanted geodesic $\gamma_{p,q} \subset  (H,d_{\sima})$ that links $p$ and $q$.
\qed
\end{remark}

\subsubsection{Codimension Estimate for $\Sigma_{\sima} \subset (H,d_{\sima})$} \label{sdm}

The (partial) regularity theory for any almost minimizer $H \in {\cal{G}}$ within some smooth ambient manifold $M^{n+1}$ says that $H$ is smooth except for a singular set $\Sigma_H$ that has Hausdorff codimension $\ge 8$ relative to $M^{n+1}$, cf.~\cite{F},  \cite[Ch.\,11]{Gi}. We extend this estimate to the singular set of $(H,d_{\sima})$. However, this is \emph{not} an obvious application of Federer's estimate for $(H,d_H)$ cf. Remark \ref{sma} below. We rather imitate Federer's  argument in the class of minimal factor geometries.

\begin{remark}[Bishop Deformers] \label{sma}  The identity from $(H, d_H)$ to $(H, d_{\sima})$ is not Lipschitz regular since $\Phi$ diverges towards $\Sigma$. However, the Hausdorff dimension is only a bi-Lipschitz invariant. To illustrate this issue recall that there is an (actually H\"{o}lder regular) homeomorphism $\phi: \R^2 \ra \R^2$ mapping $S^1$ to the Koch snowflake $K \subset \R^2$ with Hausdorff dimension $\ln4/ \ln 3$, cf.~\cite{Bi}. It is conceivable that the deformations of \cite[Th.1.1]{Bi} can be adjusted to show that if $\Sigma \subset (H, d_H)$ has Hausdorff dimension $a$, for some $0<a<n$, that, for any $b$ with $a<b<n$, there is some \emph{non-minimal growth} solution $\omega_b>0$ of $L_{H,\lambda} \phi=0$ so that the singular set of the completion of  $(H \setminus \Sigma,\omega_c^{4/(n-2)} \cdot g_H)$ has Hausdorff dimension $\ge b$.\qed
\end{remark}

\vspace{-0.3cm}

We recall some basic concepts and formulate them for the metric space $(H,d_{\sima})$.

\begin{definition}[Hausdorff Measure and Dimension]\label{haus} For some $H \in {\cal{G}}$, let $A \subset (H,d_{\sima})$, $k \in [0,\infty)$ and $\delta \in (0,\infty]$. Then we set
\begin{equation}\label{hkd}
{\mathbb{H}}^\delta_k(A):= \inf  \Big\{ \sum_i \diam(S_i)^k \, \Big| \,  A \subset \bigcup_i S_i, S_i \subset H, \diam(S_i) < \delta \Big\},
\end{equation}
\vspace{-0.4cm}
\begin{equation}\label{hk}
{\mathbb{H}}_k(A):= \lim_{\delta \ra 0} {\mathbb{H}}^\delta_k(A) = \sup_\delta {\mathbb{H}}^\delta_k(A).
\end{equation}
${\mathbb{H}}_k(A)$  is the $k$-dimensional \textbf{Hausdorff measure} of A. The infimum of all $k$ so that ${\mathbb{H}}_k(A)=0$ is the \textbf{Hausdorff dimension} $\dim_{\mathbb{H}} (A)=\dim_{\mathbb{H}} (A \subset (H,d_{\sima}))$ of $A$ as a subset of $(H,d_{\sima})$.
\end{definition}

\begin{remark} \label{ball} 1. In the literature the definition (\ref{hkd}) oftentimes contains a multiplicative gauging constant to keep  (\ref{hkd}) consistent with the Riemannian volume of  Euclidean unit balls. We are primarily interested in the Hausdorff dimension of $\Sigma$ and drop this constant.\\
2. For the same reason we henceforth only consider \textbf{open ball covers}, that is, in (\ref{hkd}) we assume that $S_i=B_{r_i}(p_i)$ for some $p_i \in H$ and $r_i>0$. The resulting  measures (occasionally called spherical Hausdorff measures) are equivalent since any $S_i$ is contained in an open ball $B_{diam (S_i)}(z_i)$, for some suitably chosen $z_i \in H$, and $diam(B_{diam (S_i)}(z_i)) \le 2 \cdot diam (S_i)$.
\qed
\end{remark}

\begin{proposition}[Basic Properties of ${\mathbb{H}}^\infty_k(A)$]\label{hausi} For every $A \subset (H,d_{\sima})$, we have ${\mathbb{H}}^\infty_k(A)=0$ if and only if ${\mathbb{H}}_k(A)=0$.
For ${\mathbb{H}}_k$-almost all $x \in A$, we have
\begin{equation}\label{hki}
\limsup_{r \ra 0} \,  {\mathbb{H}}^\infty_k(A \cap B_r(x))/r^k \ge 1.
\end{equation}
\end{proposition}

\textbf{Proof} \,  For subsets $A \subset \R^n$ this is  \cite[Lemma 11.2 and Proposition 11.3]{Gi}. These results are direct consequences of the definitions of
${\mathbb{H}}^\infty_k(A)$ and ${\mathbb{H}}_k(A)$. They do not use the regular
Euclidean structure of the underlying space and equally apply to  $(H,d_{\sima})$. \qed

We also notice that for $A$ compact and  $\delta \in (0,\infty]$, we have ${\mathbb{H}}^\delta_k(A)<\infty$  even when $A$ has Hausdorff dimension $>k$. The definition of ${\mathbb{H}}^\delta_k$ also readily implies the existence of open neighborhoods $U(A,k,\eta,\delta)$ of $A$, for any $\eta>0 $, so that ${\mathbb{H}}^\delta_k(U) \le  {\mathbb{H}}^\delta_k(A) + \eta$. Moreover, in non-compact cases with ${\mathbb{H}}^\delta_k(A)=\infty$, we also write ${\mathbb{H}}^\delta_k(A)>0$ to keep the notation consistent. \\

\textbf{Hausdorff dimension of $\Sigma$} \, We extend Federer's estimate for the dimension of $\Sigma$ relative to $(H,d_H)$ to the case where $\Sigma \cong \Sigma_{\sima}$ is interpreted as a subset of $(H,d_{\sima})$.

\begin{theorem}[Partial Regularity of $(H,d_{\sima})$]\label{haus} The Hausdorff dimension of $\Sigma$ relative to the minimal factor metric $(H,d_{\sima})$ is $\le n-7$.
\end{theorem}

 To this end, we note that \cite[Ch.3]{L3}, summarized in E.2, and the distance estimates Prop.~\ref{lcmg}(i), providing the extension to the metric completions, translate into a blow-up theory for minimal factor metrics that strongly resembles that for the original almost minimizing geometries.
\begin{theorem}\emph{\textbf{(Blow-Ups)}}\label{blooa}
For $H \in {\cal{G}}$  we consider $(H, d_{\sima}(\Phi_H))$ and any singular point $p \in \Sigma_H$. Then we get the following \textbf{blow-up invariance:}
Any sequence  $(H, \tau_i \cdot d_{\sima}(\Phi_H))$ scaled around $p$ by some sequence $\tau_i \ra \infty$, $i \ra \infty$, subconverges and the limit of any converging subsequence is $(C, d_{\sima}(\Phi_C))$ for some tangent cone $C$ of $H$ in $p$.
\end{theorem}

We only need to specify the notion of convergence. It has two layers:  a subsequence of $H_i= \tau_i \cdot H$, for some sequence $\tau_i \ra \infty$, for $i \ra \infty$, around a
 given singular point $x \in \Sigma \subset  H^n$ converges to an area minimizing tangent cone $C^n \subset \R^{n+1}$. The flat norm convergence becomes a compact $C^k$-convergence over regular subsets of $C$ expressed in terms of a $C^k$-convergence of $\D$-maps as in D.2, section \ref{bcn}. Then the sections $\Phi_H \circ \Gamma_i$ compactly $C^k$-converge to  $\Phi_C$ on $C$ after normalizing the value of the $\Phi_H \circ \Gamma_i$  in a common base point  in $C \setminus \Sigma_C$.
The fact that the Martin boundary of $L_{H,\lambda}$ for $\lambda < \lambda^{\bp}_H$ on $H \in {\cal{H}}^{\R}_n$ has exactly one point at infinity and again  Prop.~\ref{lcmg}(i) show:
 \begin{theorem}\emph{\textbf{(Euclidean Factors)}} \label{euclf}  For any non--totally geodesic $H \in {\cal{H}}^{\R}_n$  there is a \textbf{unique}$^\cs$ space $(H, d_{\sima}(\Phi_H))$. For $C \in \mathcal{SC}_{n}$, the associated space  $(C, d_{\sima}(\Phi_C))$ is invariant under scaling around $0 \in C$, that is, it is again a cone.
 \end{theorem}

In the following we also use details from the proofs of \ref{blooa} and \ref{euclf} we cite when we need them. We start with the following auxiliary result.

\begin{lemma}[Corresponding Balls]\label{cbc} For $H \in {\cal{G}}$ and some $p \in \Sigma_H$ we consider a tangent cone $C$ in $p$ we identify with $0 \in \R^{n+1}$. For any ball $B^C_r(q) \subset B_1(0) \cap C$ of radius $r \in (0,1)$ in $(C,d_{\sima})$  and any $\delta >0$, there is a ball $B^{H_j}_r(q_i)\subset B_1(0) \subset (\tau_j \cdot H,d_{\sima})$ so that, considered as subsets of the original spaces $(C,d_C)$ and  $(\tau_j \cdot H,d_{\tau_j \cdot H})$ locally embedded in $\R^{n+1}$, we have
\begin{equation}\label{de1}
B^{H_j}_r(q_j) \ra B^C_r(q) \mm{ in flat and, from this, in Hausdorff norm, for } j \ra \infty.
\end{equation}
We call the balls $B^{H_j}_r(q_j)$ asymptotically \textbf{corresponding} to $B^C_r(q)$.
\end{lemma}

\textbf{Proof} \,  From D.1 in section \ref{bcn}, we have a flat norm convergence of some open sets $O_i \subset \tau_j \cdot H$ to $B^C_r(q)$. This also implies Hausdorff  convergence from lower volume estimate for the difference set we have from \cite[Prop.5.14]{Gi}. This applies to $C$ and $H$ since they are (asymptotically) area minimizing in $\R^{n+1}$. For the (almost) minimizing geometries on $C$ and $H$ we have for any given $\eta>0$ and $i$ large enough that  the $O_i$ are almost isometric to $B^C_r(q)$ via $\D$-maps outside an $\eta$--distance tube $U_\eta(\Sigma_C)$ around $\Sigma_C$ cf. D.2.\\
For the deformed geometries $(C,d_{\sima})$ and $(H,d_{\sima})$ we recall from E.1 and E.2  that the
minimal growth solutions  converge smoothly to that on $C$ outside such  $\eta$--distance tubes. Now we reuse the length estimate (\ref{upest}) and remark \ref{bhp3}.2 to see that the tube size of the conformally deformed $\eta$--distance tube also, and uniformly in $j$, shrinks to zero when $\eta \ra 0$. Similarly the size of the complements the $\D$-map images of these distance tubes in
$O_i$  shrinks to zero when $i \ra \infty$. From this there are (not necessarily singular) points $q_i \in H_j$ so that (\ref{de1}) holds for the size of the difference set between $O_i$ and $B^{H_j}_r(q_j)$ converges to  zero when $i \ra \infty$. \qed

We use this to upper estimate the ${\mathbb{H}}^\infty_k$-measure of $\Sigma_H \subset (H,d_{\sima})$.

\begin{proposition}[Measure under Blow-Ups]\label{mup}  For $H \in {\cal{G}}$ converging under scaling by some sequence $\tau_j \ra \infty$, for $j \ra \infty$, to a tangent cone $C$ of $H$ in $p \in \Sigma$, which we identify with $0 \in C$, and any radius $R>0$, relative to the minimal factor metrics, we have
\begin{equation}\label{lim1}
{\mathbb{H}}^\infty_k(\Sigma_C \cap \overline{B_R(0)}) \ge \limsup_j {\mathbb{H}}^\infty_k(\tau_j \cdot \Sigma_H \cap \overline{B_R(0)}).
\end{equation}
\end{proposition}

\textbf{Proof} \, We  cover $\Sigma_H \cap \overline{B_1(0)}$ by finitely many balls $B_i \subset (C,d_{\sima})$ so that
\begin{equation}\label{lim2}
{\mathbb{H}}^\infty_k(\Sigma_C \cap \overline{B_1(0)}) >  \sum_i \diam(B_i)^k - \ve.
\end{equation}
 We start with a variant of an
argument used in \cite{F},  \cite[Ch.\,11]{Gi}. Let $\tau_j \cdot H$, scaled around a basepoint $p \in H$, compactly converge to a tangent cone $C$. Then a cover of $\Sigma \cap  K \subset C \cap  K$, for some compact $ K \subset \R^{n+1}$,  by open subsets of $\R^{n+1}$ also covers $\Sigma_{\tau_j \cdot H} \cap K$, for $j$ large enough and the diameter bounds for the covering sets (asympotically) carry over from $C$ to $\tau_j \cdot H$. In the case of $(H,d_{\sima})$ we need to evaluate the diameter with respect to the minimal factor metrics.\\
From Lemma~\ref{cbc} we find for any $\delta >0$ some $j_\delta >0$ so that for $j\in(0,j_\delta)$ there is a family of balls $B_i^\delta \subset \tau_j \cdot (H,d_{\sima})$ covering $\tau_j \cdot \Sigma_H \cap \overline{B_1(0)}$ with $\diam(B^\delta_i) \le (1 + \delta) \cdot \diam(B_i)$, hence
\begin{equation}\label{lim3}
{\mathbb{H}}^\infty_k(\tau_j \cdot \Sigma_H \cap \overline{B_1(0)})  \le \sum_i \diam(B^\delta_i)^k.
\end{equation}
Summarizing, we have $\limsup_j {\mathbb{H}}^\infty_k(\tau_j \cdot \Sigma_H \cap \overline{B_1(0)}) \le (1 + \delta)^k \cdot \left({\mathbb{H}}^\infty_k(\Sigma_C \cap \overline{B_1(0)})  + \ve\right)$.
For $\ve \ra 0$ and $\delta \ra 0$, the claimed estimate  (\ref{lim1}) follows.   \qed

\begin{corollary}[Cone Reduction]\label{cre}  For $H \in {\cal{G}}$, assume that  ${\mathbb{H}}_k(\Sigma_H)>0$, for some $k$. Then there  exists a tangent cone $C$ in ${\mathbb{H}}_k$-almost every point $x \in \Sigma_H$  such that ${\mathbb{H}}_k(\Sigma_C)>0$.
\end{corollary}

\textbf{Proof} \, From (\ref{hki}) in Prop.~\ref{hausi} we can find for ${\mathbb{H}}_k$-almost every point $x \in \Sigma_H$   a sequence of radii $r_i \ra 0$ for $i \ra \infty$ such that ${\mathbb{H}}^\infty_k(A \cap B_{r_i}(x)) \ge r_i^k$.  Then we get the claim from  \ref{mup}. \qed

It is obvious that the dimensions of the singular sets of an area minimizing cone $C^n \subset \R^{n+1}$ and that of the product cone $\R \times C^n \subset \R^{n+2}$ differ by one. In the area minimizing case, this is used in the inductive reduction argument for the codimension estimate \cite[Th.~11.8]{Gi}. Minimal factor metrics on $\R^m \times C^{n-m}$ are no longer Riemannian products. We recall that the eigenfunction $\Phi_{\R^m \times C^{n-m}}$ with minimal growth towards any point of $\Sigma$ is unique$^\cs$, i.e. up to multiplication by a positive constant, from \cite[Th.~3]{L2}. Hence, $\Phi_{\R^m \times C^{n-m}}$ reflects the symmetries of $\R^m \times C^{n-m}$. Following \cite[Prop.~4.6]{L3}, we can write $\Phi_{\R^m \times C^{n-m}}$ in cylindrical coordinates $x=(z, r,\omega) \in \R^m \times \R^{>0} \times (\p B_1(0) \cap C \setminus \Sigma_C)$  for  $r=r(x)=\dist(x,\R^m \times \{0\})$:
\[\Phi_{\R^m \times C^{n-m}}(\omega,r,z) = \psi(\omega) \cdot r^{\alpha_+} \mm{ for some } \alpha_+ <0 \mm{ on }  \R^m \times C^{n-m} \setminus \Sigma_{\R^m \times C^{n-m}}.\]

Thus $(\R^m \times C^{n-m},d_{\sima})$ is invariant under \emph{translations} in $\R^m$-direction and under \emph{scalings} around points in $\R^m \times \{0\}$. From this we determine the Hausdorff dimension of $\R^m \times \{0\}$.

\begin{lemma}[Hausdorff Dimension of Axes]\label{hausx}  Within $(\R^m \times C^{n-m}, d_{\sima})$, $m \ge 1$,  the $m$-planes $\R^m \times \{y\}$, $y \in C$ have Hausdorff dimension $m$. Moreover, for any $k>0$ and $y \in C$ we have ${\mathbb{H}}^\infty_k( [0,1]^m \times \{0\})>0$ if and only if  ${\mathbb{H}}^\infty_k( [0,1]^m \times \{y\})>0$ .
\end{lemma}

\textbf{Proof} \, We use Prop.~\ref{lcmg}  and the translation invariance of $(\R^m \times C^{n-m}, d_{\sima})$ in $\R^m$-direction to choose $\rho>0$ so that the balls $B_{\rho}(p_i)$,  $p_i \in \Z^m$, measured relative $(\R^m \times C^{n-m}, d_{\sima})$, cover $\R^m \times \{0\} \in \R^m \times C^{n-m}$. Now we consider the Euclidean lattices $2^{-j} \cdot \Z^m \subset \R^m$, $j =1,2,...$ and observe that due to the scaling invariance of $(\R^m \times C^{n-m}, d_{\sima})$,  the $B_{2^{-j} \cdot \rho}(2^{-j} \cdot p_i)$,  $p_i \in \Z^m$, cover $\R^m \times \{0\} \in \R^m \times C^{n-m}$. The restriction of the lattice $2^{-j} \cdot \Z^m$ to the unit cube $[0,1]^m \subset \R^m$ contains $2^{j \cdot m}$ points (up to lower orders along $\p  ([0,1]^m)$). For these balls, we get for any $\alpha >m$:
\begin{equation}\label{konv1}
{\mathbb{H}}^\infty_\alpha([0,1]^m) \le \sum_{\{p_i \in [0,1]^m \cap 2^{-j} \cdot \Z^m \}}  \diam(B_{2^{-j} \cdot \rho}(2^{-m} \cdot p_i))^\alpha  =  2^{j \cdot m} \cdot 2^{-j \cdot \alpha} \cdot \rho^\alpha \ra 0 ,\mm{ for } j \ra \infty.
\end{equation}
Thus $\dim_{\mathbb{H}} ([0,1]^m \subset (H,d_{\sima})) \le m$. In turn, since $\alpha_+ <0$ there is some $\delta_0$ so that for $\delta \in (0,\delta_0)$ any sets of diameter  $\le \delta$ that intersects
$\R^m \times \{0\}$ also has diameter $\le \delta$ when computed relative to  $d_H$. Thus we infer from the area minimizing case where $\dim_{\mathbb{H}} ([0,1]^m \subset (H,d_H)) = m$
that  $\dim_{\mathbb{H}} ([0,1]^m \subset (H,d_{\sima})) \ge m$.
More generally, the translation that maps $\R^m \times \{0\}$ onto $\R^m \times \{y\}$, $y \in C$ is  \emph{bi-Lipschitz} in terms of $d_{\sima}$. To see this, we use again that $d_{\sima}$ is $\R^m$-translation invariant and that the $d_{\sima}$-distance of any two points $x,y  \in \R^m \times \{0\}$ remains finite. This follows again from \ref{glen} (and \ref{db}). From a comparison of ball covers, defined as above, we see that, up to this constant (to the power of $m$) we get the same Hausdorff measure estimates as for  $[0,1]^m  \times \{0\} \subset  \R^m \times \{0\}$ also for $[0,1]^m  \times \{y\} \subset  \R^m \times \{y\}$ and vice versa.  \qed

\begin{lemma}[Radial Singularities]\label{radix}  For $C^{n-m} \in \mathcal{SC}_{n-m}$, assume that $0 \varsubsetneq \Sigma_C$, that is, $C^{n-m}$ is singular not only in $0$ and that for some $k>0$: ${\mathbb{H}}^\infty_k( \Sigma_{\R^m \times C^{n-m}})>0$. Then there is some ball $B \subset \R^m \times C^{n-m}$ with $\overline{B} \cap (\R^m \times \{0\}) \v$ so that ${\mathbb{H}}^\infty_k(B \cap \Sigma_{\R^m \times C^{n-m}})>0$.
\end{lemma}

\textbf{Proof} \, Assume there is no such ball. Then we have, using a suitable countable ball cover, that ${\mathbb{H}}^\infty_k(\Sigma_{\R^m \times C^{n-m}} \setminus \R^m \times \{0\})=0$. This means ${\mathbb{H}}^\infty_k(\R^m \times \{0\})>0$, but from this, Lemma~\ref{hausx} also shows that ${\mathbb{H}}^\infty_k([0,1]^m\times \{0\})>0$ for any $y \in  \Sigma_{[0,1]^m\times C^{n-m}}$,  a contradiction.\qed

\textbf{Proof of Theorem \ref{haus}} \, From Lemma~\ref{radix} and Cor.~\ref{cre}, there are a point $p \in  B \cap \Sigma_{\R^m \times C^{n-m}} \setminus \R^m \times \{0\}$ and a tangent cone $C^*$ in $p$  such that ${\mathbb{H}}_k(\Sigma_{C^*})>0$. Since $p \notin \R^m \times \{0\}$ we know that $C^*$ can be written as $C^* =\R^{m+1} \times C^{n-m-1}$. We iterate this argument until we reach some cone $C^{\cs} =\R^{m+l} \times C^{n-l}$ where $C^{n-l}$ is singular only in $0$. The value ${n-l}$ may depend on the chosen sequence of blow-up points  but we know from the fact that hypersurfaces  in ${\cal{H}}_n$
for $n \le 6$ are regular that $n-l \ge 7$. Since we have ${\mathbb{H}}_k(\Sigma_{C^{\cs}})>0$, we get from \ref{hausx} that $k \le n-7$. \qed

\setcounter{section}{3}
\renewcommand{\thesubsection}{\thesection}
\subsection{Ahlfors Regularity and Semmes Families} \label{dom}

We estimate the eccentricity of $d_{\sima}$-distance balls relative to $d_H$-distance balls. This is used to derive the Ahlfors regularity of $(H,d_{\sima},\mu_{\sima})$ and, in the next chapter, also to study the canonical Semmes families relative to $(H,d_{\sima},\mu_{\sima})$.

\subsubsection{Eccentricity of $(H,d_{\sima})$} \label{es}

The qualitative result \ref{lcmg} says that for any $r>0$, $p \in H$, there are $\kappa \ge 1$, $\rho>0$ with
\begin{equation}\label{incl}
(B_\rho(p),d_{\sima}) \subset (B_r(p), d_H) \subset (B_{\kappa \cdot \rho}(p),d_{\sima}).
\end{equation}
We enhance the arguments to control the relations between these radii quantitatively. The \textbf{local eccentricity} $\boldsymbol\vartheta$ of $d_{\sima}$ relative to $d_H$ for the ball $(B_r(p), d_H)$ is
\begin{equation}\label{kappa}
\vartheta(H,p,r,\Phi):=\inf \{\kappa \ge 1 \,|\,  (\ref{incl}) \mm{ holds for at least one }\rho>0\}.
\end{equation}
This infimum actually is a minimum since $(B_\rho(p),d_{\sima})$ and  $(B_r(p), d_H)$ are open.  In turn, this $\vartheta$ determines a unique $\varrho(\vartheta)>0$ that satisfies (\ref{incl}). The parameters $\vartheta$ and $\varrho$ depend differently on the gauging of $\Phi$, that is, on choosing a positive multiple of $\Phi$:
\begin{equation}\label{gauge}
\vartheta(\lambda \cdot \Phi) = \vartheta(\Phi) \mm{ but } \varrho(\lambda \cdot \Phi) = \lambda^{2/(n-2)} \cdot \varrho(\Phi), \mm{ for } \lambda >0.
\end{equation}

\begin{proposition}[Eccentricity of $(H,d_{\sima})$]\label{radii}   We consider the two cases where
\[\mm{ \emph{(i)} } H \in {\cal{H}}^{\R}_n \mm{ is non--totally geodesic}, \, \mm{ \emph{(ii)} } H \in {\cal{G}}^c_n \mm{ is singular.}\]
Then we have a common upper bound for $\vartheta(H,p,r,\Phi)$ for any minimal factor metric:
\begin{equation}\label{kap}
\Theta_n \ge \vartheta(H,p,r,\Phi), \mm{ depending only on } n,
\end{equation}
for any $r>0$ in case \emph{(i)} and for sufficiently small $r>0$ in case \emph{(ii)}.
\end{proposition}

\textbf{Proof} \,  For $H \in {\cal{H}}^{\R}_n$, $\tau >0$ and $v \in H$, we also have $\tau \cdot (H - v) \in {\cal{H}}^{\R}_n$.  Thus it is enough to consider that $p=0$ and $r=1$. Now we assume there is a compactly converging sequence of pointed spaces $(H_i,0)$ with limit  $(H_\infty, 0)$ for $H_i,H_\infty \in {\cal{H}}^{\R}_n$ equipped with, via \D-maps, also converging minimal growth solutions $\Phi_i$ such that $\vartheta[i]:=\vartheta(H_i,0,1,\Phi_i)\ra \infty$.\\
We select regular points $q_i  \in \p  B_{1/2}(0) \subset  H_i$ with $a_n \ge \bp(q_i) \ge b_n >0$, for constants $a_n >b_n>0$,  and $q_i \ra q_\infty \in \p  B_{1/2}(0) \subset H_\infty$, for $i \ra \infty$. This can done due to the naturality of $\bp$ and the compactness of ${\cal{H}}^{\R}_n$.\\
We choose $s_i>0$ so that $s_i \cdot \Phi_i(q_i)=1$. From  elliptic theory the $s_i \cdot \Phi_i$ compactly subconverge to a solution $\Phi_\infty>0$ of $L_{H_\infty,\lambda} \phi=0$ on $H_\infty \setminus \Sigma_{H_\infty}$ with $\Phi_\infty(q_\infty)=1$.
Since $q_\infty$ is a regular point, there is a radius $\eta \in (0,\min \{1,\delta_{\bp}(p_i)\}/4)$, for large $i$ and then independent of $i$, so that $B_{2 \cdot \eta}(q_i)$ is regular and, via $\D$-maps, nearly isometric to $B_{2 \cdot \eta}(q_\infty)$ in $C^3$-norm. Then, we have Harnack inequalities for positive solutions of $L_{H_i,\lambda} \phi=0$ on $B_{2 \cdot \eta}(q_i)$ with constants independent of $i$.  From this $s_i \cdot \Phi_i(q_i)=1$ implies uniform lower estimates $k>0$, for any $i$:
\begin{equation}\label{gcbg}
k \le \int_{B_{\eta}(q_i)}  s_i \cdot \Phi_i \, dV \le \int_{B_2(p_i) \cap H_i \setminus \Sigma} s_i \cdot \Phi_i \, dV.
\end{equation}
For $H_i \in {\cal{H}}^{\R}_n$
we have $\Ric_M, \scal_M, \tr A_H \equiv 0$ and the Bombieri--Giusti inequality (\ref{bgha}) yields
\begin{equation}\label{ab}
\inf_{B_2(0) \cap H_i \setminus \Sigma_i} s_i \cdot \Phi_i \ge  l, \mm{ for a constant } l=l(L_{H_i,\lambda})>0, \mm{ independent of } i.
\end{equation}
As in the proof of Prop.~\ref{lcmg}, this implies a lower positive estimate for the  $d_{\sima}(s_i \cdot \Phi_i )$-distance $d_i$ of  $(\p B_1(0), d_{H_i})$ to $0$: $d_i \ge \Delta$, for some $\Delta>0$ which can be chosen independently of  $i$. \\
To disprove that $\vartheta[i] \cdot \Delta \ra \infty$, for $i \ra \infty$, we consider any  $z \in \p B_1(0)$. We recall from Prop.~\ref{gphi} that due to  $a_n \ge \bp(q_i) \ge b_n >0$ and $s_i \cdot \Phi_i(q_i)=1$ we have  some $\xi>0$ with $\xi=\xi(n)$ so that for the minimal Green's function $G_i$ on $H_i$:
\begin{equation}\label{oba1}
s_i \cdot \Phi_i(z) \le \xi \cdot G_i(z,q_i), \mm{ for any  } z \in \overline{B_1(0)} \cap H_i  \setminus \Sigma_{H_i} \subset \overline{B_2(q_i)} \cap H_i  \setminus \Sigma_{H_i}.
\end{equation}
Moreover, we may assume that $B_{\eta}(q_i)$ is the ball $B=B_{2 \cdot \sigma}(q_i)$ in the hyperbolic picture of Cor.~\ref{glen}  and that from elliptic
 estimates $\Phi_i(z) \le c_\eta \cdot \Phi_i(q_i)$, for some $c_\eta \ge 1$, depending only on $\eta$, and any  $z \in  B_{\eta}(q_i)$.   From this,  (\ref{oba1}) and integrating $G^{2/(n-2)}_i(\cdot,q_i)$  as in Cor.~\ref{glen}  (\ref{upest}), outside $B_{\eta}(q_i)$, along a hyperbolic geodesic arc from $q_i$ to $z$ we get  a common upper bound $\Theta_n^* \cdot \Delta$ for the $\vartheta[i] \cdot \delta_i \ge \vartheta[i] \cdot \Delta$ contradicting the assumption. \\
Finally we reduce  case (ii), that is, $H \in {\cal{G}}^c_n$, to that of  $H \in {\cal{H}}^{\R}_n$. Here we claim that there is a radius $r_H>0$ so that $\Theta_n:=\Theta_n^* +1$ is an upper bound for $\vartheta(H,p,r,\Phi)$, for any $p \in H$, $r \in (0,r_H)$. Otherwise there is a sequence of points $p_i \in H$ and radii $r_i>0$, with $r_i \ra 0$, so that  $\vartheta[i]:=\vartheta(H,p_i,r_i,\Phi)\ra \infty$.  We use the scaling invariance of $L_{H,\lambda} \phi=0$ and scale $g_H$ to $r^{-2} \cdot g_H$ and thereby $(B_r(p),g_H)$ to $(B_1(p), r^{-2} \cdot g_H)$.
Then there is a compactly converging subsequence of the pointed spaces $H_i:=(H, r_i^{-1} \cdot d_H,p_i)$, which contain $(B_1(p_i), r_i^{-1} \cdot d_H,p_i)$, with limit pointed space $(H_\infty, d_{H_\infty},0)$  for $H_\infty \in {\cal{H}}^{\R}_n$. Now we can repeat the argument of case (i), applied to these $H_i$, and see that  $\Theta_n:=\Theta_n^* +1$ upper bounds $\vartheta(H,p,r,\Phi)$, for $r \in (0,r_H)$,  $r_H>0$ small enough. \qed

\subsubsection{Ahlfors $n$-Regularity of $(H,d_{\sima},\mu_{\sima})$} \label{met}

The Hausdorff estimate Th.\ref{haus} suggests a canonical extension of the Riemannian volume measure $\Phi^{2 \cdot n/(n-2)}\cdot \mu_H$ on $H  \setminus \Sigma$ to a measure $\mu_{\sima}$ on $(H,d_{\sima})$, where $\mu_H$ is the $n$-dimensional Hausdorff measure on $(H^n,d_H) \subset (M^{n+1},g_M)$. In turn, $\mu_H$ is the extension of the
Riemannian volume on $(H^n \setminus \Sigma,g_H) \subset (M^{n+1},g_M)$ using that also ${\mathbb{H}}^n(\Sigma)=0$ relative to $(H^n,d_H)$.
\begin{definition}\emph{\textbf{(Minimal Factor Measures $\mu_{\sima}$)}}\label{mms2}\,  For any $H \in {\cal{G}}_n$ equipped with a minimal factor metric $\Phi^{4/(n-2)} \cdot g_H$, we define the minimal factor measure $\mu_{\sima}$ on $H$ by
\begin{equation}\label{meas}
\mu_{\sima}(E):=\int_{E \setminus \Sigma_H} \Phi^{2 \cdot n/(n-2)}\cdot d\mu_H \mm{ for any Borel set } E \subset H.
\end{equation}
\end{definition}

We establish a number of regularity properties of $\mu_{\sima}$, in particular, we will see that this is a locally finite Borel measure on $(H,d_{\sima})$  making  $(H,d_{\sima},\mu_{\sima})$ a \textbf{metric measure space} cf. Rm.\ref{regumu} below.
The local finiteness of the measure follows from the following volume estimates.

\begin{proposition}[Finiteness of $\mu_{\sima}$]\label{vgr0} For any non--totally geodesic  $H \in \cal{G}$ equipped with some minimal factor metric $d_{\sima}(\Phi)$  we have:
\begin{itemize}
  \item for $H \in {\cal{G}}^c_n$, the total volume is finite: $\mu_{\sima}(H) < \infty$,
  \item for $H \in {\cal{H}}^{\R}_n$, any $q\in H$ and any $r>0$: $\mu_{\sima}(B_r(q),d_{\sima}) <\infty$.
\end{itemize}
\end{proposition}

\textbf{Proof} \,  We first show for $H \in {\cal{G}}^c_n$ that $\mu_{\sima}(H,d_{\sima}) < \infty$. This amount to prove that $|\Phi|_{L^{2 \cdot n/(n-2)}} <\infty$. We derive this from the minimal growth condition for $\Phi$ towards $\Sigma$. Note that the optimal estimate for \emph{general} (super)solutions $v>0$ of $L_{H,\lambda}\, \phi= 0$, from the Bombieri--Giusti Harnack inequality,
is $|v|_{L^{p}} <\infty$, $p < n/(n-2)$.\\
Our supersolution $\Phi>0$ properly solves $L_{H,\lambda}\, \phi= 0$ on some small neighborhood $U$ of $\Sigma$. Thus we only need to verify that $\Vol(U \setminus \Sigma,\Phi^{2 \cdot n/(n-2)}\cdot g_H) < \infty$  since $K := (H \setminus \Sigma) \setminus U$ is compact and regular and, hence, $\Vol(K ,\Phi^{2 \cdot n/(n-2)}\cdot g_H) < \infty$. On $U$ we can compare
$\Phi>0$ with the minimal Green's function $G$. We know from the boundary Harnack inequality \ref{mbhsq} that for some base point $p \in K$, with  $\dist_{\bp^*}(p,U) \ge 2 \cdot \sigma$, there is a constant $c\ge 1$ so that
\begin{equation}\label{gpg}
c^{-1} \cdot G(\cdot , p) \le \Phi \le c \cdot G(\cdot , p) \mm{ on } U \setminus \Sigma.\end{equation}

 From \cite[Prop.~3.12, Step 2]{L2}   $G(\cdot , p)$ minimizes the variational integral
\begin{equation}\label{diri}
J_{U}(f):=  \int_{U \setminus \Sigma}   | \nabla_H f |^2 + V_\lambda \cdot  f^2 \, d V, \mm{ for } V_\lambda:= \frac{n-2}{4 (n-1)} \cdot \scal_H - \lambda \cdot \bp^2,
\end{equation}
running over all $f \in H^{1,2}_{\bp}(H \setminus \Sigma)$ with $f|_{\p U} = G(\cdot , p)$ in the trace sense. In particular, using simple test functions we see that $J_{U}(G(\cdot,p)) <\infty$.  From this and the \si-adaptedness of $L_{H,\lambda}$, we have
 for   $C_F:=\int_{\overline{V}} |\nabla F|^2(x)  + V_\lambda \cdot F(x)^2  \,  dV$, for some $C^{2,\alpha}$-extension $F$ of $G(\cdot, p)|_{\p V}$ to $V$,
\small\begin{equation}\label{v1}
\int_{U \setminus \Sigma} \bp^2(x) \cdot G(x, p)^2 \, dV   \le   (\lambda^{\bp}_{L,H}-\lambda)^{-1} \cdot  \int_{U \setminus \Sigma} |\nabla G(\cdot, p)|^2(x)  + V_\lambda \cdot G(x, p)^2  \,  dV + C_F
 <\infty.
\end{equation}
\normalsize
  From (\ref{gre}), (\ref{ges}), (\ref{smot}) and since $\dist_{\bp^*}(p,U) \ge 2 \cdot \sigma$ we get some $\alpha^\circ, \beta^\circ >0$ so that for $x \in U \setminus \Sigma$:
\small \begin{equation}\label{g}
G(x,p) \le \bp^{(n-2)/2}(p)  \cdot  \bp^{(n-2)/2}(x) \cdot \beta^\circ  \cdot \exp(-\alpha^\circ \cdot d_{\bp}(x,p)).
\end{equation}
\normalsize
With this inequality and (\ref{v1}), we get, for some $c ^\circ>0$,
\small\begin{equation}\label{vq1}
 \int_{U \setminus \Sigma} G(x , p)^{4/(n-2) + 2} \,dV \le   c ^\circ \int_{U \setminus \Sigma} \bp^2(x) \,  G(x , p)^2 \,dV  < \infty.
\end{equation}
\normalsize

 For the volume element $dV(g_H)$ of $g_H$, we have  $dV(\Phi^{4/(n-2)} \cdot g_H)=\Phi^{2 \cdot n/(n-2)} \cdot dV(g_H)$. From this, writing $2 \cdot n/(n-2)=4/(n-2) + 2$, we have:
\small\begin{equation}\label{ga}
\Vol(U \setminus \Sigma, \Phi^{4/(n-2)} \cdot g_H) =
\int_{U \setminus \Sigma} \Phi^{2 \cdot n/(n-2)} (x) \,dV \le c^{2 \cdot n/(n-2)} \cdot \int_{U \setminus \Sigma} G(x , p)^{4/(n-2) + 2} \,dV < \infty.
\end{equation}
\normalsize

Thus  for $H \in {\cal{G}}^c_n$ we have $\mu_{\sima}(H,d_{\sima}) < \infty$. For $H \in {\cal{H}}^{\R}_n$ the localization of this argument to balls shows   that $\mu_{\sima}(B_r(q),d_{\sima}) <\infty$. \qed

We refine this finiteness result to volume growth estimates.
To this end we set
\begin{equation}\label{abbr}
\mu_{\sima}(B_r(q))= \Vol(B_r(q),d_{\sima})  = \int_{B_r(q) \cap H \setminus \Sigma} \Phi^{2 \cdot n/(n-2)}\cdot d\mu_H,
\end{equation}
where we write $\Vol(B_r(q),d_{\sima})$ to notationally simplify considerations where we measure distances and volumes with respect to different metrics. These mixed measurements are used to stepwise employ compactness arguments not directly applicable to $(H,d_{\sima},\mu_{\sima})$.

\begin{theorem}[Ahlfors Regularity of Minimal Factors]\label{vgr} For $H \in {\cal{G}}_n$, the metric measure space $(H,d_{\sima},\mu_{\sima})$ is \textbf{Ahlfors $\boldsymbol{n}$-regular}. That is, there are constants $A(H,\Phi),B(H,\Phi)>0$, so that for any $r \in [0,\diam(H,d_{\sima}))$ and any $q\in H$:
\begin{equation}\label{ahl}
 A \cdot r^n \le \mu_{\sima}(B_r(q),d_{\sima}) \le B \cdot r^n.
\end{equation}
For $H \in {\cal{H}}^{\R}_n$ the constants only depend on the dimension, that is, we have $A(n),B(n)>0$. More generally, for any $H \in {\cal{G}}_n$, there is some small $r_{H,\Phi}>0$ so that (\ref{ahl}) holds for these constants $A(n),B(n)>0$ provided $q\in \Sigma_H$ and $r \in (0,r_{H,\Phi})$.
\end{theorem}

\begin{corollary}\emph{\textbf{(Doubling and  Volume Decay)}} \label{dv}  For any $H \in {\cal{G}}_n$, there is a $C(H,\Phi)>0$,  and $C(n)>0$ for $H \in {\cal{H}}^{\R}_n$, so that
for radii and volumina relative to $(H,d_{\sima},\mu_{\sima})$:
\begin{enumerate}
  \item $\mu_{\sima}$ is \textbf{doubling}: for any $q \in H$ and $r \in [0,\diam(H,d_{\sima}))$:
  \begin{equation}\label{dou}
  \mu_{\sima}(B_{2 \cdot r}(q)) \le C \cdot \mu_{\sima}(B_{r}(q)).
  \end{equation}
  \item  For balls $B_* \subset B \subset H$, we have a \textbf{relative lower volume decay} of order $n$. In different terms, the  upper regularity dimension of $\mu_{\sima}$ is at most $n$:
  \begin{equation}\label{volgro}
  \diam(B_*)^n/\diam(B)^n \le C \cdot \mu_{\sima}(B_*)/\mu_{\sima}(B).
  \end{equation}
\end{enumerate}
\end{corollary}

\textbf{Proof of \ref{vgr}} \, We use Prop.~\ref{radii} to treat radius and volume estimates from separate compactness results for the spaces and the eigenfunctions.

\begin{itemize}[leftmargin=*]
\item To this end we introduce a \textbf{radial gauge} of $d_{\sima}(\Phi)$ for $H \in \cal{G}$ in a given $p \in H$:
\begin{equation}\label{incl1}
(B_1(p),d_{\sima}) \subset (B_1(p), d_H) \subset (B_{\kappa_0}(p),d_{\sima}),
\end{equation}
where we replace $\Phi$ for some suitable multiple $k \cdot \Phi$, $k>0$. Under this gauge it is enough to estimate $\Vol\big((B_1(p), d_{H}),d_{\sima}(\Phi)\big)$, that is, the unit ball relative to $d_H$ but with volume measured relative to $d_{\sima}$ for some $\Phi$ satisfying  (\ref{incl}). When $(B_1(p),d_{\sima})$ is given and needs to remain unchanged (e.g.\ this happens in (\ref{voltr}) below), there is an $a>0$ with
\begin{equation}\label{incl2}
(B_1(p), \Phi^{4/(n-2)} \cdot g_H) \subset (B_1(p), a^2\cdot g_H) \subset (B_{\kappa_0}(p), \Phi^{4/(n-2)} \cdot   g_H).
\end{equation}
Then we reinterpret $(B_1(p), \Phi^{4/(n-2)} \cdot g_H)$ as  $(B_1(p),  (a^{-2} \cdot  \Phi^{4/(n-2)}) \cdot a^2\cdot  g_H)$ and
$(B_{\kappa_0}(p), \Phi^{4/(n-2)} \cdot   g_H)$ as $(B_{\kappa_0}(p), (a^{-2} \cdot \Phi^{4/(n-2)}) \cdot  a^2\cdot   g_H)$. This yields a radial gauge not relative to $H$ but to $a \cdot H$ that equally belongs to $\cal{G}$.
\item  For $H \in {\cal{H}}^{\R}_n$ a radial gauge relative to $a \cdot H$, for some $a >0$, suffices to get the volume estimates without gauging since we derive uniform estimates valid in the gauged case for all $H \in {\cal{H}}^{\R}_n$. For $H \in {\cal{G}}^c_n$ we use the Bombieri--Giusti $L^1$--Harnack inequality (\ref{bgha}) to get a positive lower estimate $b>0$ for $\Phi$ on $H \setminus \Sigma$ and this means that we only need to ensure  (\ref{incl2}) for $a \ge b$.
\end{itemize}

We split the proof into five Claims. In Claims 1 and 2 we assume a radial gauge (\ref{incl1}), whereas in Claims 3--5 we drop it and use (\ref{incl2}) and the associated reduction to the gauged case.\\

\noindent \textbf{Claim 1.} \emph{The $d_{\sima}$-volume of the $d_H$--unit ball in $H \in {\cal{H}}^{\R}_n$ satisfies}
\begin{equation}\label{cc}
 a^*_n \le  \Vol\big((B_1(0), d_{H}),d_{\sima}(\Phi)\big), \textit{ for some $a^*_n>0$ depending only on } n.
\end{equation}

\noindent\textbf{Proof of Claim 1.} We assume there were a sequence
\begin{equation}\label{ass}
H_i \in {\cal{H}}^{\R}_n \mm{ with } \Vol\big((B_1(0), d_{H_i}),d_{\sima}(\Phi_i)\big) \ra 0 \mm{ for } i \ra \infty.
\end{equation}
We start with the case $\Phi_i \equiv a_i \cdot G_i(\cdot,p_i)$ for $p_i  \notin (B_1(0), d_{H_i}) \subset H_i$, where $G_i$ is the minimal Green's function on $H_i \setminus \Sigma_{H_i}$,  so that via $\D$-maps $p_i \ra p_\infty$  for some regular point  $p_\infty \notin (B_1(0), d_{H_\infty}) \subset H_\infty$, and where the $a_i>0$ are chosen so that $\Phi_i$ satisfies the gauge  (\ref{incl1}).
We know from \cite[Prop.~3.12]{L2} that the $G_i(\cdot,p_i)$ converge compactly to $G_\infty(\cdot,p_\infty)$, for the minimal Green's function $G_\infty$ on $H_\infty \setminus \Sigma_{H_\infty}$.
From this and the finiteness of $ \Vol\big((B_1(0), d_{H_\infty}),d_{\sima}(G_\infty)\big)$, from \ref{vgr0}, we infer that our assumption (\ref{ass}) implies that $a_i \ra 0$ for $i \ra \infty$.\\
Now we recall from Prop.~\ref{hu} and Prop.~\ref{doob} (\ref{ges}) that the
constants $\alpha, \beta >0$ and $\sigma>0$ for the estimate $G_i^{\sima}(x,y)\le \beta  \cdot \exp(-\alpha \cdot d_{\bp^*}(x,y))$, which we have for $d_{\bp^*}(x,y) \ge 2 \cdot \sigma$, only depend on $n$. From (\ref{upest}) we get a common finite upper bound for the $d_{\sima}$-length of hyperbolic geodesic rays from $p_i$ to $0$. Since $a_i \ra 0$ for $i \ra \infty$, this shows that $(B_1(0),d_{\sima}(\Phi_i)) \nsubseteq (B_1(0), d_{H_i})$ for $i$ large enough, contradicting the chosen radial gauge.\\
This argument extends to more general $\Phi_i$ satisfying (\ref{ass}) since we know that under scalings and blow-ups of a given $H \in \cal{G}$, and similarly on $ {\cal{H}}^{\R}_n$, the boundary Harnack inequality \ref{mbhsq} applies with a \emph{common} Harnack constant that is independent of $i$, cf.\ Remark \ref{bhp3}.2. \qed

\noindent \textbf{Claim 2.}  \emph{The volume of the $d_{\sima}$-unit ball in $H \in {\cal{H}}^{\R}_n$ satisfies
\begin{equation}\label{cl2}
\Vol\big((B_1(0), d_{H}),d_{\sima}(\Phi)\big) \le b^*_n, \mm{ for some } b^*_n>0 \mm{ depending only on } n.
\end{equation} }
\noindent \textbf{Proof of Claim 2.} This time we assume there were a sequence
\begin{equation}\label{ass2}
H_i \in {\cal{H}}^{\R}_n \mm{ with } \Vol\big((B_1(0), d_{H_i}),d_{\sima}(\Phi_i)\big) \ra \infty, \mm{ for } i \ra \infty.
\end{equation}
Again, it is enough to consider the case $\Phi_i \equiv a_i \cdot G_i(\cdot,p_i)$, for $p_i  \notin (B_1(0), d_{H_i}) \subset H_i$, where $G_i$ is the minimal Green's function on $H_i \setminus \Sigma_{H_i}$  so that via $\D$-maps $p_i \ra p_\infty$  for some regular point  $p_\infty \notin (B_1(0), d_{H_\infty}) \subset H_\infty$ and where the $a_i>0$ are chosen so that $\Phi_i$ satisfies the  gauge  (\ref{incl1}). We use again that the $G_i(\cdot,p_i)$ converge compactly to $G_\infty(\cdot,p_\infty)$. This time we additionally use, from the proof of \ref{vgr0},  that $J_{H_i \setminus (B_{\sigma}(p_i),d_{\bp^*})}(G_i)$ upper bounds $\Vol(B_1(0), G_i^{4/(n-2)} \cdot g_H)$. To upper bound $J_{H_i \setminus (B_{\sigma}(p_i),d_{\bp^*})} (G_i)$ we can choose cut-off functions $\phi_i \ge 0$ with $\phi_i = G_i$ on $\p B_{\sigma}(p_i)$ and $\phi_i \equiv 0$ outside $B_{2 \cdot \sigma}(p_i)$ so that  $J_{H_i \setminus (B_{\sigma}(p_i),d_{\bp^*})}(\phi_i) \le c$, where $c>0$ does not depend on $i$.   \\
Then  (\ref{ass}) implies that $a_i \ra \infty$, for $i \ra \infty$ and we now show that this contradicts the second inclusion of the radial gauge (\ref{incl1}). To this end we note  that the
compact $\D$-map convergence $G_i(\cdot,p_i) \ra G_\infty(\cdot,p_\infty)$ shows that the $L^1$-norm of $G_i(\cdot,p_i)$ on $(B_1(0), d_{H_i})$ remains positively lower bounded.  The Bombieri--Giusti $L^1$--Harnack inequality (\ref{bgha}) and the argument of Prop.~\ref{lcmg}\,(ii)  therefore show that $ (B_1(0), d_{H_i})\nsubseteq (B_{\kappa_0}(0),d_{\sima})$ for $i$ large enough.\qed

\noindent \textbf{Claim 3.}   \emph{For constants $a_n^*,b_n^* >0$ depending only on $n$, we have for any  $H \in {\cal{H}}^{\R}_n$:
\begin{equation}\label{cl3}
a_n^* \cdot r^n \le  \Vol(B_r(q),d_{\sima}) \le b_n^* \cdot r^n, \mm{ for any } r>0.
\end{equation}}
\noindent \textbf{Proof of Claim 3.}
To determine the volume growth rate of balls in $H$ in terms of $r$, we use (\ref{incl2}), as explained above, to see that the volume estimates from claims 1 and 2 hold also without the radial gauge. That is, we have for any $H \in {\cal{H}}^{\R}_n$:
\begin{equation}\label{imma}
a^*_n \le \Vol\big((B_1(0), d_{H}),d_{\sima}(\Phi)\big) \le b^*_n.
\end{equation}
Similarly, for $H \in {\cal{H}}^{\R}_n$ we also have $r^{-1} \cdot H \in {\cal{H}}^{\R}_n$ and we apply the unit ball estimate to $r^{-1} \cdot H \in {\cal{H}}^{\R}_n$ and then rescale $r^{-1} \cdot H$ to $H$. Then the identity
\begin{equation}\label{voltr}
\Vol(B_r(0), \Phi^{4/(n-2)} \cdot g_H) = r^n \cdot \Vol(B_1(0), \Phi^{4/(n-2)} \cdot  r^{-2} \cdot g_H),
\end{equation}
shows:\, $a^*_n \cdot r^n \le  \Vol(B_r(0), \Phi^{4/(n-2)} \cdot g_H)  \le b^*_n \cdot r^n.$\qed

\noindent \textbf{Claim 4.}   \emph{For $H \in {\cal{G}}^c_n$ and $r \in (0,r_H)$, for a suitably small $r_{H,\Phi}>0$, we have for $q \in \Sigma_H$:
\[a_n \cdot r^n \le  \Vol(B_r(q),d_{\sima}) \le b_n \cdot r^n,\mm{ for }a_n,b_n >0 \mm{ depending only on } n.\]}
\noindent \textbf{Proof of Claim 4.} We first consider one fixed $H \in {\cal{G}}^c_n$ and assume a radial gauge. Then there are constants $b>a>0$ and some small $r>0$ so that: $a \le  \Vol(B_1(q), \Phi^{4/(n-2)} \cdot  r^{-2} \cdot g_H) \le b$, for any $q \in \Sigma_H$ and  $r \in (0,r_H)$. Otherwise we had a converging sequence of points $q_i \in \Sigma_H$ and of radii $r_i \ra 0$, for $i \ra \infty$ so that these volumina would either converge to $0$ or $\infty$. Both cases can be ruled out as in Claim 1 and 2 for  $H \in {\cal{H}}^{\R}_n$ above.\\
Moreover, the constants $b>a>0$  can be chosen to depend only on the dimension. Otherwise we had  a compactly converging sequence $r_i^{-1} \cdot H_i \in {\cal{G}}^c_n$ with $r_i \ra 0$, for $i \ra \infty$, and of unit balls $(B_1(q_i), \Phi^{4/(n-2)} \cdot  r_i^{-2} \cdot g_{H_i})$ satisfying a radial gauge, so that, again, the associated volumina converge either to $0$ or to $\infty$. Both cases can be excluded as before. \\
Finally, we recall that $\Phi \ge \overline{b}>0$ for some constant $\overline{b}>0$ from the  $L^1$--Harnack inequality (\ref{bgha}). Then  (\ref{voltr}) in the argument for Claim 3 applies to all $r \in (0,r_{H,\Phi})$, possibly after replacing $r_{H,\Phi}$ by $b \cdot r_{H,\Phi}$, since we only need to use $\overline{a} \ge \overline{b}/r$ in (\ref{incl2}) to find a radial gauge for any
$(B_1(q_i), \Phi^{4/(n-2)} \cdot  r_i^{-2} \cdot g_{H_i})$. This shows that there are $\overline{a}_n, \overline{b}_n >0$ depending only on $n$ so that
$\overline{a}_n \cdot r^n \le  \Vol(B_r(q),d_{\sima}) \le \overline{b}_n \cdot r^n$ for any $H \in {\cal{G}}^c_n$ and $r \in (0,r_{H,\Phi})$. \qed

\noindent \textbf{Claim 5.}   \emph{For any $H \in {\cal{G}}^c_n$ there are constants $a(H,\Phi),b(H,\Phi)>0$ so that for any $r \in [0,\diam(H,d_{\sima})]$ and any $q\in H$:
$ a \cdot r^n \le \mu_{\sima}(B_r(q),d_{\sima}) \le b \cdot r^n.$}\\

\noindent \textbf{Proof of Claim 5.}  This readily follows from claim 4. For $I:= [r_H,\diam(H,d_{\sima})]$ we define
$a:= \min \{ a_n , \inf \{ \Vol(B_{r}(q),d_{\sima})/r^n \, | \, r \in I\} \}$,  $b:= \max \{ b_n , \sup \{ \Vol(B_{r}(q),d_{\sima})/r^n \, | \, r \in I \} \}.$
 \qed

\begin{remark}[Regularity of  $\mu_{\sima}$]\label{regumu}  As a consequence of the estimates in Prop.~\ref{vgr} and \ref{vgr0}, we see that $\mu_{\sima}$ is an \textbf{outer regular measure}, that is, we have for Borel subsets $E \subset H$:
\begin{equation}\label{outme}
\mu_{\sima}(E)= \inf \{\mu_{\sima}(A) \,|\, E \subset A , A \subset H \mm{ open} \}.
\end{equation}
From Prop.~\ref{vgr0}  we  know that $\mu_{\sima}(H) < \infty$ for $H \in {\cal{G}}^c_n$. Thus for any $\ve >0$ there is a neighborhood $U_\ve$ of $\Sigma$ in $H$ so that $\mu_{\sima}(U_\ve \setminus \Sigma) < \ve$. This also holds for non-compact $H \in {\cal{G}}_n$ from $\mu_{\sima}(B_r(q)) < \infty$ using suitable ball covers  of $\Sigma$.
From this we see that $\mu_{\sima}$ is a \textbf{Borel measure} on $(H,d_{\sima})$ cf.\cite[pp.~62--64]{H-T}. \qed
\end{remark}

\begin{remark}[Ahlfors Regularity of $H \in \cal{G}$] We note in passing that the original almost minimizers $H \in {\cal{G}}_n$ are Ahlfors $n$-regular as well.   For  $H \in {\cal{H}}^{\R}_n$ this is \cite[Prop.~5.14 and Rm.~5.15]{Gi} and there are constants $w_n > v_n>0$ depending only on $n$ so that  $v_n \le  \Vol\big((B_1(0), d_H),d_H\big) \le w_n \mm{ for any } H \in {\cal{H}}^{\R}_n.$
For $H \in {\cal{G}}^c_n$ the Ahlfors $n$-regularity follows similarly from the almost optimal isoperimetric inequality $|\Vol\big((B_r(p), d_H),d_H\big)-\Vol\big(P_r\big)|\le K \cdot r^{n+2 \cdot \alpha}$, for some $\alpha \in (0,1)$, $K>0$, and sufficiently small $r>0$, where $P_r$ is an area minimizing Plateau solution with boundary data $H$ along $\p B_r(p)$ and $c_n$ is the Euclidean volume of the unit ball. \qed
\end{remark}

\subsubsection{Poincar\'{e} and Sobolev Inequalities} \label{sobo}

We start with some estimates for the minimal Green's function on the twisted \si-double cones we introduced in Prop.~\ref{sem0}. Henceforth we use a fixed size parameter $\boldsymbol{d \le \min\{1,c\}/2}$ from \ref{sem0}, Step 2 and recall that $d$ is independent of the chosen $p,q \in H$, that is, it depends only on $H \in {\cal{G}}$ and for $H \in {\cal{H}}^{\R}_n$ it only depends on $n$.\\
We consider the canonical Semmes families $\Gamma_{p,q}$ and their envelopes $\textbf{E}_{p,q}$ for some hyperbolic geodesic $\gamma_{p,q}$. Due to the symmetry of the definition in $p$ versus $q$ it will be enough to analyze the part/side closer to $p$ than to $q$. To make this precise, we introduce some terminology.
For $\gamma \in \Gamma_{p,q}$ we consider those $z \in \gamma$ closer to $p$ than to $q$, measured in $\gamma$--arc length and define:
\begin{equation}\label{tuu}
\gamma^+[z]:= \mm{subcurve of } \gamma \in \Gamma_{p,q} \mm{ from } p  \mm{ to }  z, \mm{ and we have } l_H(\gamma^+[z]) = l_{min}(\gamma_{p,q}(z)).
\end{equation}
For the midpoint $m_\gamma$ of each such $\gamma$ we define the \textbf{half-curve family} and the \textbf{half-envelope}
\begin{equation}\label{plu}
\textstyle \Gamma^+_{p,q}[d]:=\{\gamma[m_\gamma]\,|\, \gamma \in \Gamma_{p,q}[d]\} \mm{ and } \textbf{E}^+_{p,q}[d]  := \bigcup_{z \in \gamma_{p,m} \setminus \{p\}} B_{d \cdot l_{min}(\gamma_{p,q}(z))/c}(z).
\end{equation}
The counterparts starting from $q$ are denoted by $\Gamma^-_{p,q}[d]$ and $\textbf{E}^-_{p,q}[d]$. We have $ \textbf{E}^+_{p,q}[d] \cap \textbf{E}^-_{p,q}[d] = B_{d \cdot l(\gamma_{p,q}(m))/c}(m)$, i.e. in  some arguments the points in $B_{d \cdot l(\gamma_{p,q}(m))/c}(m)$ are counted twice. \\ The particular choice of our core curve $\gamma_{p,q}$  gives us a good control over the analysis on the $ \textbf{E}^+_{p,q}[d]$.
We start with a variant of the Harnack inequality that controls the supersolutions $\Phi>0$ of Def.~\ref{msge} on $ \textbf{E}^+_{p,q}[d]$ transversally to the core geodesic $\gamma_{p,q}$.

 \begin{lemma}[Transversal Harnack Inequalities]\label{hat} There are constants $C(H,\Phi)>0$, and $C(n)>0$ for $H \in {\cal{H}}^{\R}_n$,  so that for $B_{d \cdot l_H(\gamma^+[z])/c}(z) \subset \emph{\textbf{E}}^+_{p,q}[d]$ and any $p,q \in H$:
 \begin{equation}\label{har}
\textstyle \sup_{B_{d \cdot l_H(\gamma^+[z])/c}(z)} \Phi \le C \cdot \inf_{B_{d \cdot l_H(\gamma^+[z])/c}(z)} \Phi.
\end{equation}
 \end{lemma}

\textbf{Proof} \, We first prove this for general positive solutions $u>0$ of $L_{H,\lambda} \, \phi=0$ on  $\textbf{E}_{p,q}[1] \subset H \setminus \Sigma$.  We
can scale any of these ball $B_{d \cdot l_H(\gamma^+[z])/c}(z)$ to unit size where the underlying geometry becomes uniformly bounded in $C^3$-norm, cf.~(\ref{pull}) in Prop.~\ref{sem0}, Step 2,  independently of $z$ and also of $H$.  After scaling around any such $z$, the exponential map pull-backs of $L_{H,\lambda}= -\Delta  +\frac{n-2}{4 (n-1)} \cdot \scal_H- \lambda \cdot \bp^2$ to the Euclidean unit ball in the tangent space have uniformly bounded coefficients independent of $z$. From this the Harnack inequality holds for positive solutions of the pull-back equations with the same Harnack constant on any of these unit balls. This relation is invariant under scalings and it survives the exponential map  transfer and rescaling back to $B_{d \cdot l_H(\gamma^+[z])/c}(z)$.\\
For $H \in {\cal{H}}^{\R}_n$ this gives  (\ref{har}) for a constant that  merely depends on $n$. When $H \in {\cal{G}}^c$  there is a compact set $K \subset H  \setminus \Sigma$ so that  $\Phi>0$ is a solution on $H  \setminus (\Sigma \cup K)$ where we can apply the argument to balls in $\textbf{E}_{p,q}[d]$ disjoint from $K$.
All balls with a non-empty intersection with $K$ belong to another still compact subset $K^* \subset H  \setminus \Sigma$ where we find a constant satisfying (\ref{har}) right from the continuity of $\Phi>0$ on $H  \setminus \Sigma$.
\qed

 \begin{corollary}[Path Integral Estimate]\label{nochma} There is a constant $k(H,\Phi)>0$, with $k(n)>0$ for $H \in {\cal{H}}^{\R}_n$, so that for any $p, q \in H$:
\begin{equation}\label{ib}
l_{\sima}(\gamma^+[z]) = \int_{\gamma^+[z]}\Phi^{2/(n-2)} ds \le  k \cdot l_{H}(\gamma^+[z]) \cdot \Phi^{2/(n-2)}(z),  \mm{ for any } \gamma \in \Gamma^+_{p,q}[d].
\end{equation}
\end{corollary}

\textbf{Proof} \, By Lemma \ref{hat} it is enough to consider the subcurves $\gamma^+[z]$ of the core $\gamma_{p,q}$. Since (\ref{ib}) is scaling invariant, we may assume that $l_{H}(\gamma^+[z])=1$. We may multiply the inequality by a constant so that $\Phi(z)=1$. The $c$-\si-uniformity of $(H,d_H)$ implies that $\bp(z) \le c$ and that  $d_H(p,z) \le l_{H}(\gamma^+[z]) \le c \cdot d_H(p,z)$ and we can apply Prop.~\ref{gphi} to get a constant $\xi>0$ with
$\Phi(x) \le \xi \cdot G(x,z), \mm{ for any  } x, z \in H \setminus \Sigma.$
From (\ref{upest}) we know that the contributions of integrals $\int_{\gamma^+[z]}G(\cdot,z)^{2/(n-2)} ds$ outside $B_{\delta_{\bp}(z)}(z)$ are uniformly upper bounded. We have $\Phi(x) \le c^*$ for $x \in B_{\delta_{\bp}(z)}(z)$, for some $c^*(H,\Phi)>0$. Namely, for $H \in {\cal{G}}^c_n$ this follows for any such ball that intersects the  compact set $K$ where $\Phi$ is only a supersolution. But for these balls such an upper estimate follows from the continuity of $\Phi$. For balls away from $K$, and for $H \in {\cal{H}}^{\R}_n$, we get the bound from elliptic estimates starting from $\Phi(z)=1$ since, in the hyperbolic picture, the balls have bounded geometry,  a fixed radius and the operator has uniformly bounded coefficients. For $H \in {\cal{H}}^{\R}_n$ this also shows that $c^*(n)>0$. Thus we get some $k^*>0$ with the asserted dependencies and $\int_{\gamma^+[z]}\Phi^{2/(n-2)} ds \le  k^*= l_{H}(\gamma^+[z]) \cdot \Phi^{2/(n-2)}(z)$.   \qed

Now we show that the canonical Semmes families of curves $\Gamma_{p,q}[d]$ with the probability measure $\sigma_{p,q}[d]$ on $(H,d_H,\mu_H)$  we have defined in Prop.~\ref{sem0} are still  Semmes families in $(H,d_{\sima},\mu_{\sima})$. This and the volume relations \ref{dv} imply Poincar\'{e}, Sobolev and isoperimetric inequalities for $(H,d_{\sima},\mu_{\sima})$.

\begin{theorem} [Semmes Families on $(H,d_{\sima},\mu_{\sima})$]\label{semin} For $H \in {\cal{G}}_n$, there is some constant $C_{\sima}(H,\Phi)>0$,   $C_{\sima}(n)>0$ for $H \in {\cal{H}}^{\R}_n$,
so that for   $p,q \in H$, the family $\Gamma_{p,q}$ and the probability measure $\sigma_{p,q}$ on  $\Gamma_{p,q}$, from \ref{sem0}, satisfy the two Semmes axioms relative to $(H,d_{\sima},\mu_{\sima})$:
\begin{enumerate}
 \item For any $\gamma \in \Gamma_{p,q}$: $l_{\sima}(\gamma|_{[s,t]}) < C_{\sima}\cdot d_{\sima}(\gamma(s),\gamma(t))$, for  $s,t \in I_\gamma$.
 \item For any Borel set $A \subset X$, the assignment $\gamma \mapsto l_{\sima}(\gamma \cap A)$  is $\sigma$-measurable with
 \small \begin{equation}\label{tcu}
 \int_{\Gamma_{p,q}} l_{\sima}(\gamma \cap A) \, d \sigma(\gamma) \le  C_{\sima} \cdot \int_{A_{C_{\sima},p,q}} \left(\frac{d_{\sima}(p,z)}{\mu_{\sima}(B_{d_{\sima}(p,z)}(p))} + \frac{d_{\sima}(q,z)}{\mu_{\sima}(B_{d_{\sima}(q,z)}(q))}\right) d \mu_{\sima}(z)
 \end{equation}
 \normalsize
 for $A_{C_{\sima},p,q}:=(B_{C_{\sima} \cdot d_{\sima}(p,q)}(p) \cup B_{C_{\sima} \cdot d_{\sima}(p,q)}(q))\cap A$.
\end{enumerate}
 \end{theorem}
\vspace{0.1cm}
\textbf{Proof of Property (i)}\, In the case of  $H \in {\cal{H}}^{\R}_n$ we first prove that  there is a constant $c^*(n)>0$ so that $l_{\sima}(\gamma_{x,y}) < c^*(n)\cdot d_{\sima}(x,y)$ for an arbitrary hyperbolic geodesic arc $\gamma_{x,y} \subset H$ linking two points $x,y \in H$ with $d_H(x,y)=1$.
For this we choose the midpoint $m \in \gamma_{x,y}$ of this  $c$-\si-uniform curve in $(H,d_H)$, measured in terms of curve  length relative to $d_H$. Now we use the gauge $\Phi(m)=1$ and apply (\ref{ib}) to the two subcurves of $\gamma_{x,y}$ starting from $m$. Since $\gamma_{x,y}$ is $c$-\si-uniform we have a length estimate
$l_{\sima}(\gamma_{x,y}) \le l_n$, for some $l_n>0$ depending only on $n$. The $c$-\si-uniformity also shows that there is a ball  $B_{r_n}(m)$ of radius $r_n>0$ where we have $\bp \le a_n$ and thus we find a uniform Harnack estimate $b_n>0$ so that $\Phi(m) \ge b_n$. The Bombieri--Giusti Harnack inequality then gives, as in (\ref{gcbg}) and (\ref{ab}),  a lower estimate $e_n>0$ for $\Phi$ on $B_{2 \cdot c}(m)$ and from this we infer a lower estimate $e^*_n>0$ for $d_{\sima}(x,y)$, that is, we have $l_{\sima}(\gamma_{x,y}) \le l_n/e^*_n \cdot d_{\sima}(x,y)$.\\
In particular, this applies to any subcurve of the core geodesic we have in any of our families $\Gamma_{p,q}[d]$. In turn, the transversal Harnack inequality (\ref{har}) yields a constant $C_0(n) \ge  l_n/e^*_n>0$ so that for any other $\gamma \in \Gamma_{p,q}[d]$ we have $l_{\sima}(\gamma|_{[s,t]}) < C_0\cdot d_{\sima}(\gamma(s),\gamma(t))$, for  $s,t \in I_\gamma$.\\
For $H \in {\cal{G}}^c_n$ the subset where $\Phi$ is not a proper solution belongs to some compact  $K \subset H \setminus \Sigma$. Since any hyperbolic geodesic arc $\gamma_{x,y} \subset H$ is a $c$-\si-uniform curve in $(H,d_H)$, we have $l_H(\gamma_{x,y}) \le c \cdot d_H(x,y)$. The contributions on $K$ of the conformal deformation by the upper and lower positively bounded function $\Phi$ merely alter $c$ to another constant that depends on the chosen $H$ and $\Phi$. We combine this with the argument for case (i) outside $K$ to infer the claim for $H \in {\cal{G}}^c_n$. \qed

\textbf{Proof of Property (ii)}\, By Prop.~\ref{vgr},  Lemma \ref{nochma} and since both sides of (\ref{tcu}) result from smooth deformations of $(H,d_H)$, they are still finite Borel measures:

\small
 \[\mu_1(A):= \int_{\Gamma_{p,q}} l_{\sima}(\gamma \cap A) \, d \sigma(\gamma)\, \mm{  and  } \, \mu_2(A):= \int_{A_{C,p,q}} \left(\frac{d_{\sima}(p,z)}{\mu_{\sima}(B_{d_{\sima}(p,z)}(p))} + \frac{d_{\sima}(q,z)}{\mu_{\sima}(B_{d_{\sima}(q,z)}(q))}\right) d \mu_{\sima}(z).\]
\normalsize
 To derive inequality (\ref{tcu})  we make a series of simplifications.
\begin{itemize}[leftmargin=*]
\item We have $\overline{\textbf{E}_d(p,q)} \cap \Sigma \subset \{p,q\}$. Thus to check  (\ref{tcu}) we only need to consider Borel sets $A \subset H \setminus \Sigma$.
  \item From the $\sigma$-additivity and the regularity of the $\mu_i$ we only need to show  $\mu_1(B) \le C \cdot \mu_2(B)$ for arbitrarily small balls $B =B_{\ve(x)}(x)$ for some $\ve(x)>0$ with $B_{3 \cdot \ve(x)}(x)\subset H \setminus \Sigma$ for any $x \in H \setminus \Sigma$, for some $C>0$  independent of $B$ and of $x$.
\item    We may assume that for any $x \in H \setminus \Sigma$  the ball $B=B_{\ve(x)}(x)$ is small enough so that
\begin{equation}\label{ggg}
1/2 \cdot \Phi^{\frac{2}{n-2}} (x) \le \Phi^{\frac{2}{n-2}} (y) \le 2 \cdot  \Phi^{\frac{2}{n-2}} (x),\mm{ for } y \in B.
\end{equation}
\end{itemize}
Then we have for $A:= B$:
\small
\begin{equation}\label{tcu0}
1/2 \cdot \int_{\Gamma_{p,q}}   l_{\sima}(\gamma \cap A)  \, d \sigma(\gamma) \le  \int_{\Gamma_{p,q}}  \Phi^{\frac{2}{n-2}} (x) \cdot l_{H}(\gamma \cap A)  \, d \sigma(\gamma) \le ...
 \end{equation}
\normalsize
Since $\Gamma_{p,q}$ is a Semmes family relative to $(H,d_H,\mu_H)$, we have from Prop.~\ref{sem0}:
  \small   \begin{equation}\label{tcu1a}
\Phi^{\frac{2}{n-2}} (x) \cdot \int_{\Gamma_{p,q}}  l_{H}(\gamma \cap A)  \, d \sigma(\gamma)\le \Phi^{\frac{2}{n-2}} (x) \cdot  C \cdot \int_{A_{C,p,q}}   \frac{d(q,z)}{\mu_H(B_{d(q,z)}(q))} d \mu_H(z) \le ...
 \end{equation}
 \normalsize
Using (\ref{ggg}) again and $d(\gamma(s),\gamma(t)) \le l_H(\gamma|_{[s,t]}) < c \cdot d(\gamma(s),\gamma(t))$, for any $\gamma \in \Gamma_{p,q}$, $s,t \in I_\gamma$, and also that $a \cdot r^n \le  \mu_H(B_r(q)) \le b \cdot r^n$,  we get a constant $C_1>0$, depending only on $(H,\Phi)$, or depending only on $n$ for $H \in {\cal{H}}^{\R}_n$, with
\small \begin{equation}\label{tcu2}
2 \cdot C \cdot \int_{A_{C,p,q}} \Phi^{\frac{2}{n-2}}(z)  \cdot  \frac{d(q,z)}{\mu(B_{d(q,z)}(q))}  d \mu_{H}(z) \le   C_1 \cdot \int_{A_{C,p,q}} \Phi^{\frac{2}{n-2}}(z)  \cdot  \frac{l_H(\gamma_q[z]) }{l_H(\gamma_q[z])^n}  d \mu_{H}(z) \le ...
 \end{equation}
\normalsize
Now we apply Lemma~\ref{nochma}, $ l_{H}(\gamma_q[z])^{-1} \le  k \cdot \Phi^{\frac{2}{n-2}}(z)/l_{\sima}(\gamma_q[z])$ for any $z \in H \setminus \Sigma \cap B_\rho(q)$:
  \small  \begin{equation}\label{d2g}
C_1  \cdot k^{n-1}  \cdot \int_{A_{C,p,q}} \Phi^{\frac{2}{n-2}}(z)  \cdot \left(\Phi^{\frac{2}{n-2}}(z) \big/l_{\sima}(\gamma_q[z])\right)^{n-1} d \mu_{H}(z) = ...
  \end{equation}
  \normalsize
 Finally, we get from $d_{\sima}(\gamma(s),\gamma(t)) \le l_{\sima}(\gamma|_{[s,t]}) < C_0 \cdot d_{\sima}(\gamma(s),\gamma(t))$ for any $\gamma \in \Gamma_{p,q}$, $s,t \in I_\gamma$, the Ahlfors regularity of $(H,d_{\sima},\mu_{\sima})$ saying that $a \cdot r^n \le  \mu_\sima(B_r(q),d_{\sima}) \le b \cdot r^n$ (Theorem~\ref{vgr}) and the eccentricity estimate of Prop.~\ref{radii} some new constant $C_2>0$, with the same dependencies, so that:
   \small \begin{equation}\label{tcu3}
 C_1 \cdot k^{n-1}  \cdot  \int_{A_{C,p,q}} \frac{l_{\sima}(\gamma_q[z])}{l_{\sima}(\gamma_q[z])^n} \cdot \Phi^{2 n/(n-2)}(z)  \cdot d \mu_{H}(z)  \le  C_2 \cdot \int_{A_{C_2,p,q}} \frac{d_{\sima}(q,z)}{\mu_{\sima}(B_{d_{\sima}(q,z)}(q))}d \mu_{\sima}(z).
 \end{equation}
 \normalsize
Thus we can choose $C_{\sima}$ to be the maximum of $C_0$ and $C_2$. \qed

A standard application of Prop.~\ref{semin} and the doubling property is the following weak $(1,1)$--Poincar\'{e} inequality  \cite[Ch.14.2]{H-T}. To state it, we recall that a measurable function $w \ge 0$ on  $(H , d_{\sima})$ is an \emph{upper gradient} of a measurable function $u$ if  $|u(x)- u(y)| \le  \int_c w(s) ds$ holds for all rectifiable curves $c$ joining $x$ to $y$, for any pair $x,y \in H$.

\begin{proposition} \emph{\textbf{(Poincar\'{e} Inequality I)}}\label{poiic} For any $H \in {\cal{G}}$, there are $C_0(H,\Phi) >0$, $\gamma_0(H,\Phi) \ge 1$, depending only on $n$ for $H \in {\cal{H}}^{\R}_n$, so that for any pair of concentric balls $B \subset  \gamma_0 \cdot B\subset (H,d_\sima)$, for any function $u$ on $H$, integrable on bounded balls, and any upper gradient $w$ of $u$ we get:
\begin{equation}\label{poinm0}
\fint_B |u-u_B| \,  d \mu_{\sima} \le C_0 \cdot \diam(B) \cdot \fint_{\gamma_0 \cdot B} w \, d \mu_{\sima}, \mm{ for } f_B:=\fint_B f  \,  d \mu_{\sima} := \int_B f \, d \mu_{\sima}/\mu_{\sima}(B).
\end{equation}
\end{proposition}

The volume decay property \emph{of order n} in Cor.~\ref{dv}(ii) allows us to improve this Poincar\'{e} inequality to the following Sobolev inequality.

\begin{corollary} \emph{\textbf{(Sobolev Inequality)}}\label{sobbc} For any $H \in {\cal{G}}$, there is a constant $C_1(H,\Phi)>0$, depending only on $n$ for $H \in {\cal{H}}^{\R}_n$,  so that for some open ball $B \subset H$, an $L^1$-function $u$ on $B$ and any upper gradient $w$ of $u$ on $B$, we have
\begin{equation}\label{ii2}
\Big(\fint_B |u-u_B|^{n/(n-1)} \,  d \mu_{\sima}\Big)^{(n-1)/n}  \le C_1 \cdot \diam(B) \cdot \fint_B  w \, d \mu_{\sima}.
\end{equation}
\end{corollary}

\textbf{Proof of \ref{poiic} and \ref{sobbc}} \,  Cor.~\ref{dv} and Prop.~\ref{semin}  imply the underlying Poincar\'{e} inequality \ref{poiic}, see \cite[14.2, p.~396]{H-T}. From this we get the refinement to the Sobolev inequality  \ref{sobbc}   from  \cite[Th.~9.1.15(i)]{H-T}, for $p=1$ and $Q=n$, see also \cite[Th.~4.5 and Rm.~4.6]{M} and {\cite{Se}. One first uses cut-off functions to restrict the support to $H \setminus \Sigma$ and then one extends the inequalities to $H$  using that the Hausdorff dimension of $\Sigma \subset  (H,d_{\sima})$ is $\le n-3$.\qed

In turn, Cor.\ref{sobbc} and the doubling property of $(H,d_{\sima},\mu_{\sima})$ show, cf.~\cite[Remark 9.1.19]{H-T}, that Prop.\ref{poiic} can be improved so that we can drop the scaling factor $\gamma_0$.
\begin{corollary} \emph{\textbf{(Poincar\'{e} Inequality II)}}\label{poiic2}  For any $H \in {\cal{G}}$, there are $C_1>0$, with the same dependencies as in  \ref{poiic}, so that  for any function $u$ on $H$, integrable on bounded balls, and any upper gradient $w$ of $u$, we get
\begin{equation}\label{poinm}
\int_B |u-u_B| \,  d \mu_{\sima} \le C_1 \cdot \diam(B) \cdot \int_{B} w \, d \mu_{\sima}, \mm{ for any ball } B\subset (H,d_\sima).
\end{equation}
\end{corollary}

\subsubsection{Oriented Minimal Boundaries} \label{omb}

We extend the Riemannian hypersurface area i.e. the element $d\mu^{n-1}_{\sima}=\Phi^{2 \cdot n/(n-2)}\cdot d \mu^{n-1}_H$, on $H \setminus \Sigma$ to $(H,d_{\sima},\mu_{\sima})$ to formulate the isoperimetric inequality and the concept of oriented minimal boundaries.  To this end we employ the BV (= bounded variations) approach of Ambrosio \cite{A} and Miranda \cite{M} on complete metric spaces with a doubling measure supporting a Poincar\'{e} inequality. By the results in the last two sections this theory applies to $(H,d_{\sima},\mu_{\sima})$. We reformulate \cite[Definition 4.1]{M} as follows:

\begin{definition}\emph{\textbf{(Perimeters in $(H,d_{\sima},\mu_{\sima})$)}}\label{mmsh}\,  For some Borel set $E$ and an open set $\Omega$ in $(H,d_{\sima},\mu_{\sima})$, for $H \in {\cal{G}}_n$, we define the \textbf{perimeter} $\mu^{n-1}_{\sima}(\p E \cap \Omega)$, as
\begin{equation}\label{du}
\inf \big\{\liminf_{k \ra \infty} \int_\Omega|\nabla u_k| \, d\mu_{\sima} \,\Big|\,  L^1_{\mathrm{loc}}(\Omega)\mm{-converging sequences } u_k \ra \chi_E, u_k \in \mathrm{Lip_{loc}}(\Omega)\big\},
\end{equation}
where $\chi_E$ is the characteristic function of $E$ and, for $u \in \mathrm{Lip_{loc}}(\Omega)$, we use the particular upper gradient  $|\nabla u|(x) :=\liminf_{\varrho \ra 0} \sup_{y \in\overline{ B _{\varrho}(x)}}  |u(x)-u(y)|/\varrho$ cf.\cite[p.982]{M}. We call $E$ a  \textbf{Caccioppoli set} in $(H,d_{\sima},\mu_{\sima})$ provided $\mu^{n-1}_{\sima}(\p E \cap \Omega) < \infty$, for any bounded $\Omega$.
\end{definition}

For a smoothly bounded open set $E \subset \R^n$ the perimeter is the hypersurface area of the boundary $\p E$ in $\Omega$  that equals its $(n-1)$-dimensional Hausdorff measure, see \cite[Example 1.4]{Gi}. In general, only the expression $\mu^{n-1}_{\sima}(\p E \cap \Omega)$ is relevant and well-defined.\\
 The perimeter (\ref{du}) satisfies the coarea formula \cite[Prop.~4.2]{M} and thus we find many non-trivial Caccioppoli sets e.g. \cite[Cor.~4.4]{M}, we have
$\mu^{n-1}_{\sima}(\p B_r(q)) < \infty$, for almost any $r>0$ and $q \in H$. From the lower semi-continuity of perimeters and the compactness of the BV-function space in the $L_\mathrm{loc}^1$-function space  \cite[Prop.3.6 and 3.7]{M}  we get, as in \cite[Th.1.20]{Gi}:
\begin{proposition}\emph{\textbf{(Plateau Problems in $(H,d_{\sima},\mu_{\sima})$)}} Let $\Omega \subset H$ be a bounded open and orientable
set and let $A \subset H$ be a Caccioppoli set. Then there exists a set $E\subset H$ coinciding with
$A$ outside $\Omega$ and such that
\begin{equation}\label{mininh}
\mu^{n-1}_{\sima}(\p E \cap \Omega) \le \mu^{n-1}_{\sima}(\p F \cap \Omega)
\end{equation}
for every Borel set $F  \subset H$ with $F = A$ outside $\Omega$.
\end{proposition}
The classical regularity theory of \cite[Ch.8]{Gi} applies in the manifold $H \setminus \Sigma$. It shows that such a minimizer $E$ can be assumed to be an open subset of $\Omega$ with boundary $\p E$ and so that $\p E \cap \Omega \setminus \Sigma$ is an area minimizing hypersurface smooth outside a set of Hausdorff dimension $\le n-7$. In this case we call  $\p E$ an \textbf{oriented minimal boundary}.\\

From the Poincar\'{e} inequality (\ref{poinm}) and the Ahlfors $n$-regularity we have  \cite[Remark 4.6]{M}:

\begin{corollary} [Isoperimetric Inequality]\label{iipc}   For $H \in {\cal{G}}$ there is a constant $\gamma(H,\Phi)>0,$ depending only on $n$ when $H \in {\cal{H}}^{\R}_n$, so that for any Caccioppoli set $U \subset H$:
\begin{equation}\label{ii23}
\min  \{ \mu_{\sima}(B_{\rho} \cap U),  \mu_{\sima} (B_{\rho} \setminus U)\}^{(n-1)/n} \le \gamma \cdot \mu^{n-1}_{\sima}(B_{\rho} \cap \p U), \mm{ for any } \rho>0,
\end{equation}
\end{corollary}

From this and again the Ahlfors regularity, we have a counterpart of Euclidean volume growth estimates in \cite[Prop.~5.14]{Gi}  for area minimizing hypersurfaces in $(H,d_{\sima},\mu_{\sima})$.

\begin{proposition}[Volume Growth of Area Minimizers] \label{vol}  For $(H,d_{\sima},\mu_{\sima})$, some open  subset $\Omega \subset H$ and an oriented minimal boundary $L^{n-1} \subset \Omega$ bounding an open set $L^+ \subset \Omega$ there are constants $\kappa, \kappa^+(H,\Phi)  >0$, so that for any $p \in L$:
\begin{equation}\label{est}
\kappa \cdot r^{n-1}  \le \mu^{n-1}_{\sima}(L \cap B_r(p)) \,\mm{ and }\, \kappa^+\cdot r^n  \le  \mu_{\sima}(L^+ \cap B_r(p)),
\end{equation}
for $r \in [0, (A/B)^{1/n} \cdot \dist(p,\p \Omega)/4)$, where $0<A<B$ are the Ahlfors constants. For $H \in {\cal{H}}^{\R}_n$, $\kappa, \kappa^+  >0$ depend only on $n$.
\end{proposition}

\textbf{Proof} \,  We start with the inequality for $L^+$. Since $L$ is area minimizing, we get
\small \begin{equation}\label{i1}
\mu^{n-1}_{\sima}(L \cap B_r(p)) \le \mu^{n-1}_{\sima}(L^+ \cap  \p B_r(p)).
\end{equation}
\normalsize
Since $r \mapsto \mu_{\sima}(L^+ \cap  B_r(p))$ is nondecreasing and bounded, it is differentiable almost everywhere on $\R^{>0}$. For almost any $r>0$ we therefore get the following two inequalities:
\small\begin{equation}\label{i2}
\mu^{n-1}_{\sima}(\p (L^+ \cap B_r(p))) =  \mu^{n-1}_{\sima}(L \cap B_r(p)) + \mu^{n-1}_{\sima}(L^+ \cap  \p B_r(p)).
\end{equation}
\normalsize
\vspace{-0.7cm}
\small\begin{equation}\label{i3}
\mu^{n-1}_{\sima}(\p (L^+ \cap B_r(p))) \le  2  \cdot  \mu^{n-1}_{\sima}(L^+ \cap \p B_r(p)) = 2  \cdot \frac{\p}{\p r}  \mu_{\sima}(L^+ \cap  B_r(p)).
\end{equation}
\normalsize
From the Ahlfors regularity, we notice that for $U= L^+ \cap B_r(p)$, $r \in [0, (A/B)^{1/n} \cdot \dist(p,\p \Omega)/4)$ and $\rho =\dist(p,\p \Omega)$ we have for $s \ge (B/A)^{1/n}$ and $q \in H$: $\mu_{\sima}(B_r(p) )\le B \cdot r^n  \le \mu_{\sima}(B_{s \cdot r}(q))$.
From this we have for $B_r(p)  \cap B_{s \cdot r}(q) \v,  B_r(p) \cup B_{s \cdot r}(q) \subset B_{\rho}(p):$
\begin{equation}\label{iu2}
 \mu_{\sima}(B_{\rho}(p) \cap U) \le \mu_{\sima} (B_{\rho}(p) \setminus U)
\end{equation}
and the isoperimetric inequality Cor.~\ref{iipc} for  $U=L^+ \cap B_r(p) \subset B_\rho(p)$ shows:
\small\begin{equation}\label{i4}
 \mu_{\sima}(L^+ \cap B_r(p))^{(n-1)/n} \le 2 \cdot \gamma \cdot \frac{\p}{\p r}  \mu_{\sima}(L^+ \cap B_r(p)).
\end{equation}
\normalsize
Integration gives the lower bound: $\mu_{\sima}(L^+ \cap B_r(p)) \ge \kappa^+ \cdot r^n$, for some $\kappa^+ >0$. The same estimate applies to $L^-=H \setminus \overline{L^+}$. From this we get the inequality for $\mu^{n-1}_{\sima}(L \cap B_r(p))$ using again the isoperimetric inequality: $\gamma \cdot \mu^{n-1}_{\sima}(L \cap B_{r}) \ge ({\kappa^+ \cdot r^n})^{(n-1)/n}$.
\qed

\footnotesize
\renewcommand{\refname}{\fontsize{14}{0}\selectfont

\end{document}